\documentclass[12pt]{amsart}
\usepackage[utf8]{inputenc}
\usepackage[T1]{fontenc}
\usepackage{lmodern}
\usepackage[english]{babel}
\usepackage{amsmath,mathtools,amsthm,amssymb,amsfonts}
\usepackage{mathrsfs}
\usepackage{graphicx}
\usepackage{enumitem}
\usepackage{tikz} 
\usepackage{pgfplots}
\usepackage{esvect}
\usepackage{pdfpages}
\usetikzlibrary{shapes,arrows}
\usetikzlibrary{calc}
\usetikzlibrary{shapes.geometric}
\usetikzlibrary{shapes.arrows}
\usetikzlibrary{arrows.meta}
\usetikzlibrary{decorations.markings}
\usetikzlibrary{fit}
\usetikzlibrary{patterns}
\usetikzlibrary{hobby}
\usepgfplotslibrary{patchplots}
\pgfplotsset{compat=1.11}
\usepackage{hyperref}
\usepackage{nicefrac}
\usepackage{cancel}
\usepackage{esint}
\usepackage{color}
\usepackage{etoolbox}
\theoremstyle{plain}
\newtheorem{theorem}{Theorem}[section]
\newtheorem{corollary}[theorem]{Corollary}
\newtheorem{conjecture}[theorem]{Conjecture}
\newtheorem{proposition}[theorem]{Proposition}
\newtheorem{lemma}[theorem]{Lemma}
\theoremstyle{definition}
\newtheorem{definition}[theorem]{Definition}
\newtheorem{remark}[theorem]{Remark}
\newtheorem*{theorem*}{Theorem}
\newtheorem*{remark*}{Remark}
\newtheorem*{thm*}{Theorem}
\newtheorem*{conjecture*}{Conjecture}

\newcommand{\definedas}{\mathrel{\raise.095ex\hbox{\rm :}\mkern-5.2mu=}}
\newcommand{\asdefined}{\mathrel{=\mkern-5.2mu\raise.095ex\hbox{\rm :}}}

\newcommand{\R}{\mathbb{R}}
\newcommand{\abs}[1]{\ensuremath{\left\vert #1\right\vert}}

\newcommand{\diver}{\operatorname{div}}
\newcommand{\tr}{\operatorname{tr}}
\newcommand{\scal}{\operatorname{R}}
\newcommand{\RN}{$\Lambda$-Reissner--Nordstr\"om manifold}
\newcommand{\DEC}{DEC}
\newcommand{\sDEC}{strict DEC}
\newcommand{\ECR}{electrically charged Riemannian manifold}
\AtBeginEnvironment{itemize}{\apptocmd{\item}{\phantomsection}{}{\errmessage{couldn't patch item}}}
\allowdisplaybreaks

\title[Electrically charged extensions with prescribed asymptotics]{Constructing electrically charged Riemannian manifolds with minimal boundary, prescribed asymptotics, and controlled mass}
\author[{Cabrera Pacheco}]{Armando J. {Cabrera Pacheco}}
\address{Department of Mathematics, Universit\"at T\"ubingen, 72076 T\"ubingen, Germany}
\email{cabrera@math.uni-tuebingen.de}

\author[Cederbaum]{Carla Cederbaum}
\address{Department of Mathematics, Universit\"at T\"ubingen, 72076 T\"ubingen, Germany}
\email{cederbaum@math.uni-tuebingen.de}

\author[Gehring]{Penelope Gehring}
\address{Max-Planck Institute for Gravitational Physics, 14476 Potsdam, Germany}
\email{penelope.gehring@aei.mpg.de}

\author[Pe\~nuela Diaz]{Alejandro Pe\~nuela Diaz}
\address{Max-Planck Institute for Gravitational Physics, 14476 Potsdam, Germany}
\email{alejandro.penuela@aei.mpg.de}

\begin{document}
\begin{abstract}
In \cite{pri}, Mantoulidis and Schoen constructed $3$-dimensional asymptotically Euclidean manifolds with non-negative scalar curvature whose ADM mass can be made arbitrarily close to the optimal value of the Riemannian Penrose Inequality, while the intrinsic geometry of the outermost minimal surface can be ``far away'' from being round. The resulting manifolds, called \emph{extensions}, are geometrically not ``close'' to a spatial Schwarzschild manifold. This suggests instability of the Riemannian Penrose Inequality as discussed in~\cite{CQG+,survey}. Their construction was later adapted to $n+1$ dimensions by Cabrera Pacheco and Miao~\cite{sec}. In recent papers by Alaee, Cabrera Pacheco, and Cederbaum~\cite{impor} and by Cabrera Pacheco, Cederbaum, and McCormick \cite{hyper}, a similar construction was performed for asymptotically Euclidean \ECR s and for asymptotically hyperbolic Riemannian manifolds, respectively, obtaining $3$-dimensional extensions.

This paper combines and generalizes all the aforementioned results by constructing suitable asymptotically hyperbolic or asymptotically Euclidean extensions with electric charge in $n+1$ dimensions for $n\geq2$. We study in detail the sub-extremality of these manifolds,  consider  the so far unstudied case of extremality in  extensions with electric charge and allow more general conditions for the metric of our extensions.

Besides suggesting instability of a naturally conjecturally generalized Riemannian Penrose Inequality, the constructed extensions give insights into an ad hoc generalized notion of Bartnik mass, similar to the Bartnik mass estimate for minimal surfaces proven by Mantoulidis and Schoen via their extensions, and unifying the Bartnik mass estimates in the various scenarios mentioned above.
\end{abstract}
	
\maketitle
	
\section{Introduction and Results}
The celebrated Riemannian Penrose Inequality states that every asymptotically Euclidean Riemannian $3$-manifold $(M^3,\gamma)$ with ADM mass $m$, non-negative scalar curvature $\scal(\gamma)\geq0$, and strictly outward minimizing minimal boundary $\partial M$ satisfies
\begin{align}\label{eq:RPIintro}
m \geq \sqrt{\frac{\abs{\partial M}}{16 \pi}}, 
\end{align}
where $\abs{\partial M}$ denotes the area of $\partial M$ with respect to the induced area measure.

Equality holds if and only if $(M^3,\gamma)$ is isometric to the spatial Schwarzschild manifold of mass $m$ up to its minimal inner boundary. The Riemannian Penrose Inequality, including the rigidity statement, was proven independently by Huisken and Ilmanen~\cite{RPI} assuming that the boundary $\Sigma$ is connected and by Bray~\cite{RPI2} in the general case. 

Recent years have seen an increasing interest in understanding the \emph{stability} or \emph{qualitative rigidity} of geometric inequalities such as the Riemannian Penrose Inequality~\eqref{eq:RPIintro}. This is to say, if a manifold $(M^{3},\gamma)$ satisfies the assumptions of the Riemannian Penrose Inequality and if it almost saturates~\eqref{eq:RPIintro}, does this imply that $(M^{3},\gamma)$ is necessarily geometrically ``close'' to the spatial Schwarzschild manifold of mass $m$? 

This very subtle topic was first addressed by Lee and Sormani~\cite{Lee-Sormani} in the rotationally symmetric case. They found sufficient conditions for stability in the bi-Lipschitz topology; at the same time, they identified a rotationally symmetric sequence of Riemannian manifolds satisfying the assumptions of the Riemannian Penrose Inequality, all with the same minimal boundary area $\abs{\partial M}$, with masses that saturate the inequality in the limit. Yet, this sequence converges to a spatial Schwarzschild manifold of mass $m=\sqrt{\frac{\vert\partial M\vert}{16 \pi}}$ with a round cylinder attached to the minimal inner boundary -- a Riemannian manifold manifestly non-isometric to the Schwarzschild manifold itself. This sequence indicates an instability of the Riemannian Penrose Inequality.

More recently, Allen \cite{Allen} proved stability of the Riemannian Penrose Inequality~\eqref{eq:RPIintro} in a Sobolev topology assuming the manifolds in a sequence as above are uniformly foliated by a smooth, embedded Inverse Mean Curvature Flow.

Moreover, Mantoulidis and Schoen's work~\cite{pri} -- initially tailored to computing the Bartnik mass of minimal surfaces, see below -- can also be interpreted as suggesting instability of the Riemannian Penrose Inequality~\cite{CQG+,survey}. Indeed, for a given Riemannian metric $g_o$ on the $2$-sphere $\mathbb{S}^{2}$ satisfying a certain eigenvalue condition or in particular having positive Gau{\ss} curvature, and for a given mass parameter $m$ satisfying $m>\sqrt{\frac{\vert\mathbb{S}^{2}\vert_{g_o}}{16\pi}}$, they construct an asymptotically Euclidean Riemannian $3$-manifold $(M^{3},\gamma_{m})$ with inner boundary $\partial M$ which has non-negative scalar curvature $\scal(\gamma_{m})\geq0$ and ADM mass~$m$, such that the boundary $\partial M$ is isometric to $(\mathbb{S}^{2},g_o)$, strictly outward minimizing, and minimal. Their explicit construction as well as the fact that the boundary $(\mathbb{S}^{2},g_o)$ need not be round suggest that these manifolds will generically not be geometrically ``close'' to the spatial Schwarzschild manifold of mass $m$, even if $m$ is close to the optimal value $\sqrt{\frac{\vert\mathbb{S}^{2}\vert_{g_o}}{16\pi}}$ in the Riemannian Penrose Inequality~\eqref{eq:RPIintro}. Many authors have since refined and/or generalized this construction. We direct the reader to \cite{survey} for a survey, and due to its relevance for potential future work, we highlight the work of Miao and Xie~\cite{Miao1} which provides a new perspective.

The purpose of this paper is to adapt and generalize the Mantoulidis--Schoen construction~\cite{pri} simultaneously to higher dimensions, to the presence of an electric field, and to both asymptotically Euclidean and asymptotically hyperbolic Riemannian manifolds, unifying and extending the higher dimensional generalizations by Cabrera~Pacheco and Miao~\cite{sec}, the $3$-dimensional generalization in the presence of an electric field by Alaee, Cabrera~Pacheco, and Cederbaum~\cite{impor}, and the $3$-dimensional asymptotically hyperbolic generalization by Cabrera~Pacheco, Cederbaum, and McCormick~\cite{hyper}. 

The results presented here can be viewed as indicating instability of the known or conjectured Riemannian Penrose Inequalities in higher dimensions, in the presence of an electric field, a non-positive cosmological constant, and for the two considered types of asymptotic assumptions described in Section~\ref{sec:prelims}.

We obtain the following main result:
\begin{theorem*}[Existence of Extensions]
Let $n\geq2$, let $g_o$ be a smooth Riemannian metric on $\mathbb{S}^{n}$ and let $r_o$ denote the volume radius of $g_o$, $\abs{\mathbb{S}^n}_{g_0}\asdefined\omega_n r_o^n$, where $\omega_{n}=\vert \mathbb{S}^{n}\vert_{g_{*}}$ denotes the volume of the unit round sphere $(\mathbb{S}^{n},g_{*})$. Assume that one of the following holds
\begin{enumerate}
\item $n\geq2$, $\lambda_{1}(-\triangle_{g_o}+\tfrac{1}{2}\scal(g_{o}))>0$, and $(\mathbb{S}^{n},g_{o})$ is conformal to $(\mathbb{S}^{n},g_{*})$,\\ where $\lambda_{1}(-\triangle_{g_o}+\tfrac{1}{2}\scal(g_{o}))$ denotes the first eigenvalue of $-\triangle_{g_o}+\tfrac{1}{2}\scal(g_o)$, or
\item $n=2$ and $K(g_o)>-\tau$ on $\mathbb{S}^{2}$ for some $\tau>0$ or
\item $2\leq n\leq3$  and $\scal(g_o)>0$ on $\mathbb{S}^{n}$, where $\scal(g_o)$ denotes the scalar curvature of $g_o$, or
\item $n>3$, $\scal(g_o)>0$ on $\mathbb{S}^{n}$, and $(\mathbb{S}^n,g_o)$ has pointwise $\nicefrac{1}{4}$-pinched sectional curvature.
\end{enumerate}
Here, $\scal(g_{o})$ and $K(g_{o})$ denote the scalar and Gaussian curvatures of $g_{o}$, respectively. Then there exists a constant $\kappa \geq 0$ with $\kappa<\lambda_{1}(-\triangle_{g_o}+\tfrac{1}{2}\scal(g_{o}))$ in Case~(1), $\kappa\leq\tau$ in Case~(2), and $2\kappa<\scal(g_o)$ in Cases~(3) and~(4) such that, for any constant $\Lambda\leq0$ and any constant $q\in\R$ satisfying
\begin{align*}
\begin{split}
\begin{array}{lll}
\frac{n(n-1)q^2}{r_{o}^{2n}}&<2(\kappa-\Lambda)&\text{ in Cases~(1), (3), and (4), and}\\
\frac{q^2}{r_{o}^{4}}&<-(\kappa+\Lambda)&\text{ in Case~(2)},
\end{array}
\end{split}
\end{align*}
and any constant $m\in\mathbb{R}$ satisfying
\begin{align*}
m>m_o&\definedas\frac{r_o^{n-1}}{2}\left(1+\frac{q^2}{r_o^{2(n-1)}} -\frac{2\Lambda r_o^{2}}{n(n+1)}\right),
\end{align*}
there exists an $(n+1)$-dimensional \ECR\ $(M,\gamma,E,\Lambda)$ which is geodesically complete up to its strictly outward minimizing minimal boundary $\partial M$, satisfies the Dominant Energy Condition
\begin{align*}
\scal(\gamma)\geq 2\Lambda+n(n-1)\abs{E}_\gamma^2,
\end{align*}
and has the following properties
\begin{enumerate}[label=(\roman*)]
\item $(\partial M,g)$ is isometric to $(\mathbb{S}^n,g_o)$, where $g$ is the metric induced on $\partial M$ by $\gamma$, and
\item there exists a compact subset $C\subset M$ and a radius $r_{C}>r_{+}$ such that $(M\setminus C,\gamma,E,\Lambda)$ is isometric to the subset $((r_{C},\infty)\times\mathbb{S}^n,\gamma_{m,q,\Lambda},E_{m,q,\Lambda},\Lambda)$ of the sub-extremal \RN\ of mass $m$, charge $q$, and cosmological constant $\Lambda$, with $r_{+}=r_{+}(m,q,\Lambda)$ denoting the volume radius of its horizon inner boundary.
\end{enumerate}
In particular, the (conjectured) Riemannian Penrose Inequality~\eqref{eq:RPI} in $(M,\gamma,E,\Lambda)$ reduces to $m_{o}\leq m$. In other words, the total mass $m$ can be chosen arbitrarily close to the optimal value in the (conjectured) Riemannian Penrose Inequality~\eqref{eq:RPI}.
\end{theorem*}
Note that, for $n> 2$, this theorem allows for more general conditions on the metric $g_o$ as were used previously in \cite{sec}. This is made possible by our use of different flows and paths of metrics. For more details see Remark \ref{metrirema}.

As we have briefly mentioned above, Mantoulidis and Schoen~\cite{pri} initially constructed Riemannian manifolds with prescribed minimal, strictly outward minimizing inner boundary --- called ``extensions'' --- in order to compute the Bartnik mass of minimal Riemannian $2$-spheres. These estimates have been generalized to ad hoc generalizations of Bartnik mass in the presence of an electric field by Alaee, Cabrera~Pacheco, and Cederbaum~\cite{impor} and to the asymptotically hyperbolic context by Cabrera~Pacheco, Cederbaum, and McCormick~\cite{hyper}, both in $3$ dimensions. Introducing an ad hoc generalization of Bartnik mass (see Section~\ref{sec:discussion}), we obtain the following unifying generalization of these results.

\begin{theorem*}[Generalized Bartnik Mass Estimate]
Under the same assumptions as in the previous theorem, there exist constants $\kappa\geq 0$, $\Lambda\leq0$, and $q\in\R$ satisfying the same properties mentioned in the previous theorem, such that the minimal (generalized) Bartnik data $(\mathbb{S}^{n},g_{o},H_{o},q,\Lambda)$ have $\mathcal{A}(\mathbb{S}^{n},g_{o},H_{o},q,\Lambda)\neq\emptyset$ and
\begin{align*}
\mathcal{B}(\mathbb{S}^{n},g_{o},H_{o}=0,q,\Lambda)&\leq\mathcal{M}(\mathbb{S}^{n},g_{o},H_{o}=0,q,\Lambda).
\end{align*}
Moreover, the \RN\ of mass $\mathcal{M}(\mathbb{S}^{n},g_{o},H_{o}=0,q,\Lambda)$, charge $q$, and cosmological constant $\Lambda$ is sub-extremal. Here, $\mathcal{A}$, $\mathcal{B}$, and $\mathcal{M}$ denote the class of admissible extensions, the generalized Bartnik mass functional (see Section~\ref{sec:discussion}) and the generalized Hawking mass functional (see Definition~\ref{def:Hawking}), respectively.
\end{theorem*}

This article is organized as follows: After introducing preliminaries such as the \RN s and their properties, general \ECR s, a notion of electric charge, and a generalized notion of Hawking mass in Section~\ref{sec:prelims}, we will construct collar extensions in Section~\ref{sec:collar}. This section also contains an analysis of when such collar extensions satisfy the \DEC\ as well as estimates on the generalized Hawking mass at the far end of the collars. Next, in Section~\ref{sec:gluing}, we prove technical lemmas that will later allow us to smoothly glue a given collar extension to a piece of a \RN. The following Section~\ref{sec:main} combines these technical tools to prove existence of extensions almost saturating the (conjectured) Riemannian Penrose Inequality. We will also provide more context for our result in view of the (in-)stability of the (conjectured) Riemannian Penrose Inequality. Finally, in Section~\ref{sec:discussion}, we will introduce an ad hoc notion of generalized Bartnik mass and reinterpret our main result as a Bartnik mass estimate for minimal hypersurfaces.

\paragraph*{\emph{Acknowledgements.}} We would like to thank Gerhard Huisken, Klaus Kr\"oncke, Stephen McCormick, Jan Metzger, Anna Sakovich, and Markus Wolff for interesting and helpful discussions. Aspects of the results presented here were contained in the MSc theses~\cite{G,P}.
	
This research is supported by the International Max Planck Research School for Mathematical and Physical Aspects of Gravitation, Cosmology and Quantum Field Theory and by the focus program on Geometry at Infinity (Deutsche Forschungsgemeinschaft, SPP 2026). The work of Armando Cabrera~Pacheco was supported by the Carl Zeiss Foundation. The work of Carla Cederbaum was supported by the Institutional Strategy of the University of T\"ubingen (Deutsche Forschungsgemeinschaft, ZUK 63).

\section{Preliminaries} \label{sec:prelims} 
Before we define the central notions needed for this work, let us dedicate some time to spatial \RN s, a central tool and model for the construction of extensions to be presented in this paper. We recall that the $(n+1)$-dimensional \emph{(spatial) Schwarzschild manifold $(M_{m}=(r_+, \infty) \times \mathbb{S}^n,\gamma_{m})$}, $n\geq2$, with Riemannian metric $\gamma_{m}$ given~by
\begin{align}
\gamma_{m}&=\left( 1-\frac{2m}{r^{n-1}}\right)^{-1} dr^2 +r^2g_*,
\end{align}
where $g_*$ denotes the standard round metric and $r_{+}\definedas(2m)^{\frac{1}{n-1}}$ if $m>0$, while $r_{+}\definedas0$ if $m\leq0$, is one of the most important examples of an asymptotically Euclidean manifold, central both in geometric analysis and general relativity. When its mass parameter $m\in\R$ is positive, it models the $\lbrace t=0\rbrace$-time-slice of a static, rotationally symmetric black hole in vacuum, isolated from external influences; for $m=0$, it coincides with Euclidean space in spherical coordinates, away from the origin; for $m<0$ it has a singularity at its center and is considered unphysical. 

The spatial Schwarzschild manifolds can be adapted to include an electric charge parameter $q\in\R$; the resulting manifolds are called \emph{(spatial) Reissner--Nordstr\"om manifolds}. In addition, one can combine this generalization with a (negative) cosmological constant $\Lambda<0$ as in the well-known \emph{(spatial) anti-de Sitter manifold} $(M_{\Lambda}=(0,\infty)\times\mathbb{S}^{n},\gamma_{\Lambda})$ with Riemannian metric $\gamma_{\Lambda}$ given by
\begin{align}
\gamma_{\Lambda}&=\left( 1-\frac{2\Lambda r^{2}}{n(n+1)}\right)^{-1} dr^2 +r^2g_*
\end{align}
 or in other words as in hyperbolic space with \emph{hyperbolic radius $\sqrt{-\tfrac{n(n+1)}{2\Lambda}}$. Many authors only consider the special case $\Lambda=\Lambda_{n}=-\tfrac{n(n+1)}{2}$ with hyperbolic radius $1$.} Of course, including a negative cosmological constant $\Lambda$ will change the asymptotics of the metric from being asymptotically Euclidean to asymptotically hyperbolic, see below. All these generalizations can be combined into one which we will call the \emph{(spatial) \RN s}, allowing both $\Lambda=0$ and $\Lambda<0$, which we will now define.

\subsection{Spatial $\Lambda$-Reissner--Nordstr\"om Manifolds}\label{sec:RN}
\begin{definition}[Spatial $\Lambda$-Reissner--Nordstr\"om Manifolds]\label{RN}
Let $m,q,\Lambda\in\mathbb{R}$ be constants with $\Lambda\leq 0$ and let $n\geq2$. Using the abbreviations
\begin{align}\label{eq:defp}
p_{m,q,\Lambda}(r)&\definedas1-\frac{2m}{r^{n-1} } +\frac{q^2}{r^{2(n-1)}}-\frac{2\Lambda r^2}{n(n+1)},
\end{align} 
for $r>0$ and $r_+\definedas\max\lbrace{\text{largest positive root of }p_{m,q,\Lambda},0\rbrace}$, the $(n+1)$-dimensional \emph{(spatial) \RN\ $(M_{m,q,\Lambda}= (r_+, \infty) \times \mathbb{S}^n,\gamma_{m,q,\Lambda})$} is given by
\begin{align}\label{metricg}
\gamma_{m,q,\Lambda}&= \frac{dr^2}{ p_{m,q,\Lambda}(r)} +r^2g_*,
\end{align}
where $g_*$ denotes the standard round metric on $\mathbb{S}^{n}$. Its \emph{electric field} is the smooth vector field $E_{m,q,\Lambda}\in\Gamma(TM_{m,q,\Lambda})$ given by 	
\begin{align}
E_{m,q,\Lambda}&=\frac{q\left( p_{m,q,\Lambda}(r)\right)^\frac{1}{2}}{r^n}  \partial_r .
\end{align}
The parameters $m,q,\Lambda$ are called the \emph{mass}, the \emph{(electric) charge}, and the \emph{cosmological constant} of the \RN, respectively.\footnote{Note that we are suppressing the dependence on the dimensional parameter $n$ in the names of the auxiliary function $p_{m,q,\Lambda}$, the metric, the electric field, and the manifold.}
\end{definition}

To see why $r_{+}\geq0$ is well-defined note that the auxiliary function $p_{m,q,\Lambda}$ can be rewritten as a polynomial in $r$ divided by $r^{2(n-1)}$. Hence, $p_{m,q,\Lambda}$ has finitely many positive roots (if any); if it has no positive roots, then $r_{+}=0$ is to be understood.

The following facts are well-known properties of \RN s and can be verified by straightforward computations.
\begin{proposition}\label{prop:RN}
Let $m,q,\Lambda\in\mathbb{R}$ be constants with $\Lambda\leq 0$ and let $n\geq2$. Then $E_{m,q,\Lambda}$ is divergence-free,
\begin{align}
\diver_{\gamma_{m,q,\Lambda}}E_{m,q,\Lambda}&=0
\end{align}
on $M_{m,q,\Lambda}$, and the scalar curvature of $(M_{m,q,\Lambda},\gamma_{m,q,\Lambda})$ is given by 
\begin{align}\label{eq:exactDEC}
\scal(\gamma_{m,q,\Lambda})= 2\Lambda+n(n-1) |E_{m,q,\Lambda}|^2_{\gamma_{m,q,\Lambda}},
\end{align}
where $\vert\cdot\vert_{\gamma_{m,q,\Lambda}}$ denotes the length computed with respect to $\gamma_{m,q,\Lambda}$ and $\scal(\gamma_{m,q,\Lambda})$ denotes the scalar curvature of $\gamma_{m,q,\Lambda}$.
\end{proposition}

\begin{remark}\label{rem:special}
As described above, the definition of \RN s specializes to many more well-known examples of manifolds for suitable choices of parameters:
\begin{enumerate}
\item For $q=\Lambda=0$, one recovers the (higher dimensional, spatial) Schwarzschild manifold of mass $m$.
\item For $\Lambda=0$, one recovers the (higher dimensional, spatial) Reissner--Nordstr\"om manifold of mass $m$ and electric charge $q$.
\item For $m=q=0$, $\Lambda<0$, one recovers the $(n+1)$-dimensional hyperbolic space of hyperbolic radius $\sqrt{-\tfrac{n(n+1)}{2\Lambda}}$ in spherical coordinates, away from the origin.
\item For $q=0$ and $\Lambda=\Lambda_{n}=-\frac{n(n+1)}{2}$, one recovers the \emph{(spatial) anti-de Sitter--Schwarz\-schild manifold} of mass $m$.
\end{enumerate}
\end{remark}

In the Schwarzschild case, it is well-known that only choices of mass parameter $m\geq0$ correspond to physically reasonable solutions (of the vacuum Einstein equations), modeling the vacuum exterior region of a static, rotationally symmetric black hole when $m>0$. A similar, but more involved distinction exists for \RN s; here, one distinguishes between the physically reasonable \emph{sub-extremal case},  modeling the electro-vacuum exterior region of a ``non-degenerate'' static, rotationally symmetric black hole, the \emph{extremal case}, where the black hole degenerates in a well-defined
sense, and the \emph{super-extremal case} not modeling a black hole of any sorts. The last case is typically considered to be unphysical (except for Euclidean and hyperbolic space, see Remark~\ref{rem:specialsub}). This is made precise in the following definition.

\begin{definition}[Extremality, Sub-Extremality, and Super-Extremality]\label{def:extremality}
Let $m,q,\Lambda\in\mathbb{R}$ be constants with $\Lambda\leq 0$ and let $n\geq2$. The $(n+1)$-dimensional spatial \RN\ of mass $m$, charge $q$, and cosmological constant $\Lambda$ is called
\begin{itemize}
\item \emph{sub-extremal} if $r_{+}>0$ and $p_{m,q,\Lambda}'(r_+)>0$,
\item \emph{extremal} if $r_{+}>0$ and $p_{m,q,\Lambda}'(r_+)=0$, and
\item \emph{super-extremal} if $r_{+}=0$,
\end{itemize}
where we recall that $r_+=\max\lbrace{\text{largest root of }p_{m,q,\Lambda},0\rbrace}$. Here and in what follows, $'$ is a short-hand for $\frac{d}{dr}$.
\end{definition}
	
To see why this definition covers all possible cases, recall that the auxiliary function $p_{m,q,\Lambda}$ can be rewritten as a polynomial in $r$ with positive leading order coefficient, divided by $r^{2(n-1)}$. Hence, $p_{m,q,\Lambda}>0$ on $(r_{+},\infty)$ and if $r_{+}>0$, we must have $p_{m,q,\Lambda}'(r_{+})\geq0$.

\begin{remark}\label{rem:specialsub}
In the specializations discussed in Remark~\ref{rem:special}, we find these conditions to signify the following:
\begin{enumerate}
\item In the Schwarzschild case $q=\Lambda=0$, $m>0$ corresponds to being sub-extremal, while $m\leq0$ implies being super-extremal. No extremal cases occur.
\item In the Reissner--Nordstr\"om case $\Lambda=0$, the only possible roots of $p_{m,q,0}$ are $r_\pm=\left(m\pm\sqrt{m^2-q^2}\right)^{1/(n-1)}$ with $r_{+}\geq r_{-}$ being real and positive if and only if $q^2\leq m^{2}$ and $m>0$, otherwise the Reissner--Nordstr\"om manifold is super-extremal. For $q^2\leq m^{2}$ and $m>0$, $p_{m,q,0}'(r_+)$ vanishes if and only if $r_+=\left(\frac{q^2}{m}\right)^\frac{1}{n-1}$, which is equivalent to $m=\abs{q}$ and to $r_-=r_+$. Hence, the Reissner--Nordstr\"om manifold is sub-extremal if and only if $m>\abs{q}$ and extremal if and only if $m=\vert q\vert>0$.
\item The hyperbolic space case $m=q=0$, $\Lambda<0$, is always super-extremal.
\item In the anti-de Sitter--Schwarzschild case $q=0$, $\Lambda=\Lambda_{n}=-\frac{n(n+1)}{2}$, $p_{m,0,\Lambda_{n}}$ has a unique positive root $r_+$ and $p'_{m,0,\Lambda_{n}}>0$ is positive throughout if $m>0$. Hence these manifolds are always sub-extremal if $m>0$. If $m\leq0$, they are super-extremal.
\end{enumerate}
\end{remark}

Before we discuss more observations on sub-extremality and extremality for \RN s in more detail, let us introduce another auxiliary  function which will be very helpful in this paper. Let $q,\Lambda\in\R$ be constants with $\Lambda\leq0$ and let $n\geq2$. We set
\begin{align}\label{eq:hqLamba}
h_{q,\Lambda}(r)\definedas r^{2(n-1)}-q^{2}-\frac{2\Lambda r^{2n}}{n(n-1)}
\end{align}
for $r>0$.

\begin{lemma}[Properties of Roots of $p_{m,q,\Lambda}$]\label{lem:h}
Let $m,q,\Lambda\in\R$ be constants with $\Lambda\leq0$ and let $n\geq2$. Then $h_{q,\Lambda}\colon\R^{+}\to\R$ is strictly increasing with image $h_{q,\Lambda}(\R^{+})=(-q^{2},\infty)$ and, for any root $r_{o}>0$ of $p_{m,q,\Lambda}$, we have
\begin{align}\label{eq:hp'}
p_{m,q,\Lambda}'(r_{o})&=\frac{(n-1)h_{q,\Lambda}(r_{o})}{r_{o}^{2n-1}}.
\end{align}
Moreover, if $h_{q,\Lambda}(r_{o})=0$ then $p_{m,q,\Lambda}''(r_{o})>0$.
\end{lemma}
\begin{proof}
Strict monotonicity of $h_{q,\Lambda}$ and $h_{q,\Lambda}(\R^{+})=(-q^{2},\infty)$ are obvious. To see~\eqref{eq:hp'}, we compute that $p_{m,q,\Lambda}(r_{o})=0$ is equivalent to $\frac{2m}{r_{o}^{n-1}}=1+\frac{q^{2}}{r_{o}^{2(n-1)}}-\frac{2\Lambda r_{o}^{2}}{n(n+1)}$ from which we get
\begin{align*}
p_{m,q,\Lambda}'(r_{o})&=\frac{2(n-1)m}{r_{o}^{n}}-\frac{2(n-1)q^{2}}{r_{o}^{2n-1}}-\frac{4\Lambda r_{o}}{n(n+1)}\\
&=\frac{n-1}{r_{o}}-\frac{(n-1)q^{2}}{r_{o}^{2n-1}}-\frac{2\Lambda r_{o}}{n}=\frac{(n-1)h_{q,\Lambda}(r_{o})}{r_{o}^{2n-1}}.
\end{align*}
Moreover, using first the above characterization of $p_{m,q,\Lambda}(r_{o})=0$ and later the fact that $h_{q,\Lambda}(r_{o})=0$ is equivalent to $q^{2}=r_{o}^{2(n-1)}-\frac{2\Lambda r_{o}^{2n}}{n(n-1)}$, we find
\begin{align*}
p_{m,q,\Lambda}''(r_{o})&=-\frac{2n(n-1)m}{r_{o}^{n+1}}+\frac{2(n-1)(2n-1)q^{2}}{r_{o}^{2n}}-\frac{4\Lambda}{n(n+1)}\\
&=-\frac{n(n-1)}{r_{o}^{2}}+\frac{(n-1)(3n-2)q^{2}}{r_{o}^{2n}}+\frac{2(n-2)\Lambda}{n}=\frac{2(n-1)^{2}}{r_{o}^{2}}-4\Lambda
\end{align*}
if $h_{q,\Lambda}(r_{o})=0$, which is manifestly positive.
\end{proof}

The following observations on sub-extremality and extremality will be helpful later.
\begin{lemma}[Characterization of Sub-Extremality and Extremality]\label{lem:sub}
Let $m,q,\Lambda\in\mathbb{R}$ be constants with $\Lambda\leq 0$ and let $n\geq2$. If the $(n+1)$-dimensional \RN\ of mass $m$, charge $q$, and cosmological constant $\Lambda$ is sub-extremal then $m>\vert q\vert$ and $p_{m,q,\Lambda}$ has at most two positive roots $0<r_{-}<r_{+}$. The second root $r_{-}$ exists unless $q=0$ and $m>0$. In addition, if a second root $r_{-}$ exists one has $p_{m,q,\Lambda}'(r_{-})<0$. If the $(n+1)$-dimensional \RN\ of mass $m$, charge $q$, and cosmological constant $\Lambda$ is extremal then $q\neq0$ and $p_{m,q,\Lambda}$ has precisely one positive root $r_{+}$. 

Moreover, the $(n+1)$-dimensional \RN\ of mass $m$, charge $q$, and cosmological constant $\Lambda$ is sub-extremal if and only if $h_{q,\Lambda}(r_{+})>0$ and extremal if and only if $h_{q,\Lambda}(r_{+})=0$.
\end{lemma}

\begin{proof}
We can rewrite
\begin{align*}
p_{m,q,\Lambda}(r)&=\frac{2(\vert q\vert-m)}{r^{n-1}}+\underbrace{\left(1-\frac{\vert q\vert}{r^{n-1}}\right)^{2}-\frac{2\Lambda r^{2}}{n(n+1)}}_{\geq0\text{ for }r>0},
\end{align*}
hence existence of a (largest) positive root $r_{+}$ necessitates $m\geq \vert q \vert$, with $m=\vert q\vert$ implying $r_{+}^{n-1}=\vert q\vert$ and $\Lambda=0$. If $r_{+}^{n-1}=\vert q\vert$, and $\Lambda=0$, we find $h_{q,0}(r_{+})=h_{q,0}(\vert q\vert^{\frac{1}{n-1}})=0$ so that \eqref{eq:hp'} tells us that $m>\vert q\vert$ in the sub-extremal case. From~\eqref{eq:hp'}, we also get that $p_{m,q,\Lambda}'(r_{+})>0$ is equivalent to $h_{q,\Lambda}(r_{+})>0$ while $p_{m,q,\Lambda}'(r_{+})=0$ is equivalent to $h_{q,\Lambda}(r_{+})=0$. As $h_{q,\Lambda}$ is manifestly positive for $q=0$, we find that extremality implies $q\neq0$.

Now, in the extremal case, Lemma~\ref{lem:h} tells us that $p_{m,q,\Lambda}''(r_{o})>0$ for any root $r_{o}>0$ of $p_{m,q,\Lambda}$ with $p_{m,q,\Lambda}'(r_{o})=0$ which applies in particular to the largest root $r_{+}$. Hence $p_{m,q,\Lambda}>0$ holds on some open neighborhood $(r_{+}-\varepsilon,r_{+}+\varepsilon)\setminus\lbrace r_{+}\rbrace$, $\varepsilon>0$. By strict monotonicity of $h_{q,\Lambda}$ and~\eqref{eq:hp'}, it follows that any root $0<r_{o}<r_{+}$ would have $p_{m,q,\Lambda}'(r_{o})<0$ so that $p_{m,q,\Lambda}<0$ on some interval $(r_{o},r_{o}+\delta)$, $\delta>0$. By continuity there would hence need to be another root $r_{o}<r_{1}<r_{+}$ of $p_{m,q,\Lambda}$ near which $p_{m,q,\Lambda}$ is increasing from negative to positive values. By~\eqref{eq:hp'}, this implies that $h_{q,\Lambda}(r_{1})\geq0$, a contradiction to $r_{1}<r_{+}$, $h_{q,\Lambda}(r_{+})=0$, and strict monotonicity of $h_{q,\Lambda}$. This shows that $r_{+}$ is indeed the only positive root of $p_{m,q,\Lambda}$.

On the other hand, in the sub-extremal case, we have $h_{q,\Lambda}(r_{+})>0$ and $p_{m,q,\Lambda}<0$ on $(r_{+}-\varepsilon,r_{+})$ for some $\varepsilon>0$. Let $0<r_{-}<r_{+}$ be a root of $p_{m,q,\Lambda}$. Now suppose first that $h_{q,\Lambda}(r_{-})>0$. Arguing by Lemma~\ref{lem:h}, this gives that $p_{m,q,\Lambda}>0$ on $(r_{-},r_{-}+\delta)$ for some $\delta>0$. Hence by continuity, $p_{m,q,\Lambda}$ must have another root $r_{-}<r_{1}<r_{+}$ of $p_{m,q,\Lambda}$ near which $p_{m,q,\Lambda}$ is decreasing from positive to negative values so that $h_{q,\Lambda}(r_{1})\leq0$. This contradicts the strict monotonicity of $h_{q,\Lambda}$. Next, suppose that $h_{q,\Lambda}(r_{-})=0$. Arguing as in the extremal case, this would again imply the existence of another root $r_{-}<r_{1}<r_{+}$ of $p_{m,q,\Lambda}$ near which $p_{m,q,\Lambda}$ is decreasing from positive to negative values, a contradiction. Hence we find $h_{q,\Lambda}(r_{-})<0$ for any root of $p_{m,q,\Lambda}$ which satisfies $0<r_{-}<r_{+}$. Repeating the above arguments shows that there cannot be other roots $0<r_{o}\neq r_{-}, r_{+}$ of $p_{m,q,\Lambda}$. Applying~\eqref{eq:hp'} one last time, we conclude that $p_{m,q,\Lambda}'(r_{-})<0$ as claimed. By continuity, we can conclude furthermore that $r_{-}$ exists if and only if $p_{m,q,\Lambda}(r)>0$ near $r=0$, hence if and only if $q\neq0$ or $m\leq0$.
\end{proof}

\begin{lemma}[Inherited Sub-Extremality]\label{lem:largerm}
Let $n\geq 2$ and let $m$, $q$, $\Lambda\in\mathbb{R}$ be constants with $\Lambda\leq 0$. If the $(n+1)$-dimensional \RN\ $(M_{m,q,\Lambda},\gamma_{m,q,\Lambda})$ is sub-extremal then the $(n+1)$-dimensional \RN\ $(M_{\widetilde{m},\widetilde{q},\Lambda},\gamma_{\widetilde{m},\widetilde{q},\Lambda})$ is sub-extremal for any $\widetilde{m}\geq m$, $\widetilde{q}^{\,2} \leq q^{2}$, as well. Moreover, if $(M_{m,q,\Lambda},\gamma_{m,q,\Lambda})$ is extremal, then $(M_{\widetilde{m},\widetilde{q},\Lambda},\gamma_{\widetilde{m},\widetilde{q},\Lambda})$ is sub-extremal whenever $\widetilde{m}>m$, $\widetilde{q}^{\,2}\leq q^{2}$ or $\widetilde{m}\geq m$, $\widetilde{q}^{\,2}<q^{2}$.

Moreover, if $q\neq0$, there is a unique $m_{q,\Lambda}\in\R$ such that the $(n+1)$-dimensional \RN\ of mass $m$, charge $q$, and cosmological constant $\Lambda$ is sub-extremal for all $m>m_{q,\Lambda}$, extremal for $m=m_{q,\Lambda}$, and super-extremal for all $m<m_{q,\Lambda}$. If $q=0$, we set $m_{q,\Lambda}\definedas0$; then the $(n+1)$-dimensional \RN\ of mass $m$, charge $q$, and cosmological constant $\Lambda$ is sub-extremal for all $m>m_{q,\Lambda}=0$ and super-extremal for all $m\leq m_{q,\Lambda}=0$.
\end{lemma}

\begin{proof}
From $\widetilde{m}\geq m$, $\widetilde{q}^{\,2} \leq q^{2}$, we deduce that $p_{\widetilde{m},\widetilde{q},\Lambda}\leq p_{m,q,\Lambda}$ on $(0,\infty)$ so that the largest zero $\widetilde{r}_{+}$ of $p_{\widetilde{m},\widetilde{q}s,\Lambda}$ satisfies $\widetilde{r}_{+}\geq r_{+}>0$, with $r_{+}$ denoting the largest zero of $p_{m,q,\Lambda}$. This shows that $(M_{\widetilde{m},\widetilde{q},\Lambda},\gamma_{\widetilde{m},\widetilde{q},\Lambda})$ is not super-extremal. Now suppose it were extremal, i.e., suppose that $p_{\widetilde{m},\widetilde{q},\Lambda}'(\widetilde{r}_{+})=0$. By a direct computation, this gives $\widetilde{q}^{2}=\widetilde{m}\widetilde{r}_{+}^{\,n-1}-\tfrac{2\Lambda\widetilde{r}_{+}^{\,2n}}{n(n-1)(n+1)}$, while sub-extremality or extremality of $(M_{m,q,\Lambda},\gamma_{m,q,\Lambda})$, i.e., $p_{m,q,\Lambda}'(r_{+})\geq0$ is equivalent to $q^{2}\leq m r_{+}^{n-1}-\tfrac{2\Lambda r_{+}^{2n}}{n(n-1)(n+1)}$, with strict inequality in the sub-extremal case. Using $\widetilde{m}\geq m$, $\widetilde{q}^{\,2} \leq q^{2}$, $\widetilde{r}_{+}\geq r_{+}>0$, and $\Lambda\leq0$, this gives $0\geq mr_{+}^{n-1}-\widetilde{m}\widetilde{r_{+}}^{\,n-1}\geq-\tfrac{2\Lambda}{n(n-1)(n+1)}(\widetilde{r}_{+}^{\,2n}-r_{+}^{2n})\geq0$ which is a contradiction unless $\widetilde{m}=m$, $\widetilde{q}^{\,2} = q^{2}$, and $(M_{m,q,\Lambda},\gamma_{m,q,\Lambda})$ is extremal. Hence $(M_{\widetilde{m},\widetilde{q},\Lambda},\gamma_{\widetilde{m},\widetilde{q},\Lambda})$ is sub-extremal, as claimed.

Moreover, if $q\neq0$, Lemma~\ref{lem:h} and the Intermediate Value Theorem tell us that $h_{q,\Lambda}$ has a unique zero $r_{q,\Lambda}\in(0,\infty)$. We set $m_{q,\Lambda}\definedas \tfrac{r_{q,\Lambda}^{n-1}}{2}\left(1+\tfrac{q^{2}}{r_{q,\Lambda}^{2(n-1)}}-\tfrac{2\Lambda r_{q,\Lambda}^{2}}{n(n+1)}\right)$, so that $p_{m_{q,\Lambda},q,\Lambda}(r_{q,\Lambda})=0$ and $p_{m_{q,\Lambda},q,\Lambda}'(r_{q,\Lambda})=0$. From this, the claims about sub-extremality, extremality, and super-extremality follow from what we already showed in this proof. If, on the other hand, $q=0$, then $p_{m,0,\Lambda}$ is positive for $m\leq0$, and has $p_{m,0,\Lambda}(r)\to-\infty$ for $r\searrow0$ and $p_{m,0,\Lambda}>0$ for $m>0$ which proves the last claim, again by the Intermediate Value Theorem.
\end{proof}

Sub-extremal \RN s can be smoothly extended to their ``horizons'' $\lbrace r_{+}\rbrace\times\mathbb{S}^{n}$ as the following proposition asserts.

\begin{proposition}[Extension to Boundary, Radial Coordinate $s$, Radial Profile]\label{prop:extendRN}
Let $n\geq2$, and let $m,q,\Lambda\in\mathbb{R}$ be constants with $\Lambda\leq 0$ chosen such that the $(n+1)$-dimensional \RN\ $(M_{m,q,\Lambda},\gamma_{m,q,\Lambda})$ is sub-extremal. Changing the radial variable $r\in(r_{+},\infty)$ to 
\begin{align}\label{eq:radialtransfo}
s(r)&\definedas \int_{r_+}^r \left( p_{m,q,\Lambda}(t)\right)^{-\frac{1}{2}} dt,
\end{align}
$\gamma_{m,q,\Lambda}$ smoothly extends to $\overline{M_{m,q,\Lambda}}\approx [0,\infty)\times\mathbb{S}^{n}$ where it can be expressed as 
\begin{align}\label{eq:newgamma}
\gamma_{m,q,\Lambda}&=ds^2 +u_{m,q,\Lambda}(s)^2 g_*, 
\end{align}
with \emph{radial profile} $u_{m,q,\Lambda}\colon [0, \infty) \to [r_+, \infty)$ satisfying
\begin{enumerate}[label=(\roman*)]
\item $u_{m,q,\Lambda}(0)=r_+$
\item $u'_{m,q,\Lambda}=  \left(p_{m,q,\Lambda}\circ u_{m,q,\Lambda}\right)^\frac{1}{2}$, in particular $u'_{m,q,\Lambda}>0$ on $(0,\infty)$, $u'_{m,q,\Lambda}(0)=0$, and
\item $u''_{m,q,\Lambda}= \frac{1}{2} (p_{m,q,\Lambda}'\circ u_{m,q,\Lambda})$ and in particular $u''_{m,q,\Lambda}>0$ on $[0,\infty)$,
\end{enumerate}
where now $'$ abbreviates $\frac{d}{ds}$. Moreover, one finds that the electric field $E_{m,q,\Lambda}$ also smoothly extends to $\overline{M_{m,q,\Lambda}}\approx [0,\infty)\times\mathbb{S}^{n}$ and reads
\begin{align}\label{eq:newE}
E_{m,q,\Lambda}&=\frac{q}{u_{m,q,\Lambda}(s)^n} \partial_s
\end{align}
with respect to the new radial coordinate $s$.
\end{proposition}

\begin{proof}
By the facts that $p_{m,q,\Lambda}>0$ on $(r_{+},\infty)\times\mathbb{S}^{n}$ and $p_{m,q,\Lambda}'(r_{+})>0$ , the radial coordinate change is well-defined by Taylor's Formula. We clearly have $\lim_{r\searrow r_{+}}s(r)=0$ while $\limsup_{r\nearrow\infty}s(r)=\infty$ follows from the fact that  $p_{m,q,\Lambda}$ can be viewed as a polynomial in $r$ with positive leading order coefficient, divided by $r^{2(n-1)}$. Transforming the radial coordinate in $\gamma_{m,q,\Lambda}$ shows~\eqref{eq:newgamma} for $s>0$, with $u_{m,q,\Lambda}$ given by the inverse of the function $s\colon(r_{+},\infty)\to(0,\infty)$ defined in~\eqref{eq:radialtransfo} which can naturally be smoothly extended to $s(r_{+})\definedas0$. This shows~\emph{(i)}. By chain rule, one hence finds $u_{m,q,\Lambda}'=(s'\circ u_{m,q,\Lambda})^{-1}$ on $[0,\infty)$ which shows~\emph{(ii)}. Using the explicit form of the auxiliary function $p_{m,q,\Lambda}$, one finds~\emph{(iii)} from~\emph{(ii)}, while $u''_{m,q,\Lambda}(0)>0$ follows from the same computation and the sub-extremality condition $p'_{m,q,\Lambda}(r_+)>0$. The transformation and extension of the electric field is straightforward.
\end{proof}

\begin{remark}\label{rem:notextend}
The procedure described in Proposition~\ref{prop:extendRN} does not work for extremal \RN s because $(p_{m,q,\Lambda})^{-\frac{1}{2}}$ is not integrable near $r_{+}$ in this case and hence the coordinate change in~\eqref{eq:radialtransfo} would fail. In the super-extremal case, $(p_{m,q,\Lambda})^{-\frac{1}{2}}\leq 2$ near $r=0$ so integrability is not an issue. In both cases, we set
\begin{align}\label{eq:radialtransfomu}
s_{\mu}(r)&\definedas \int_{\mu}^r \left( p_{m,q,\Lambda}(t)\right)^{-\frac{1}{2}} dt,
\end{align}
for any $\mu>r_{+}$ in the extremal and any $\mu>0$ in the super-extremal case. Doing so gives~\eqref{eq:newgamma}, \eqref{eq:newE}, \emph{(ii)}, and \emph{(iii)} for the corresponding radial profile $u_{m,q,\Lambda}^{\mu}$ on $[0,\infty)\times\mathbb{S}^{n}\approx[\mu,\infty)\times\mathbb{S}^{n}\subset M_{m,q,\Lambda}$, and of course $u_{m,q,\Lambda}^{\mu}(0)=\mu$, ${u_{m,q,\Lambda}^{\mu}}'(0)=(p_{m,q,\Lambda}(\mu))^{\frac{1}{2}}>0$.
\end{remark}

\begin{proposition}[(Outermost Black Hole) Horizon]\label{prop:RNhorizon}
Let $n\geq2$ and let $m,q,\Lambda\in\mathbb{R}$ be constants with $\Lambda\leq 0$  chosen such that the $(n+1)$-dimensional \RN\ of mass $m$, charge $q$, and cosmological constant $\Lambda$ is sub-extremal. Then the coordinate sphere $\lbrace s=0\rbrace\times\mathbb{S}^{n}\approx \partial M_{m,q,\Lambda}=\lbrace r_{+}\rbrace\times\mathbb{S}^{n}$ is a \emph{minimal hypersurface}, i.e., its mean curvature $H$ satisfies $H=0$, and it is \emph{outward minimizing} in the sense\footnote{See Definition~\ref{def:horizons} for a more precise definition of outward minimizing minimal hypersurfaces and Theorem~\ref{thm:outwardminfoliation} for the assertion that a positive constant mean curvature foliation implies being outward minimizing.} that all other coordinate spheres $\lbrace {r\rbrace}\times\mathbb{S}^{n}$ for $r>r_{+}$ have positive constant mean curvature. It is called the \emph{(outermost black hole) horizon}.
\end{proposition}
\begin{proof}
Slightly abusing notation to include $r=r_{+}$, by rotational symmetry, all coordinate spheres $\lbrace {r\rbrace}\times\mathbb{S}^{n}$ for $r\geq r_{+}$ have constant mean curvature. What remains to be shown is that this mean curvature is positive for $r>r_{+}$ (or $s>0$) and vanishes for $r=r_{+}$ (or $s=0$).
Both claims follow from the more general formula~\eqref{eq:meancurvature} to be proved below, using the form of the metric~\eqref{eq:newgamma} with respect to the $s$-coordinate, see also Remark~\ref{rem:notextend}.
\end{proof}

\subsection{Time-Symmetric, Electrically Charged Initial Data Sets with Prescribed Asymptotics}
The \RN s introduced in Section~\ref{sec:RN} are specific examples of the types of manifolds appearing as extensions in the construction to be performed in this paper. Before we give precise definitions of those, let us give a little bit of context from general relativity.

In general relativity, $(n+1)$-dimensional Riemannian manifolds $(M,\gamma)$ arise as \emph{time-symmetric initial data sets} of an $(n+2)$-dimensional spacetime (or time-orientable Lorentzian manifold). Here, an \emph{initial data set} is a spacelike hypersurface $M$ in a spacetime, together with its induced Riemannian metric $\gamma$ and second fundamental form $K$. An initial data set is called \emph{time-symmetric} if $K=0$, i.e., if it is totally geodesic in the ambient spacetime. If the spacetime carries an electromagnetic field (expressed as a $2$-form $F$ on the spacetime), the initial data set also contains the electric and magnetic vector fields $E$ and $B$ induced on $M$ by the $2$-form $F$. For time-symmetric initial data sets $(M,\gamma)$, it is usually\footnote{especially for $n=2$, but also for $n>2$.} assumed that $B=0$ and we will follow this physically reasonable custom here.

An electromagnetic field is said to be source-free if its spacetime divergence vanishes; this implies in particular that $\diver_{\gamma}E=0$ on $M$ which can be reworded as the electric field $E$ having ``vanishing charge density''\footnote{It would also be physically reasonable to instead consider charge densities that have a sign on $M$. As the manifolds we construct here all naturally satisfy $\diver_{\gamma}E=0$, we will not pursue this option here.}. If an $(n+2)$-dimensional spacetime carrying an electromagnetic field and subject to a cosmological constant $\Lambda$ satisfies the so-called Dominant Energy Condition, a physical condition on the stress-energy tensor of the spacetime which can be reworded as a condition on the Ricci tensor of the spacetime, then any time-symmetric initial data set $(M,\gamma)$ with electric field $E$ in said spacetime is subject to the condition
\begin{align*}
\operatorname{R}(\gamma)\geq 2\Lambda+n(n-1)\abs{E}_\gamma^2,
\end{align*}
where $\operatorname{R}(\gamma)$ denotes the scalar curvature of $\gamma$ in $M$. We hence make the following definition.

\begin{definition}[Electrically Charged Riemannian Manifolds, Isometry, (Strict) Dominant Energy Condition]\label{def:ECRM}
Let $n\geq2$. We say that a tuple $(M,\gamma, E,\Lambda)$ is an \emph{\ECR\ (of dimension $n+1$)} if $(M,\gamma)$ is a smooth $(n+1)$-dimensional Riemannian manifold, $E$ is a smooth, divergence-free vector field on $M$, $\diver_{\gamma}E=0$, to be interpreted as the \emph{electric field}, and $\Lambda\leq0$ is a constant, interpreted as the \emph{cosmological constant}. Two \ECR s $(M_{i},\gamma_{i}, E_{i},\Lambda_{i})$, $i=1,2$, are said to be \emph{isometric} if $\Lambda_{1}=\Lambda_{2}$ and if there exists an isometry $\Phi\colon(M_{1},\gamma_{1})\to(M_{2},\gamma_{2})$ such that $\Phi_{*}E_{1}=E_{2}$. 

Furthermore, an \ECR\ $(M,\gamma,E,\Lambda)$ is said to satisfy the \emph{Dominant Energy Condition (\DEC)} if 
\begin{align}\label{eq:DEC}
\operatorname{R}(\gamma)&\geq 2\Lambda+n(n-1)\abs{E}_\gamma^2\text{ on }M
\end{align}
where $\operatorname{R}(\gamma)$ denotes the scalar curvature of $\gamma$ and $\vert E\vert_{\gamma}$ denotes the length of $E$ with respect to $\gamma$. It is said to satisfy the \emph{strict Dominant Energy Condition (\sDEC)} if 
\begin{align}\label{eq:sDEC}
\operatorname{R}(\gamma)&> 2\Lambda+n(n-1)\abs{E}_\gamma^2\text{ on }M.
\end{align}
\end{definition}

The \RN s $(M_{m,q,\Lambda},\gamma_{m,q,\Lambda},E_{m,q,\Lambda},\Lambda)$ introduced in Section~\ref{sec:RN} are standard examples of \ECR s saturating the \DEC, see Proposition~\ref{prop:RN}. The \sDEC\ leaves room for mollification, see Section~\ref{sec:gluing}.

\begin{remark}[Asymptotics]\label{rem:asymptotics}
When considering \ECR s $(M,\gamma,E,\Lambda)$ of dimension $n+1\geq3$, one usually assumes that they are subject to precise asymptotic conditions outside some compact set. More precisely, if $\Lambda=0$, one typically assumes that $(M,\gamma)$ is asymptotically Euclidean and $E$ decays to zero near infinity suitably fast. If $\Lambda<0$, one typically assumes that $(M,\gamma)$ is asymptotically hyperbolic with asymptotic hyperbolic radius $\sqrt{-\tfrac{n(n+1)}{2\Lambda}}$ and again that $E$ decays to zero near infinity suitably fast. To write down these precise asymptotic conditions is tedious, in particular in the asymptotically hyperbolic case, and in fact will not be necessary here as we will only consider \ECR s which are isometric to a suitable asymptotic piece of a \RN\ outside some compact set. 

Here, by a ``suitable asymptotic piece'' of a \RN, we mean a submanifold of the form $(\mu,\infty)\times\mathbb{S}^{n}$ with $\mu>r_{+}$ in the sub-extremal and extremal\footnote{In the extremal case, the restriction $\mu>r_{+}$ is essential as the piece $(r_{+},\mu)\times\mathbb{S}^{n}$ of an extremal \RN\ is in fact an asymptotically cylindrical end, see also Remark~\ref{rem:notextend}.} case and $\mu>0$ in the super-extremal case. It is well-known that such suitable asymptotic pieces are \emph{asymptotically Euclidean} if $\Lambda=0$ and \emph{asymptotically hyperbolic} with asymptotic hyperbolic radius $\sqrt{-\tfrac{n(n+1)}{2\Lambda}}$ if $\Lambda<0$ regardless of which precise definitions of these terms are being used. We will hence be referring to \ECR s which are isometric to a suitable asymptotic piece of a \RN\ $(M_{m,q,\Lambda},\gamma_{m,q,\Lambda},E_{m,q,\Lambda},\Lambda)$ outside some compact set as being \emph{asymptotically Euclidean} if $\Lambda=0$ and \emph{asymptotically hyperbolic} if~$\Lambda<0$.
\end{remark}

\begin{remark}[(Total) Mass]\label{rem:mass}
There are precise notions of \emph{(total) mass} for both asymptotically Euclidean and asymptotically hyperbolic Riemannian manifolds $(M,\gamma)$ of dimension $n+1\geq3$, referred to as the \emph{(ADM-)mass}, see \cite{ADM,ADMBartnik}, and the \emph{(hyperbolic) mass}, see \cite{Herzlich,Wang}, respectively. It is well-known and can easily be checked that the \RN\ $(M_{m,q,\Lambda},\gamma_{m,q,\Lambda},E_{m,q,\Lambda},\Lambda)$ has ADM mass $m$ when $\Lambda=0$ and hyperbolic mass $m$ when $\Lambda<0$. As we will only work with \ECR s $(M,\gamma,E,\Lambda)$ that are isometric to a  \RN\ outside some compact set, and as the ADM- and hyperbolic mass only depend on the asymptotic region, we can hence speak of their mass without having to give a precise definition here.
\end{remark}

\begin{remark}[(Total Electric) Charge]\label{rem:charge}
Similarly, there is a precise notion of \emph{(total electric) charge} for \ECR s $(M,\gamma,E,\Lambda)$ of dimension $n+1\geq3$, see for example \cite{Dain} and the references cited therein. It is well-known and can easily be checked that the \RN\ $(M_{m,q,\Lambda},\gamma_{m,q,\Lambda},E_{m,q,\Lambda},\Lambda)$ has electric charge $q$. As we will only work with \ECR s $(M,\gamma,E,\Lambda)$ that are isometric to a  \RN\ outside some compact set, and as the electric charge only depends on the asymptotic region, we can hence speak of their charge without having to give a precise definition here.
\end{remark}

In general relativity, one typically considers one of two kinds of \ECR s, namely either geodesically complete ones or those with what is called a \emph{horizon} boundary, a rather straightforward generalization of the horizons we have encountered in the sub-extremal \RN s, see Proposition~\ref{prop:RNhorizon}.

\begin{definition}[Outermost Black Hole Horizons, Outward Minimizing Minimal Hypersurfaces]\label{def:horizons}
Let $n\geq2$ and let $(M,\gamma,E,\Lambda)$ be an $(n+1)$-dimensional \ECR\ with compact boundary $\partial M\neq\emptyset$. Then $\partial M$ is said to consist of \emph{(black hole) horizons} if it consists of minimal hypersurfaces, i.e., if its mean curvature $H$ vanishes, $H=0$. It is said to be a \emph{(strictly) outermost (black hole) horizon} or a \emph{(strictly) outward minimizing minimal hypersurface} if its volume $\vert\partial M\vert_{g}$ is (strictly) less than the volume of any hypersurface $\Sigma$ enclosing $\partial M$, $\vert\partial M\vert_{g}\leq (<)\, \vert\Sigma\vert_{\widehat{g}}$. Here, $\abs{\cdot}_{g}$ and  $\abs{\cdot}_{\widehat{g}}$ denote the volumes computed with respect to the volume forms $dV_{g}$ and $dV_{\widehat{g}}$ of the metrics $g$ and $\widehat{g}$ induced by $\gamma$ on $\partial M$ and $\Sigma$, respectively.
\end{definition}
\newpage
It is well-known that the horizons $\lbrace{r_{+}\rbrace}\times\mathbb{S}^{n}=\partial M_{m,q,\Lambda}$ of the sub-extremal \RN s are outward minimizing minimal hypersurfaces in this sense (compare Proposition~\ref{prop:RNhorizon}). To see this is an application of a well-known theorem which will also be useful for the extensions to be constructed in this paper. As we have not been able to find a written proof of said theorem, we will give one here for completeness; the method of proof is a calibration argument. To state it, we need to give the following definition.

\begin{definition}[(Strictly) Mean Convex Hypersurfaces]
Let $n\geq2$ and let $(M,\gamma)$ be an $(n+1)$-dimensional Riemannian manifold. Let $\Sigma\subset M$ be a closed hypersurface with unit normal $\nu$ and mean curvature $H$ with respect to $\nu$. Then $\Sigma$ is said to be \emph{(strictly) mean convex} if $H\geq (>)\, 0$ everywhere on $\Sigma$.
\end{definition}

\begin{theorem}[Strictly Mean Convex Foliations and Outward Minimizing Hypersurfaces]\label{thm:outwardminfoliation}
Let $n\geq2$ and let $(M,\gamma)$ be a closed $(n+1)$-dimensional Riemannian manifold with closed, connected minimal boundary~$\partial M$. Assume that $M$ is regularly foliated by closed, connected hypersurfaces $\Sigma_{\tau}$, $\tau\in [\tau_{o},\infty)$, such that $\partial M=\Sigma_{\tau_{o}}$. Let $\nu\in\Gamma(TM)$ be the unit vector field which is normal to all $\Sigma_{\tau}$, $\tau\in I$ and has $\nu\vert_{\partial M}$ pointing into $M$. If all $\Sigma_{\tau}$, $\tau\in (\tau_{o},\infty)$, are (strictly) mean convex with respect to $\nu\vert_{\Sigma_{\tau}}$ then $\partial M=\Sigma_{\tau_{o}}$ is (strictly) outward minimizing in $(M,\gamma)$.
\end{theorem}

\begin{proof}
By the foliation assumption, the map $\tau\colon M\to[\tau_{o},\infty)\colon x\mapsto \tau_{x}$ with $x\in\Sigma_{\tau_{x}}$ is well-defined. Hence, we can define the map $\mathcal{H}\colon M\to \R\colon x\mapsto H(\Sigma_{\tau_x})(x)=\diver_\gamma (\nu)\vert_{x}$, where $H(\Sigma_{\tau_x})$ denotes the mean curvature of $\Sigma_{\tau_x}$ with respect to the unit normal $\nu\vert_{\Sigma_{\tau_{x}}}$. Using the mean convexity assumption on all $\Sigma_{\tau}$, $\tau\in[\tau_{o},\infty)$, we know that $\mathcal{H}\geq0$ in $M$.

Now let $\widehat{\Sigma}\subset M\setminus\partial M$ be an arbitrary closed, oriented hypersurface in $M$ enclosing $\partial M$ and let $\widehat{g}$ and $g$ denote the metrics induced on $\widehat{\Sigma}$ and $\partial M$, respectively. Set $\Omega\subseteq M$ to be the open domain satisfying $\partial\Omega=\widehat{\Sigma}\cup \partial M$. 
 Let $\widehat{\nu}$ be the unit normal  to $\widehat{\Sigma}$ pointing out of $\Omega$. Using the Divergence Theorem and the Cauchy--Schwarz Inequality, we obtain
\begin{align*}
0\leq\int_\Omega\mathcal{H}\, dV_\gamma
&= \int_\Omega \diver_\gamma(\nu)\,d V_\gamma = \int_{\widehat{\Sigma}} \gamma(\widehat{\nu},\nu)\,dV_{\widehat{g}}-\int_{\partial M} \gamma(\nu,\nu)\,dV_g \leq \vert\widehat{\Sigma}\vert_{\widehat{g}}-\abs{\partial M}_g,
\end{align*}
where $dV_{\gamma}$ denotes the volume measure induced on $\Omega$ by $\gamma$ and $dV_{\widehat{g}}$ and $dV_{g}$ denote the volume measures induced on $\widehat{\Sigma}$ and $\Sigma_{\tau_{o}}$, respectively.
Thus, we get $\abs{\partial M}_g \leq \vert\widehat{\Sigma}\vert_{\widehat{g}}$. Relaxing the assumption $\widehat{\Sigma}\subset M\setminus\partial M$ to $\widehat{\Sigma}\subset M$, or in other words allowing $\widehat{\Sigma}$ to touch and possibly (partially) coincide with $\partial M$, we can apply the above conclusion to $\widehat{\Sigma}_{\varepsilon}\definedas \lbrace \exp_{p}(\varepsilon\widehat{\nu}_{p})\,\vert\,p\in \widehat{\Sigma}\rbrace$ for suitably small $0<\varepsilon<\varepsilon(\widehat{\Sigma})$ and obtain $\abs{\partial M}_g \leq \vert\widehat{\Sigma}_{\varepsilon}\vert_{\widehat{g}_{\varepsilon}}$ so that $\abs{\partial M}_g \leq \vert\widehat{\Sigma}\vert_{\widehat{g}}$ follows by continuity as $\varepsilon\searrow0$. This shows that $\partial M$ is outward minimizing in $(M,\gamma)$.

In case the $\Sigma_{\tau}$, $\tau\in(\tau_{o},\infty)$, are strictly outward minimizing, we have that $\mathcal{H}>0$ on $M\setminus\partial M$, hence $\int_{\Omega}\mathcal{H}\,dV_{\gamma}>0$ unless $\Omega=\emptyset$ (or $\bigcap_{0<\varepsilon<\varepsilon_{o}} \Omega_{\varepsilon}=\emptyset$ if $\partial M\cap\widehat{\Sigma}\neq\emptyset$) or in other words unless $\widehat{\Sigma}=\partial M$. This shows that $\partial M$ is strictly outward minimizing in $(M,\gamma)$.
\end{proof}

The most important theorem on geodesically complete \ECR s is the Riemannian Positive Mass Theorem which we will state here somewhat imprecisely. In the way we word it here, it follows from the fact that the term $n(n-1)\vert E\vert_{\gamma}^{2}$ in the \DEC~\eqref{eq:DEC} is non-negative which implies the uncharged Dominant Energy Condition $\scal(\gamma)\geq 2\Lambda$ assumed in the Riemannian Positive Mass Theorem.

\begin{theorem}[Riemannian Positive Mass Theorem~\cite{SY1,Witten,Wang,Herzlich,Cai,SY3,MR3801943,Delay,HJM,SakovichPMT}]
Let $(M,\gamma)$ be a geodesically complete Riemannian manifold of dimension $n+1\geq3$, let $\Lambda\leq0$, and assume that $(M,\gamma)$ is asymptotically Euclidean if $\Lambda=0$ and asymptotically hyperbolic with asymptotic hyperbolic radius $\sqrt{-\tfrac{n(n+1)}{2\Lambda}}$ if $\Lambda<0$, respectively. Furthermore, assume that $\scal(\gamma)\geq 2\Lambda$ holds on $M$. Then the mass $m$ of $(M,\gamma)$ satisfies $m\geq0$, with equality if and only if $(M,\gamma)$ is globally isometric to Euclidean space for $\Lambda=0$ and globally isometric to the hyperbolic space of radius $\sqrt{-\tfrac{n(n+1)}{2\Lambda}}$ for $\Lambda<0$.
\end{theorem}

In the presence of an electric field $E$, the Riemannian Positive Mass Theorem can be refined to the Riemannian Mass-Charge Inequality $m\geq \abs{q}$ when $\Lambda=0$, see \cite[Sec.~3.1.1]{Dain} or \cite{Weinstein2} and the references cited therein for a precise statement. To the best knowledge of the authors, this has not yet been generalized to include the case $\Lambda<0$.

A more refined but similar inequality has been shown (or conjectured, in some cases) for Riemannian manifolds or \ECR s with outward minimizing minimal boundary. Let us first state the conjectured version of this Riemannian Penrose Inequality with charge that we will allude to in this paper.

\begin{conjecture}[Generalized Riemannian Penrose Inequality with Charge and Cosmological Constant]\label{conj:RPI}
Let $(M,\gamma,E,\Lambda)$ be an \ECR\ of dimension $n+1\geq3$ which is geodesically complete up to its closed, strictly outward minimizing minimal boundary $\partial M$. Assume that $(M,\gamma,E,\Lambda)$ is asymptotically Euclidean if $\Lambda=0$ and asymptotically hyperbolic with asymptotic hyperbolic radius $\sqrt{-\tfrac{n(n+1)}{2\Lambda}}$ if $\Lambda<0$, respectively, and denote its mass and charge by $m$ and $q$, respectively. If $q\neq0$, assume in addition that $\partial M$ is connected. Furthermore, assume that the Dominant Energy Condition~\eqref{eq:DEC} holds on $M$. Then $(M,\gamma,E,\Lambda)$ satisfies the \emph{generalized Riemannian Penrose Inequality (with charge)}
\begin{align}\label{eq:RPI}
 \frac{1}{2}\left[\left(\frac{\abs{\partial M}_{g}}{\omega_n}\right)^{\frac{n-1}{n}}+q^2\left(\frac{\abs{\partial M}_{g}}{\omega_n}\right)^{-\frac{n-1}{n}}-\frac{2\Lambda}{n(n+1)}\left(\frac{\abs{\partial M}_{g}}{\omega_n}\right)^{\frac{n+1}{n}}\right]\leq m,
\end{align}
where $\abs{\partial M}_{g}$ denotes the volume of $\partial M$ with respect to the induced metric $g$ on $\partial M$ and $\omega_n=\abs{\mathbb{S}^{n}}_{g_{*}}$ denotes the volume of the unit round sphere $(\mathbb{S}^{n},g_{*})$. Equality in~\eqref{eq:RPI} holds if and only if $(M,\gamma,E,\Lambda)$ is isometric to the necessarily sub-extremal \RN\ $(\overline{M_{m,q,\Lambda}},\gamma_{m,q,\Lambda},E_{m,q,\Lambda},\Lambda)$ of mass $m$, charge $q$, and cosmological constant $\Lambda$.
\end{conjecture}

The fact that the sub-extremal \RN s satisfy~\eqref{eq:RPI} can easily be seen from $\abs{\partial M}_{g}=\omega_{n}r_{+}^{n}$ and $p_{m,q,\Lambda}(r_{+})=0$, see Proposition~\ref{prop:extendRN}. 

\begin{remark}[Connectedness]
We impose a connectedness assumption on $\partial M$ in Conjecture~\ref{conj:RPI} in case the total charge $q\neq0$ because of a known counter-example to the $n+1=3$-dimensional Riemannian Penrose Inequality with charge with disconnected horizon boundary, see Weinstein and Yamada~\cite{WY}. See Khuri, Weinstein, and Yamada~\cite{Khuri2} for a partial extension of the Riemannian Penrose Inequality with charge (for $\Lambda=0$, $n+1=3$) to the case of disconnected boundary, using a charged Conformal Flow.
\end{remark}

\begin{remark}[Known Cases]\label{rem:RPI}
The Riemannian Penrose Inequality (with charge)~\eqref{eq:RPI} has been proven at least in the following cases: 
\begin{enumerate}
\item $\Lambda=0$, $E=0$: The first proofs were famously given for $n+1=3$ by Huisken and Ilmanen~\cite{RPI} in the case of connected boundary using (weak) Inverse Mean Curvature Flow, and by Bray~\cite{RPI2} without the connectedness assumption, using Conformal Flow. Bray and Lee~\cite{RPI3} later generalized the Conformal Flow approach to higher dimensions. Lam~\cite{Lam} gave a proof of the inequality in higher dimensions under the assumption that $(M,\gamma)$ is a graph in Euclidean space, while Huang and Wu obtained the corresponding rigidity statement \cite{HuangWu}.
\item $\Lambda=0$, $E\neq0$: Jang \cite{Jang} and Disconzi and Khuri~\cite{Discon} (see also McCormick \cite{MccormickQ} who allows for non-vanishing charge density) adapted the (weak) Inverse Mean Curvature Flow approach to include an electric field, and Khuri, Weinstein, and Yamada~\cite{Khuri2} adapted the Conformal Flow approach to include an electric field, both for $n+1=3$.
\item $\Lambda<0$, $E=0$: A proof by Dahl, Gicquaud, and Sakovich \cite{Dahl}, complemented by de Lima and Gir\~{a}o \cite{lima1}, has been given in the graphical case. Other cases have been established by de Lima, Gir\~{a}o, Loz\'orio, and Silva~\cite{lima} and by Ambrozio~\cite{Ambrozio}.
\item $\Lambda<0$, $E\neq0$: A proof was given for graphs over sub-extremal \RN s by Chen, Li, and Zhou~\cite{Chen}.
\end{enumerate}		
We refer the interested reader to the reviews by Mars \cite{Mars} and by Bray and Chru\'{s}ciel~\cite{RPIsurvey} and the articles cited above for more information on the Riemannian Penrose Inequality and its context in general relativity.
\end{remark}

\subsection{A Generalized Notion of Hawking Mass}\label{subsec:Hawking}
Before we begin constructing extensions of minimal boundary spheres, let us introduce a notion of quasi-local mass which mimics the definition of Hawking (or Geroch) mass~\cite{Hawking,Geroch} and shares some important properties with it. To derive it, we proceeded in a similar way as Disconzi and Khuri~\cite{Discon}. First, let us define the quasi-local electric charge.

\begin{definition}[Quasi-Local (Electric) Charge]\label{def:QLcharge}
Let $(M,\gamma,E,\Lambda)$ be an \ECR\ of dimension $n+1\geq3$ and let $\Sigma\subset M$ be a closed hypersurface with unit normal $\nu$. Then the \emph{quasi-local (electric) charge $\mathcal{Q}(\Sigma)$ (with respect to $\nu$)} is defined as
\begin{align}
\mathcal{Q}(\Sigma)&\definedas \frac{1}{\omega_{n}}\int_{\Sigma}\gamma(E,\nu)dV_{g},
\end{align}
where $dV_{g}$ denotes the volume measure on $\Sigma$ with respect to the induced metric $g$.
\end{definition}

\begin{remark}[Quasi-Local Charge in Rotational Symmetry]\label{rem:QLchargeSS}
Consider an \ECR\ of the form $(I \times \mathbb{S}^n , \gamma=ds^2 +f(s)^2 g_* , E, \Lambda)$, where $I$ is an interval, $f\colon I\to\R^{+}$ is a smooth, positive function, and $g_{*}$ denotes the canonical metric on $\mathbb{S}^{n}$. Then every coordinate sphere $\{s\} \times \mathbb{S}^n$, $s\in I$, has the same quasi-local charge $\mathcal{Q}(\{s\} \times \mathbb{S}^n)=\mathcal{Q}$ by the Divergence Theorem and as $\diver_{\gamma}E=0$.
\end{remark}

\begin{definition}[Generalized Hawking Mass]\label{def:Hawking}
Let $(M,\gamma,E,\Lambda)$ be an \ECR\ of dimension $n+1\geq3$. Let $\Sigma\subset M$ be a closed hypersurface with unit normal $\nu$. Let $g$ denote the induced metric on $\Sigma$ and let $H$ be the mean curvature of $\Sigma$ with respect to $\nu$. The \emph{generalized Hawking mass $\mathcal{M}(\Sigma)$ of $\Sigma$ in $(M,\gamma,E,\Lambda)$} is defined as
\begin{align}
\begin{split}
\mathcal{M}(\Sigma)& \definedas \frac{1}{2} \left( \frac{|\Sigma|_g}{\omega_n} \right)^{\frac{n-1}{n}}\times\\
&\,\left(1- \frac{1}{n^2 \omega_n} \left( \frac{\omega_n}{ |\Sigma|_g} \right)^{\frac{n-2}{n}}\! \int_{\Sigma }H^2dV_{g} + \mathcal{Q}(\Sigma)^2 \left( \frac{\omega_n}{ |\Sigma|_g} \right)^{\frac{2(n-1)}{n}}\!-\frac{2\Lambda}{n(n+1)} \left( \frac{|\Sigma|_g}{\omega_n} \right)^{\frac{2}{n}}\right),
\end{split}
\end{align}
where $dV_{g}$ denotes the volume measure on $\Sigma$.
\end{definition}

\begin{remark}[Consistency with Previous Definitions]\label{rem:consistent}
As can easily be verified, the generalized Hawking mass $\mathcal{M}$ coincides with the original definition for $n=2$ in the case $\Lambda=0$ without electric field (see \cite{Hawking,Geroch}), with the ``asymptotically hyperbolic Hawking mass'' for $n=2$ and without electric field when $\Lambda=\Lambda_{n}=-\frac{n(n+1)}{2}$ (see~\cite{Neves,hyper}), with the ``charged Hawking mass'' (see~\cite{Discon,Jang}) for $n=2$ and $\Lambda=0$, and with the higher dimensional Hawking mass for $n\geq3$, $\Lambda=0$, and without electric field (see for example~\cite{Bryden}).
\end{remark}

\begin{proposition}[Rewriting the Generalized Hawking Mass]\label{prop:Hawkingre}
Using the functions $p_{m,q,\Lambda}$ introduced in~\eqref{eq:defp} and the volume radius $r(\Sigma)\definedas \left(\tfrac{\vert\Sigma\vert_{g}}{\omega_{n}}\right)^{\frac{1}{n}}$, we can rewrite the generalized Hawking mass as
\begin{align}\label{eq:Hawkingsimple}
\mathcal{M}(\Sigma)& = \frac{r(\Sigma)^{n-1}}{2} \left( p_{0,\mathcal{Q}(\Sigma),\Lambda}(r(\Sigma))- \frac{1}{n^2 \omega_n r(\Sigma)^{n-2}} \int_{\Sigma }H^2dV_{g} \right).
\end{align}
\end{proposition}

\begin{proposition}[Quasi-Local (Sub-)Extremality of Minimal Hypersurfaces]\label{prop:QLsub}
Let $n\geq2$ and let $(M,\gamma,E,\Lambda)$ be an \ECR\ of dimension $n+1$. Let $\Sigma\subset M$ be a closed, minimal hypersurface in $(M,\gamma)$, and let $g$ denote its induced metric. Let $r(\Sigma)$ denote the volume radius of $\Sigma$, i.e., $\vert\Sigma\vert_{g}\asdefined\omega_{n}r(\Sigma)^{n}$. Then the \RN\ of mass $\mathcal{M}(\Sigma)$, charge $\mathcal{Q}(\Sigma)$, and cosmological constant $\Lambda$ is sub-extremal with largest positive root $r_{+}=r(\Sigma)$ if and only if $h_{\mathcal{Q}(\Sigma),\Lambda}(r(\Sigma))>0$ and extremal if and only if $h_{\mathcal{Q}(\Sigma),\Lambda}(r(\Sigma))=0$. 

Moreover, if it is sub-extremal then $\mathcal{M}(\Sigma)>\vert\mathcal{Q}(\Sigma)\vert$ while if it is extremal then $\mathcal{Q}(\Sigma)\neq0$.
\end{proposition}

\begin{proof}
By definition of $\mathcal{M}(\Sigma)$ and the fact that $\Sigma$ is minimal, we find
\begin{align*}
\mathcal{M}(\Sigma)&=\frac{r(\Sigma)^{n-1}}{2}\left(1+\frac{\mathcal{Q}(\Sigma)^{2}}{r(\Sigma)^{2(n-1)}}-\frac{2\Lambda r(\Sigma)^{2}}{n(n+1)}\right),
\end{align*}
hence $r(\Sigma)$ is a positive root of $p_{\mathcal{M}(\Sigma),\mathcal{Q}(\Sigma),\Lambda}$. This shows that the \RN\ of mass $\mathcal{M}(\Sigma)$, charge $\mathcal{Q}(\Sigma)$, and cosmological constant $\Lambda$ must be either sub-extremal or extremal. By Lemma~\ref{lem:sub}, we know that if it is extremal then $p_{\mathcal{M}(\Sigma),\mathcal{Q}(\Sigma),\Lambda}$ has precisely one positive root $r_{+}=r(\Sigma)$ with $h_{\mathcal{Q}(\Sigma),\Lambda}(r(\Sigma))=0$. Similarly, if it is sub-extremal, then $p_{\mathcal{M}(\Sigma),\mathcal{Q}(\Sigma),\Lambda}$ has at most two positive roots $0<r_{-}<r_{+}$ with $h_{\mathcal{Q}(\Sigma),\Lambda}(r_{-})<0<h_{\mathcal{Q}(\Sigma),\Lambda}(r_{+})$. Hence $h_{\mathcal{Q}(\Sigma),\Lambda}(r(\Sigma))>0$ implies $r(\Sigma)=r_{+}$ also in this case. The remaining claims follow directly from Lemma~\ref{lem:sub}.
\end{proof}

The generalized Hawking mass naturally simplifies when studied in rotationally symmetric \ECR s. Moreover, it has very useful monotonicity properties along the coordinate spheres in rotationally symmetric \ECR s, but also, for $n=2$, more generally along collar extensions such as those to be constructed in Section~\ref{sec:collar}. These properties will be discussed in Section~\ref{sec:Hawkingcollar}.

\section{Collar extensions}\label{sec:collar}
This section will be devoted to constructing ``collar extensions'', a central tool for constructing extensions with prescribed asymptotics. More precisely, for a given smooth metric $g_o$ on the sphere $\mathbb{S}^n$, we will construct a \emph{collar extension}, that is, a \ECR\ $(M,\gamma,E,\Lambda)$ such that $\partial M$, with its induced metric $g$, is isometric to $(\mathbb{S}^{n},g_{o})$ and is a strictly outward minimizing hypersurface in $(M,\gamma)$. We intend to use Theorem~\ref{thm:outwardminfoliation} to show that $\partial M$ is outward minimizing.

Searching for collar extensions, we make the ansatz that $M=[0,1] \times \mathbb{S}^n$ arises as a direct product, so that in particular $\partial M=\lbrace0\rbrace\times\mathbb{S}^{n}$, while $\gamma$ shall be of the form
\begin{align}\label{eq:ansatzgamma}
\gamma &=v(t,\cdot)^2 dt^2 +F(t)^2 g(t),
\end{align}
where $v\colon M=[0,1] \times \mathbb{S}^n\to\R^{+}$ and $F\colon[0,1]\to\R^{+}$ are smooth, positive functions and $\{g(t) \}_{t\in[0,1]}$ is a suitable smooth path of metrics connecting $g_o$ to a round metric on $\mathbb{S}^n$.
As usual when constructing extensions \`a la Mantoulidis and Schoen, the smooth path of metrics $\{g(t) \}_{t\in[0,1]}$ will be chosen such that it satisfies the following Conditions (E1-E3):
\begin{enumerate}\label{page:conditions}
\item[(E1)] $g(0)=g_o$ and $g(1)$ is round, \label{extension:1}
\item[(E2)] $g'(t)=0$ for $t \in [\theta, 1]$ where $0<\theta <1$, and \label{extension:2}
\item[(E3)] $\tr_{g(t)} g'(t) =0$ for $t \in [0,1]$,\label{extension:3}
\end{enumerate}
where $'$ abbreviates $\frac{\partial}{\partial t}$. Here, Condition \hyperref[extension:1]{(E1)} ensures that $\partial M=\lbrace0\rbrace\times\mathbb{S}^{n}$, with its induced metric $g(0)$, is isometric to $(\mathbb{S}^{n},g_{o})$ if we fix $F(0)=1$. Moreover, together with Condition \hyperref[extension:2]{(E2)}, it will later allow us to smoothly glue the collar 
extension to a rotationally symmetric \RN. Condition \hyperref[extension:3]{(E3)} is helpful for simplifying the formula for the scalar curvature of $(M,\gamma)$ and is thereby useful for ensuring that the \DEC~\eqref{eq:DEC} will be satisfied by the \ECR\ we construct. Moreover, it simplifies the formula for the mean curvature of the coordinate spheres $\lbrace t\rbrace\times\mathbb{S}^{n}$, see Lemma~\ref{lem:meancurvature}, and thereby helps ensuring both that $\lbrace0\rbrace\times\mathbb{S}^{n}$ is minimal and that $\lbrace t\rbrace\times\mathbb{S}^{n}$ is strictly mean convex for $t\in(0,1]$. It implies in particular that the volume $\abs{\lbrace t\rbrace\times\mathbb{S}^{n}}_{F(t)^{2}g(t)}$ is preserved (i.e., independent of $t$).

\begin{remark}\label{rem:vol0}
It can easily be seen that Condition~\hyperref[extension:3]{(E3)} is equivalent to $\frac{d}{dt}dV_{g(t)}=0$ for $t \in [0,1]$, where $dV_{g(t)}$ denotes the volume form induced on $\mathbb{S}^{n}$ by the metric $g(t)$.
\end{remark}

Of course, when $n>2$, one needs to assume or otherwise ensure a priori that $g_{o}$ can be connected to a round metric on $\mathbb{S}^{n}$ in a smooth fashion.  The assumptions we will make will depend on the dimensional parameter $n$ and will be discussed below in Sections~\ref{sec:n2},~\ref{sec:n3},~\ref{sec:nn}. Indeed, once a smooth path of metrics $\{h(\tau) \}_{\tau\in [0,T)}$, $0<T\leq\infty$, connecting $g_{o}=h(0)$ to a round metric on $\mathbb{S}^{n}$ as $\tau\nearrow T$ has been constructed, a path $\{g(t) \}_{t\in[0,1]}$ satisfying Conditions \hyperref[page:conditions]{(E1-E3)} can be constructed by rescaling the time $\tau$, adjoining the round limit metric $\lim_{\tau\nearrow T}g(\tau)$, dilating the metrics in the path smoothly depending on time to preserve their volume, and applying a diffeomorphism fixing procedure developed by Mantoulidis and Schoen for $n+1=3$ (see proof of Lemma 1.2 in \cite{pri}) and generalized to higher dimensions by Cabrera~Pacheco and Miao (see proof of Lemma 4.1 in \cite{sec}) to ensure Condition~\hyperref[extension:3]{(E3)}. Moreover, positive scalar curvature can be preserved in this procedure. We summarize this as follows.

\begin{lemma}[Normalizing Paths, \cite{pri,sec}]\label{lem:normalizing}
Let $n\geq2$ and let $g_{o}$ be a smooth metric on $\mathbb{S}^{n}$. Suppose there is a smooth path of metrics $\{h(\tau) \}_{\tau\in [0,T)}$, $0<T\leq\infty$, connecting $g_{o}=h(0)$ to a round metric on $\mathbb{S}^{n}$ as $\tau\nearrow T$. Then this path can be transformed to a smooth path of metrics $\{g(t) \}_{t\in[0,1]}$ satisfying Conditions \hyperref[page:conditions]{(E1-E3)}. If all metrics in the path $\{h(\tau) \}_{\tau\in [0,T)}$ have positive scalar curvature $\scal(h(\tau))>0$ then it can be arranged that also $\scal(g(t))>0$ for all $t\in[0,1]$.
\end{lemma}

Before we discuss the assumptions we need to make on $g_{o}$, let us quickly mention how we will construct the remaining ingredients for constructing the \ECR\ constituting our collar extension. We will use the \emph{volume radius} $r_o$ defined by $|\mathbb{S}^n|_{g_o} \asdefined \omega_n r_o^n $, where $|\mathbb{S}^n|_{g_o}$ is the volume of the unit sphere with respect to the metric $g_o$ and $\omega_n$ denotes the volume of the standard unit $n$-sphere. For the electric field $E$, we make the ansatz
\begin{align}\label{eq:ansatzE}
E&=\frac{q}{v(t,\cdot)r_{o}^{n}F(t)^{n}}\partial_{t},
\end{align}
where $q\in\R$ is an a priori arbitrary constant, so that $E$ is naturally divergence-free. When we will later glue the constructed collar extensions to suitable sub-extremal \RN s, $q$ will play the role of (total) charge and we will impose conditions on it as needed. As the glued extensions shall satisfy the \DEC~\eqref{eq:DEC}, we will need to pick a designated cosmological constant $\Lambda$ in advance.

We will now discuss the different assumptions we will make on the initial metric $g_{o}$ on $\mathbb{S}^{n}$ and collect the relevant theorems allowing to construct a smooth path of metrics $\{g(t) \}_{t\in[0,1]}$ satisfying said assumptions along the path.
	
\subsection{Case $n=2$ and a Conformal Approach for $n\geq3$}\label{sec:n2}
The first option we have is to consider paths of metrics in the set $\mathcal{M}^{+} \definedas\{g\text{ metric on } \mathbb{S}^2 : \lambda_1 (-\triangle_g + K(g)) >0 \}$, where $\lambda_1 (-\triangle_g + K(g))$ is the first eigenvalue of the operator $-\triangle_g + K(g)$ on $\mathbb{S}^{2}$, $K(g)$ is the Gaussian curvature,  $\Delta_g$ is the Laplace--Beltrami operator of $g$. The condition on $\lambda_1$ arises naturally in the theory of minimal surfaces, and its connection to general relativity comes from the fact that outermost horizons (aka outward minimizing minimal surfaces) in asymptotically Euclidean Riemannian manifolds with non-negative scalar curvature are, in fact, stable minimal spheres (for more details see \cite{impor,hyper, pri,LM}). It was shown by Mantoulidis and Schoen~\cite{pri} that $\mathcal{M}^{+}$ is path connected and that we can always find a smooth path of metrics connecting a given metric on $\mathcal{M}^{+}$ to the round metric as desired\footnote{Chau and Martens~\cite{ACAM} generalized this result to $\lambda_1(-\triangle_g+K(g))\geq 0$; we do not pursue this case here.}

\begin{lemma}[{\cite[Lemma 1.2]{pri}}]\label{le1.2}
For any $g_{o}\in \mathcal{M}^{+}$, there exists a smooth path of metrics $\lbrace g(t)\rbrace_{t\in[0,1]} \subset \mathcal{M}^{+}$ satisfying Conditions \hyperref[page:conditions]{(E1-E3)}. In particular, if $K(g_{o})>0$ then the whole path $\{g(t)\}_{t\in[0,1]}$ has $K(g(t))>0$.
\end{lemma}

The construction of this path of metrics relies on the Uniformization Theorem. In fact, it is easy to see that the same approach works for $n\geq3$ provided that one picks $g_{o}\in[g_{*}]$, where $[g_{*}]$ denotes the conformal class of the canonical metric $g_{*}$ on $\mathbb{S}^{n}$. To see this, we set
$\mathcal{M}^{+}_{[g_{*}]} \definedas\{g\in[g_{*}] : \lambda_1 (-\triangle_g + \tfrac{1}{2}\scal(g)) >0 \}$, where $\lambda_1 (-\triangle_g + \tfrac{1}{2}\scal(g))$ is the first eigenvalue of the operator $-\triangle_g + \tfrac{1}{2}\scal(g)$ on $\mathbb{S}^{n}$, with $\Delta_g$ denoting the Laplace--Beltrami operator and $\scal(g)$ the scalar curvature of $g$. Then for $g_{o}\in\mathcal{M}^{+}_{[g_{*}]}$, we know that there exists a smooth function $w\colon\mathbb{S}^{n}\to\R$ and a diffeomorphism $\Phi\colon\mathbb{S}^{n}\to\mathbb{S}^{n}$ such that $g_{o}=e^{2w}\Phi_{*}g_{*}$. Setting $h(t)\definedas e^{2(1-t)w}\Phi_{*}g_{*}$, $t\in[0,1]$, we obtain a smooth path of metrics in $[g_{*}]$ such that $h(0)=g_{o}$ and $h(1)$ is round. To see that $h(t)\in\mathcal{M}^{+}_{[g_{*}]}$ for all $t\in[0,1]$, let $f\colon\mathbb{S}^{n}\to\R$ be a smooth function and compute
\begin{align*}
&\int_{\mathbb{S}^{n}}f(-\triangle_{h(t)} + \tfrac{1}{2}\scal(h(t)))f\,dV_{h(t)}=\int_{\mathbb{S}^{n}}(\vert df\vert^{2}_{h(t)} + \tfrac{1}{2}\scal(h(t))f^{2})\,dV_{h(t)}\\
&=\int_{\mathbb{S}^{n}}e^{(n-2)(1-t)w}\left(\vert df\vert^{2}_{\Phi_{*}g_{*}} + \tfrac{n-1}{2}\left[n-2(1-t)\triangle_{\Phi_{*}g_{*}}w-(n-2)(1-t)^{2}\vert dw\vert^{2}_{\Phi_{*}g_{*}}\right]f^{2}\right)dV_{\Phi_{*}g_{*}}\\
&=\int_{\mathbb{S}^{n}}\left(\vert d(e^{\frac{(n-2)(1-t)w}{2}}f)\vert^{2}_{\Phi_{*}g_{*}}-(n-2)(1-t)(\Phi_{*}g_{*})(dw,df)e^{(n-2)(1-t)w}f\right.\\
&\quad\quad\left.+\left\lbrace-\frac{(n-2)^{2}(1-t)^{2}}{4}\vert dw\vert^{2}_{\Phi_{*}g_{*}} + \frac{n-1}{2}\left[n-(n-2)(1-t)^{2}\vert dw\vert^{2}_{\Phi_{*}g_{*}}\right]\right\rbrace e^{(n-2)(1-t)w}f^{2}\right.\\
&\quad\quad\left.+(n-1)\left[(n-2)(1-t)^{2}\vert dw\vert_{\Phi_{*}g_{*}}^{2}f^{2}+2(1-t)(\Phi_{*}g_{*})(dw,df)f\right]e^{(n-2)(1-t)w}\right)dV_{\Phi_{*}g_{*}},
\end{align*}
where $dV_{h(t)}$ denotes the volume form on $\mathbb{S}^{n}$ with respect to $h(t)$, $dV_{\Phi_{*}g_{*}}$ the one with respect to $\Phi_{*}g_{*}$, and where we have used the Divergence Theorem, $\scal(g_{*})=n(n-1)$, the standard formula for the transformation of the scalar curvature under conformal changes, and $dV_{h(t)}=e^{n(1-t)w}dV_{\Phi_{*}g_{*}}$. Now set $\psi\definedas e^{\frac{(n-2)(1-t)w}{2}}f$ and observe that this can be rewritten as
\begin{align*}
&\int_{\mathbb{S}^{n}}f(-\triangle_{h(t)} + \tfrac{1}{2}\scal(h(t)))f\,dV_{h(t)}\\
&=\int_{\mathbb{S}^{n}}\vert d\psi\vert^{2}_{\Phi_{*}g_{*}}dV_{\Phi_{*}g_{*}}+\frac{n(n-1)}{2}\int_{\mathbb{S}^{n}}\psi^{2}\,dV_{\Phi_{*}g_{*}}+(1-t)\int_{\mathbb{S}^{n}}n(\Phi_{*}g_{*})(dw,d\psi)\psi\, dV_{\Phi_{*}g_{*}}
\\
&\quad\quad-(1-t)^{2}\int_{\mathbb{S}^{n}}\frac{n(n-2)\vert dw\vert_{\Phi_{*}g_{*}}^{2}}{4}dV_{\Phi_{*}g_{*}}\\
&\asdefined B_{\psi}(t),
\end{align*} 
using again the Divergence Theorem. In order to see that $\lambda_{1}(-\triangle_{h(t)}+\tfrac{1}{2}\scal(h(t)))>0$ for all $t\in[0,1]$, it hence suffices to show that $B_{\overline{\psi}}(t)>0$ for all $t\in[0,1]$ and all smooth, non-vanishing $\overline{\psi}\colon\mathbb{S}^{n}\to\R$. Fix a smooth, non-vanishing $\overline{\psi}\colon\mathbb{S}^{n}\to\R$. Then $B_{\overline{\psi}}(1)$ is manifestly positive, while $B_{\overline{\psi}}(0)$ is positive by our assumption that $g_{o}\in\mathcal{M}^{+}_{[g_{*}]}$. On the other hand, the map $t\mapsto B_{\overline{\psi}}(t)$ is an at most quadratic polynomial in $t$ which is manifestly concave, i.e., satisfies $0<(1-t) B_{\overline{\psi}}(0)+tB_{\overline{\psi}}(1)\leq B_{\overline{\psi}}(t)$ for all $t\in[0,1]$. Hence $h(t)\in \mathcal{M}^{+}_{[g_{*}]}$ for all $t\in[0,1]$ as claimed.

In particular, if $\scal(g_{o})>0$, the above argument applies without the term $\int_{\mathbb{S}^{n}}\vert df\vert_{h(t)}^{2}\,dV_{h(t)}$ and then shows that $\scal(g(t))>0$ for all $t\in[0,1]$ (cf. \cite{Cabrera-thesis}). Combining this with Lemma~\ref{lem:normalizing}, we obtain the following conclusion.

\begin{lemma}[Paths from Uniformization, $n\geq3$]\label{lem:uniform}
For any $g_{o}\in \mathcal{M}^{+}_{[g_{*}]}$, there exists a smooth path of metrics $\lbrace g(t)\rbrace_{t\in[0,1]} \subset \mathcal{M}^{+}_{[g_{*}]}$ satisfying Conditions \hyperref[page:conditions]{(E1-E3)}.
In particular, if $\scal(g_{o})>0$ then $\scal(g(t))>0$ for all $t\in[0,1]$.
\end{lemma}

Moreover, the above comment about $\lambda_{1}>0$ being related to stable outward minimizing minimal spheres applies analogously for $n\geq3$.

Returning our attention to $n=2$, one can alternatively construct a path for an initial metric $g_{o}$ with $\scal(g)$ or in other words $K(g)$ bounded below by a negative constant as was shown by Cabrera~Pacheco, Cederbaum, and McCormick~\cite{hyper} using area-preserving Ricci flow (see \cite{Chow,Hamilton2}). This condition turned out to be particularly useful to construct the asymptotically hyperbolic extensions in \cite{hyper}.
\begin{lemma}[{\cite[Lemma 5.1]{hyper}}]\label{path2}
Given $(\mathbb{S}^2,g_o)$ with $K(g_o)>-\tau$, where $\tau\geq0$, there exists a smooth path of metrics $\{g(t)\}_{t\in[0,1]}$ on $\mathbb{S}^{2}$ satisfying Conditions \hyperref[page:conditions]{(E1-E3)} and $K(g(t))>-\tau$ for all $t\in[0,1]$.
\end{lemma}

\subsection{Case $n=3$}\label{sec:n3}
For $n=3$, we are only aware of one result allowing to construct a smooth path of metrics connecting a given metric $g_{o}$ to a round metric on $\mathbb{S}^{3}$. Let $\operatorname{Scal}^+(\mathbb{S}^n)$ denote the set of smooth metrics with positive scalar curvature on $\mathbb{S}^n$. The space $\operatorname{Scal}^+(\mathbb{S}^n)$ was shown to be path connected for $n=3$ by Marques~\cite{Mar} using Ricci flow (see \cite{Hamil,Brendle}). Cabrera~Pacheco and Miao~\cite{sec} mollified such continuous paths of metrics and obtained the following theorem.
\begin{theorem}[{\cite[Corollary 2.1]{sec}}]\label{pat3}
Given any $g_{o} \in \operatorname{Scal}^+(\mathbb{S}^3) $, there exist a smooth path $\{h(t) \}_{t\in[0,1]} $ in $\operatorname{Scal}^+(\mathbb{S}^3) $ connecting $g_{o}$ to a round metric on $\mathbb{S}^3 $.
\end{theorem}

Applying Lemma~\ref{lem:normalizing}, they then obtain the desired path.

\begin{corollary}[Existence of Paths with Desired Properties, $n=3$]\label{coro:paths3}
Given $g_{o}\in \operatorname{Scal}^+(\mathbb{S}^3)$ with $\scal(g_{o})>0$, there exists a smooth path of metrics $\{g(t)\}_{t\in[0,1]}$ on $\mathbb{S}^{3}$ satisfying Conditions~\hyperref[page:conditions]{(E1-E3)} and $\scal(g(t))>0$ for all $t\in[0,1]$.
\end{corollary}

\subsection{Case $n>3$}\label{sec:nn}
To apply Ricci flow successfully in higher dimensions $n\geq 4$, we have to assume in addition that the sphere has pointwise $\nicefrac{1}{4}$-pinched sectional curvatures in view of the Differential Sphere Theorem by Brendle and Schoen~\cite{sphere} proved using Ricci flow (see \cite{Hamil,Brendle}). Let us briefly state the definition of this pinching condition, the Differentiable Sphere Theorem, and the generalization thereof that we will exploit.

\begin{definition}[Pointwise $\nicefrac{1}{4}$-Pinched Sectional Curvature Condition]
A Riemannian manifold $(\Sigma,g)$ has \emph{pointwise $\nicefrac{1}{4}$-pinched sectional curvature} if $\Sigma$ has positive sectional curvature, $\mathcal{K}\geq0$, and if $\mathcal{K}(\pi_1) <4\mathcal{K}(\pi_2)$ holds for every pair of $2$-planes $\pi_1,\pi_2 \subset T_p\Sigma$ at every $p\in\Sigma$.
\end{definition}

Note that the fact that the sectional curvature is positive implies positive scalar curvature.

\begin{theorem}[Differentiable Sphere Theorem~\cite{sphere2}]
Let $(\Sigma,g)$ be a compact Riemannian manifold of dimension $n\geq 4$ with pointwise $\nicefrac{1}{4}$-pinched sectional curvature. Then $\Sigma$ admits a metric of constant curvature and therefore is diffeomorphic to a spherical space form.
\end{theorem}
	
\begin{theorem}[Brendle--Schoen \cite{sphere2}]\label{patn}
Let $(\Sigma,g_o)$  be a compact Riemannian manifold of dimension $n\geq 4$ with pointwise $\nicefrac{1}{4}$-pinched sectional curvature. Let $\lbrace h(\tau)\rbrace_{\tau\in [0,T)}$, denote the unique maximal solution to Ricci flow on $\Sigma$ with initial metric $g_o$. Then the rescaled path of metrics $\frac{1}{2 (n-1)(T-\tau)}h(\tau)$ smoothly converges to a metric of constant sectional curvature 1 as $\tau\rightarrow T$.
\end{theorem}

Note that the $\nicefrac{1}{4}$-pinched sectional curvature condition may not be preserved by the flow, but the positivity of the scalar curvature is preserved along the flow. Combining this with Lemma~\ref{lem:normalizing}, one gets the desired path.

\begin{corollary}[Existence of Paths with Desired Properties, $n\geq4$]\label{coro:paths}
Let $n\geq 4$. Given any metric $g_{o}$ on $\mathbb{S}^n$ with pointwise $\nicefrac{1}{4}$-pinched sectional curvature. Then there exists a smooth path of metrics $\{g(t)\}_{t\in[0,1]}$ on $\mathbb{S}^{n}$ satisfying Conditions \hyperref[page:conditions]{(E1-E3)} and $\scal(g(t))>0$ for all $t\in[0,1]$.
\end{corollary}

In \cite{sec}, Cabrera~Pacheco and Miao consider an alternative way of constructing such paths for $n\geq4$, using extrinsic geometric flows developed by Gerhardt~\cite{Gerhardt} and Urbas~\cite{Urbas}.

\subsection{Constructing the Collar Extension}
Using the paths described above, we will now construct, as in \cite{impor}, a collar extension which is an \ECR\ $(M,\gamma,E,\Lambda)$ of the form $M=[0,1] \times \mathbb{S}^n$,  $\gamma$ as in~\eqref{eq:ansatzgamma}, $E$ as in~\eqref{eq:ansatzE}, satisfying the \sDEC~\eqref{eq:sDEC}, inducing a metric isometric to a given metric $g_{o}$ on the \emph{inner boundary} $\lbrace0\rbrace\times\mathbb{S}^{n}$, and such that $\lbrace t\rbrace\times\mathbb{S}^{n}$ has positive mean curvature for $t\in(0,1]$ and vanishing mean curvature for $t=0$. For later convenience, we will slightly deviate from above by using a general compact interval $[a,b]$.

\begin{lemma}[Divergence-Free Electric Field; Outward Unit Normal]\label{lem:diver}
Let $n\geq2$ and let $\{g(t) \}_{t\in[a,b]}$ be a smooth path of metrics on $\mathbb{S}^{n}$ satisfying $\tr_{g(t)}g'(t)=0$ for all $t\in[a,b]$, with $r_{a}$ denoting the volume radius of $g(a)$, i.e., $\omega_{n}r_{a}^{n}\definedas\vert\lbrace a\rbrace\times\mathbb{S}^{n}\vert_{g(a)}$. Set $M\definedas[a,b]\times\mathbb{S}^{n}$ and let $v\colon M \to\mathbb{R}^{+}$ and $F\colon [a,b] \to \mathbb{R}^{+}$  be smooth, positive functions. Let $\gamma$ be the smooth Riemannian metric given by $\gamma = v(t,\cdot)^2 dt^2 +F(t)^2 g(t)$ on $M$, and let $q\in\R$ be a constant. Then the smooth vector field $E$ defined by
$E\definedas \frac{q}{v(t,\cdot) r_a^n F(t)^n} \partial_t$ is divergence-free, i.e., satisfies $\diver_\gamma E =0$ on $M$, and $\mathcal{Q}(\lbrace t\rbrace\times\mathbb{S}^{n})=q$ for all $t\in[a,b]$ with respect to the \emph{outward unit normal} $\nu=\frac{1}{v(t,\cdot)}\partial_{t}$.
\end{lemma}

\begin{proof}
We compute $\diver_{\gamma}E=\frac{q}{vr_{a}^{n}F^{n}}\left(-\frac{1}{vF^{n}}(vF^{n})'+\diver_{\gamma}\partial_{t}\right)=\frac{q}{2vr_{a}^{n}F^{n}}\tr_{g}g'=0$. The claim about $\mathcal{Q}$ follows from the Divergence Theorem as $E$ is divergence-free, and from
\begin{align*}
\mathcal{Q}(\lbrace a\rbrace\times\mathbb{S}^{n})&=\frac{1}{\omega_{n}}\int_{\lbrace a\rbrace\times\mathbb{S}^{n}}\gamma\left(\frac{q}{v(a,\cdot)r_{a}^{n}F(a)^{n}}\partial_{t},\frac{1}{v(a,\cdot)}\partial_{t}\right)dV_{g}
&=\frac{q\vert\lbrace a\rbrace\times\mathbb{S}^{n}\vert_{g}}{r_{a}^{n}F(a)^{n}\omega_{n}}=q.\quad\quad\quad\qedhere
\end{align*}
\end{proof}

Next, let us briefly recall the mean curvature of the coordinate spheres $\lbrace t\rbrace\times\mathbb{S}^{n}$ in $(M,\gamma)$ and the scalar curvature of $(M,\gamma)$ as in Lemma~\ref{lem:diver}.

\begin{lemma}[Mean Curvature of Coordinate Spheres]\label{lem:meancurvature}
Let $n\geq2$ and let $\{g(t) \}_{t\in[a,b]}$ be a smooth path of metrics on $\mathbb{S}^{n}$ satisfying $\tr_{g(t)}g'(t)=0$ for all $t\in[a,b]$. Set $M\definedas[a,b]\times\mathbb{S}^{n}$ and let $v\colon M \to\mathbb{R}^{+}$ and $F\colon [a,b] \to \mathbb{R}^{+}$  be smooth, positive functions. Let $\gamma$ be the smooth Riemannian metric given by $\gamma = v(t,\cdot)^2 dt^2 +F(t)^2 g(t)$
 on $M$. Then the mean curvature $H(\lbrace t\rbrace\times\mathbb{S}^{n})$ of the coordinate sphere $\lbrace t\rbrace\times\mathbb{S}^{n}$, $t\in[a,b]$, is given by
 \begin{align}\label{eq:meancurvature}
 H(\lbrace t\rbrace\times\mathbb{S}^{n})&=\frac{nF'(t)}{v(t,\cdot)F(t)}
 \end{align}
 with respect to the outward unit normal $\nu=\frac{1}{v}\partial_{t}$. In particular, the sign of $H(\lbrace t\rbrace\times\mathbb{S}^{n})$ depends solely on the sign of $F'(t)$, and $\lbrace t\rbrace\times\mathbb{S}^{n}$ is minimal if and only if $F'(t)=0$.
\end{lemma}
\begin{proof}
A straightforward computation gives
\begin{align*}
 H(\lbrace t\rbrace\times\mathbb{S}^{n})&=\frac{nF'(t)}{v(t,\cdot)F(t)}+\frac{1}{2v(t,\cdot)}\tr_{g(t)}g'(t)=\frac{nF'(t)}{v(t,\cdot)F(t)}.\qedhere
\end{align*}
\end{proof}

\begin{lemma}[Scalar Curvature]\label{lem:scalar}
Let $n\geq2$ and let $\{g(t) \}_{t\in[a,b]}$ be a smooth path of metrics on $\mathbb{S}^{n}$ satisfying $\tr_{g(t)}g'(t)=0$ for all $t\in[a,b]$. Set $M\definedas[a,b]\times\mathbb{S}^{n}$ and let $v\colon M \to\mathbb{R}^{+}$ and $F\colon [a,b] \to \mathbb{R}^{+}$  be smooth, positive functions. Let $\gamma$ be the smooth Riemannian metric given by $\gamma = v(t,\cdot)^2 dt^2 +F(t)^2 g(t)$ on $M$. Then the scalar curvature $\scal(\gamma)$ of $(M,\gamma)$ can be written as
\begin{align}\label{eq:scal}
\begin{split}
\scal(\gamma)&= \frac{\scal(g)}{F^{2}}-\frac{2\Delta_{g}v}{vF^{2}}+\frac{1}{v^2}\left(\frac{2nv'F'}{vF}-\frac{n}{F^{2}}\left((n-1)F'^{2}+2FF''\right)-\frac{1}{4}\abs{g'}^2_{g} \right),
\end{split}
\end{align}
where $'$ denotes a (partial) time derivative.
\end{lemma}
	
\begin{proof}
A proof of this formula can be found in \cite[(3.6)]{survey}.
\end{proof}

The following lemma specializes the scalar curvature formula to rotational symmetry and characterizes when rotationally symmetric \ECR s satisfy the (strict) \DEC. It will be useful in Section~\ref{sec:gluing}.
\begin{lemma}[\DEC\ in Rotational Symmetry, Radial Profile]\label{le4.1}
Let $n\geq2$, $M\definedas[a,b]\times\mathbb{S}^{n}$, and $f\colon [a,b] \to \mathbb{R}^{+}$ a smooth, positive function. Let $\gamma$ be the smooth Riemannian metric given by $\gamma = ds^2 +f(s)^2 g_{*}$ on $M$, where $g_{*}$ denotes the canonical metric on $\mathbb{S}^{n}$. Let $q\in\mathbb{R}$ and let $E$ be the smooth vector field given by $E\definedas \frac{q}{f^n} \partial_s$. Then for any $\Lambda\leq0$, the \ECR\ $(M,\gamma,E,\Lambda)$ satisfies the \DEC~\eqref{eq:DEC} if and only if
\begin{align}\label{4.1}
f'' &\leq \frac{n-1}{2f}\left( \frac{h_{q,\Lambda}\circ f}{f^{2(n-1)}}- (f')^2 \right)
\end{align}
holds on $[a,b]$, while the \sDEC~\eqref{eq:sDEC} is satisfied if and only if 
\begin{align}\label{s4.1}
f'' &< \frac{n-1}{2f}\left( \frac{h_{q,\Lambda}\circ f}{f^{2(n-1)}}- (f')^2 \right)
\end{align}
on $[a,b]$. The function $f$ is called the \emph{radial profile of $(M,\gamma,E,\Lambda)$}.
\end{lemma}
\begin{proof}
Applying Lemma~\ref{lem:diver} with $v=1$, $F=f$ and $t=s$, we immediately see that $E$ is defined as it is defined there and is hence divergence-free so that $(M,\gamma,E,\Lambda)$ is an \ECR. Similarly, applying Lemma~\ref{lem:scalar} in the rotationally symmetric case $v=1$, $F=f$, $t=s$, we obtain
\begin{align*}
\scal(\gamma) &= \frac{\scal(g_{*})}{f^{2}}-\frac{n}{f^{2}}\left((n-1)f'^{2}+2ff''\right)=\frac{n}{f^2} [(n-1)(1-f'^2)-2 ff'' ].
\end{align*}
Imposing $\scal(\gamma) \geq (>)\,2\Lambda+\frac{n(n-1)q^2}{f(s)^{2n}}$ leads to the inequality~\eqref{4.1} (respectively to~\eqref{s4.1}).
\end{proof}

Let us now turn to a sufficient condition ensuring the \sDEC~\eqref{eq:sDEC} beyond rotational symmetry. Our ansatzes for $F$ and $v$ go back to and unify those in \cite{pri,hyper,impor}.

\begin{lemma}[Collar Extension]\label{lem:collar}
Let $n\geq2$ and let $\{g(t) \}_{t\in[0,1]}$ be a smooth path of metrics on $\mathbb{S}^{n}$ satisfying Condition~\hyperref[extension:3]{(E3)}. Let $\varepsilon\in(0,1]$ and set $F\colon [0,1] \to \mathbb{R}^{+}\colon t\mapsto\sqrt{1+\varepsilon t^{2}}$. and let $\gamma$ be the smooth Riemannian metric given by $\gamma = v(t,\cdot)^2 dt^2 +F(t)^2 g(t)$ on $M\definedas[0,1]\times\mathbb{S}^{n}$ where $v\colon M\to\R^{+}$ is given by the following formulas under the following conditions
\begin{enumerate}
\item if $n\geq2$ and $\lambda_{1}(-\triangle_{g(t)}+\tfrac{1}{2}\scal(g(t)))>\kappa$ for all $t\in[0,1]$: set $v\definedas Au$, where $u(t,\cdot)$ is the positive first $L^{2}(g(t))$-normalized eigenfunction of $-\triangle_{g(t)}+\tfrac{1}{2}\scal(g(t))$, i.e., normalized such that $\int_{\mathbb{S}^{n}}u(t,\cdot)^{2}\,dV_{g(t)}=\omega_{n} r_{o}^{n}$,
\item if $n=2$ and $K(g(t))>-\kappa$ for all $t\in[0,1]$: set $v\definedas A$,
\item if $n\geq2$ and $\scal(g(t))>2\kappa$ for all $t\in[0,1]$: set $v\definedas A$,
\end{enumerate}
where $\kappa\geq0$, $A>0$ are real constants and $r_{o}$ denotes the volume radius of $\lbrace0\rbrace\times\mathbb{S}^{n}$ with respect to the volume measure of $g(0)$. Now choose $\Lambda\leq0$ and $q\in\R$ satisfying
\begin{align}\tag{C}\label{eq:Qcond}
\begin{split}
\left\{
\begin{array}{lll}
\frac{n(n-1)q^2}{r_{o}^{2n}}&<2(\kappa-\Lambda)&\text{ in Cases~(1) and (3)}\\
\frac{q^2}{r_{o}^{4}}&<-(\kappa+\Lambda)&\text{ in Case~(2)},
\end{array}
\right.
\end{split}
\end{align}
and let $E$ be the electric field defined by~\eqref{eq:ansatzE}. Then there exists a constant $A_{o}>0$ depending only on $n$, $\Lambda$, $\kappa$, $q$, and the path $\{g(t)\}_{t\in[0,1]}$ (and in particular not on $\varepsilon$), such that for every $A\geq A_{o}$, the \ECR\ $(M,\gamma,E,\Lambda)$ satisfies the \sDEC~\eqref{eq:sDEC}. Moreover, the inner boundary $\lbrace0\rbrace\times\mathbb{S}^{n}$ is minimal and the spheres $\{t \} \times \mathbb{S}^n $ for $t\in (0,1]$ are strictly mean convex with respect to the outward unit normal $\nu=\frac{1}{v}\partial_t$.
\end{lemma}
\begin{proof}
We will treat Cases~(1) through (3) separately, using~\eqref{eq:scal} and the specific definition of $v$. In any case, we will use $1\leq F(t)\leq\sqrt{2}$ and $F'(t)=\frac{\varepsilon t}{F(t)}$, $F''(t)=\frac{\varepsilon}{F(t)}-\frac{\varepsilon^{2}t^{2}}{F(t)^{3}}=\frac{\varepsilon-F'(t)^{2}}{F(t)}$ for all $t\in[0,1]$. For Case~(1), we compute
\begin{align*}
&\scal(\gamma)-2\Lambda-n(n-1)\vert E\vert^{2}_{\gamma}\\
&\quad\quad=\frac{2}{vF^{2}}\left(-\triangle_{g}v+\frac{1}{2}\scal(g)v\right)+\frac{1}{v^2}\left(\frac{2nv'F'}{vF}-\frac{n}{F^{2}}\left((n-1)F'^{2}+2FF''\right)-\frac{1}{4}\abs{g'}^2_{g} \right)\\
&\quad\quad\quad-2\Lambda-n(n-1)\vert E\vert_{\gamma}^{2}\\
&\quad\quad=\frac{2}{uF^{2}}\left(-\triangle_{g}u+\frac{1}{2}\scal(g)u\right)+\frac{1}{A^2u^{2}}\left(\frac{2nu'F'}{uF}-\frac{n}{F^{2}}\left((n-3)F'^{2}+2\varepsilon\right)-\frac{1}{4}\abs{g'}^2_{g} \right)\\
&\quad\quad\quad-2\Lambda-\frac{n(n-1)q^2}{r_{o}^{2n}F^{2n}}\\
&\quad\quad\geq \frac{2}{F^{2}}\lambda_{1}(-\triangle_{g}+\tfrac{1}{2}\scal(g))-2\Lambda-\frac{n(n-1)q^2}{r_{o}^{2n}F^{2n}}\\
&\quad\quad\quad+\frac{1}{A^2}\frac{\min\lbrace 0,\min_{M}\left(\frac{2nu'}{u}\right)\rbrace-c_{n}-\max_{M}(\frac{1}{4}\abs{g'}^2_{g})}{\max_{M}(u^{2})}\\
&\quad\quad =\frac{1}{F^{2}} \underbrace{\left(2\lambda_{1}(-\triangle_{g}+\tfrac{1}{2}\scal(g)) -2 F^{2}\Lambda - \frac{n(n-1)q^2}{r_{o}^{2n}F^{2(n-1)}}\right)}_{=: (I)}\\
&\quad\quad\quad+\frac{1}{A^2}\underbrace{\frac{\min\lbrace 0,\min_{M}\left(\frac{2nu'}{u}\right)\rbrace-c_{n}-\max_{M}(\frac{1}{4}\abs{g'}^2_{g})}{\max_{M}(u^{2})}}_{<0},
\end{align*}
with $c_{2}\definedas 4$ and $c_{n}\definedas n(n-1)$ for $n\geq3$. We estimate $(I)$ from below by $(I)>2(\kappa -\Lambda)  -\frac{n(n-1)q^2}{r_{o}^{2n}}>0 $, using $\lambda_{1}(-\triangle_{g}+\tfrac{1}{2}\scal(g))>\kappa\geq0$, $F\geq1$, and Condition~\eqref{eq:Qcond}. Recalling the choice of $u$, the second bracket only depends on $n$ and on the path $\{g(t)\}_{t\in[0,1]}$, and hence we can choose $A_{o}>0$ depending only on $n$, $\kappa$, $\Lambda$, $q$, and the path $\{g(t)\}_{t\in[0,1]}$ such that the \sDEC\ holds for all $A\geq A_{o}$ and all $\varepsilon\in(0,1]$.

Similarly, for Cases~(2) and (3), we compute
\begin{align*}
&\scal(\gamma)-2\Lambda-n(n-1)\vert E\vert^{2}_{\gamma}\\
&\quad\quad=\frac{\scal(g)}{F^{2}}+\frac{1}{A^2}\left(-\frac{n}{F^{2}}\left((n-3)F'^{2}+2\varepsilon\right)-\frac{1}{4}\abs{g'}^2_{g} \right)-2\Lambda-\frac{n(n-1)q^2}{r_{o}^{2n}F^{2n}}.\\
&\quad\quad\geq\frac{\scal(g)}{F^{2}}-2\Lambda-\frac{n(n-1)q^2}{r_{o}^{2n}F^{2n}}+\frac{1}{A^2}\left(-c_{n}-\max_{M}(\frac{1}{4}\abs{g'}^2_{g}) \right)\\
&\quad\quad =\frac{1}{F^{2}} \underbrace{\left(\scal(g) -2 F^{2}\Lambda - \frac{n(n-1)q^2}{r_{o}^{2n}F^{2(n-1)}}\right)}_{=: (I)} +\frac{1}{A^2}\underbrace{\left(-c_{n}-\max_{M}(\frac{1}{4}\abs{g'}^2_{g}) \right)}_{<0}
\end{align*}
Now, in Case~(2),  we can estimate $(I)$ from below by $(I)>-2\left(\kappa+\Lambda+\tfrac{q^2}{r_{o}^{4}}\right)>0$, using $\scal(g)=2K(g)>-2\kappa$, $F\geq1$,  and Condition~\eqref{eq:Qcond}. In Case~(3), we estimate $(I)$ from below by $(I)>2(\kappa - \Lambda) -\frac{n(n-1)q^2}{r_{o}^{2n}}>0 $, using $\scal(g)>2\kappa\geq0$, $F\geq1$, and Condition~\eqref{eq:Qcond}. In both cases, we can conclude as before. The claim that $(M,\gamma,E,\Lambda)$ constitutes an \ECR\ is a direct application of Lemma~\ref{lem:diver} while the claims about the mean curvatures of the coordinate spheres follow from Lemma~\ref{lem:meancurvature}.
\end{proof}

\begin{remark}
Note that, in Case~(2), Condition~\eqref{eq:Qcond} cannot be fulfilled for arbitrary choices of $\kappa\geq0$ and $\Lambda\leq0$; even the uncharged case $q=0$ is ruled out if $\kappa+\Lambda\geq0$.
\end{remark}

\begin{remark}\label{rem:normalized}
Note that the normalization $\fint_{\mathbb{S}^{n}}u^{2}dV_{g}=1$ of the (positive) first eigenfunction $u$ of $-\triangle_{g}+\tfrac{1}{2}\scal(g)$ used in Lemma~\ref{lem:collar} and the rest of this paper differs from the one used in previous works. It ensures that $u=1$ if $g$ is round, making Case~(1) of Lemma~\ref{lem:collar} more similar to Cases~(2)--(4).
\end{remark}

The reason why Cases~(1) and (3) and Case~(2) differ in this regard and in their treatment in the proof is that we impose a positive lower bound on $\lambda_{1}(-\triangle_{g(t)}+\tfrac{1}{2}\scal(g(t)))$ and $\scal(g(t))$ along the path $\{g(t)\}_{t\in[0,1]}$ in Cases~(1) and (3), respectively, while we impose a negative lower bound on $K(g(t))$ along the path $\{g(t)\}_{t\in[0,1]}$ in Case~(2). This is related to the results on existence of such paths described above and will become more clear in Section~\ref{sec:main}.

Before we move on, let us collect a few further properties of the collar extensions constructed in Lemma~\ref{lem:collar}.
 
\begin{lemma}[A Priori Sub-Extremality of Collar Extensions]\label{lem:subapriori}
Let $n\geq2$, $\kappa\geq0$, and let $\{g(t) \}_{t\in[0,1]}$ be a smooth path of metrics on $\mathbb{S}^{n}$ as in any of the Cases~(1), (2), or (3) of Lemma~\ref{lem:collar}, and assume in addition that $g(1)$ is round. Let $r_{o}$ denote the volume radius of $g(0)$ and assume that $r_{o}$, $q$, $\Lambda$, and $\kappa$ satisfy Condition~\eqref{eq:Qcond} of Lemma~\ref{lem:collar}. Then the \emph{quasi-local sub-extremality condition} $h_{q,\Lambda}(r_{o})>0$ holds. Moreover, $\kappa < \frac{n(n-1)}{2 r_{o}^{2}}$ holds in Cases~(1) and~(3) of Lemma~\ref{lem:collar}.
\end{lemma}
\begin{proof}
Let us first show that $\kappa < \frac{n(n-1)}{2r_{o}^{2}}$ in Cases~(1) and~(3): To see this, we choose $t=1$ and use that $\scal(g(1))= \frac{n(n-1)}{r_{o}^2}$ as $g(1)$ is a round sphere of volume radius $r_{o}$. Thus, in Case~(3), the condition $\scal(g(1))> 2\kappa$ readily implies $\kappa < \frac{n(n-1)}{2 r_{o}^{2}}$. In Case~(1), the fact that the first eigenfunction of $-\triangle_{g(1)}+\tfrac{1}{2}\scal(g(1))$ is given by $u(1,\cdot)=1$ (see Remark~\ref{rem:normalized}) shows that $\tfrac{1}{2}\scal(g(1))=\lambda_{1}(-\triangle_{g(1)}+\tfrac{1}{2}\scal(g(1))>\kappa$ implies that $\kappa < \frac{n(n-1)}{2 r_{o}^{2}}$. Combining this with the first condition in \eqref{eq:Qcond}, we hence find $h_{q,\Lambda}(r_{o})=r_{o}^{2(n-1)}-q^{2}-\frac{2\Lambda r_{o}^{2n}}{n(n-1)}>r_{o}^{2(n-1)}\left(1-\frac{2\kappa r_{0}^{2}}{n(n-1)}\right)>0$ in Cases~(1) and~(3). 

In Case (2), the second condition in~\eqref{eq:Qcond} and $\kappa\geq0$ lead to $h_{q,\Lambda}(r_{o})=r_{o}^{2}-q^{2}-\Lambda r_{o}^{4}>r_{o}^{2}+\kappa r_{o}^{4}  >0$.
\end{proof}
 
The following lemma asserts that the generalized Hawking Mass at the far end of the collar, $\lbrace1\rbrace\times\mathbb{S}^{n}$, is larger than that at the inner boundary for any collar extension constructed via Lemma~\ref{lem:collar}. For $n=2$, this is of course a consequence of the strict monotonicity of the generalized Hawking mass along such collar extensions, see Section~\ref{sec:Hawkingcollar}. However, for $n\geq3$, the monotonicity properties we will prove for the generalized Hawking mass along collar extensions in Proposition~\ref{prop:monotonicity} do not seem to apply here; but see Proposition~\ref{prop:monotonicity2}. The lemma also establishes an upper bound for the generalized Hawking mass at the far end of such collars; both will be useful later.

\begin{lemma}[Generalized Hawking Mass on Far End of Collar]\label{lem:Hawking*}
Let $n\geq2$ and let $g_{o}$, $r_{o}$, $q$, $\Lambda$, $\kappa$, $A\geq A_{o}$, and $\varepsilon\in(0,1]$ be as in Lemma~\ref{lem:collar}. Set $m_o\definedas\tfrac{r_o^{n-1}}{2}\left(1+\tfrac{q^2}{r_o^{2(n-1)}} -\tfrac{2\Lambda r_o^{2}}{n(n+1)}\right)$ and consider the collar extension constructed in Lemma~\ref{lem:collar}, denoted as
\begin{align*}
\mathcal{S}_{\varepsilon,A}=([0,1]\times\mathbb{S}^{n},\gamma_{\varepsilon,A}=(Au(t,\cdot))^{2}dt^{2}+(1+\varepsilon t^{2})g(t),E_{\varepsilon,A}=\frac{q}{Au(t,\cdot)r_{o}^{n}(1+\varepsilon t^{2})^{\frac{n}{2}}}\partial_{t},\Lambda).
\end{align*}
Assume that $g(1)$ is round. Then, for all $t\in[0,1]$, $\varepsilon\in(0,1]$, and $A\geq A_{o}$, we have
\begin{align}\label{eq:m0A}
\mathcal{M}_{\mathcal{S}_{\varepsilon,A}}(\lbrace0\rbrace\times\mathbb{S}^{n})&=m_{o},\\\label{eq:qA}
\mathcal{Q}_{\mathcal{S}_{\varepsilon,A}}(\lbrace t\rbrace\times\mathbb{S}^{n})&=q,
\end{align}
where $\mathcal{M}_{\mathcal{S}_{\varepsilon,A}}$ denotes the generalized Hawking mass functional and $\mathcal{Q}_{\mathcal{S}_{\varepsilon,A}}$ denotes the quasi-local charge functional in $\mathcal{S}_{\varepsilon,A}$. Moreover, there exists $0<\varepsilon_{o}\leq1$ depending only on $n$, $r_{o}$, $q$, $\Lambda$, and $A_{o}$ such that 
\begin{align}\label{eq:grows}
 \mathcal{M}_{\mathcal{S}_{\varepsilon,A}}(\lbrace 1\rbrace\times\mathbb{S}^n)>m_{o}
\end{align}
for all $\varepsilon\in(0,\varepsilon_{o}]$ and all $A\geq A_{o}$. Finally, we have the upper bound
\begin{align}
\mathcal{M}_{\mathcal{S}_{\varepsilon,A}}(\lbrace 1\rbrace\times\mathbb{S}^n)<(1+\varepsilon)^{\frac{n+1}{2}}m_{o}
\end{align}
for all $\varepsilon\in(0,1]$ and all $A\geq A_{o}$.
\end{lemma}
\begin{proof}
First of all, \eqref{eq:qA} follows from Lemma~\ref{lem:diver} for all $t\in[0,1]$, $\varepsilon\in(0,1]$, $A\geq A_{o}$. As a consequence of~\eqref{eq:qA} and the fact that $\partial M=\lbrace0\rbrace\times\mathbb{S}^{n}$ is minimal in $\mathcal{S}_{\varepsilon,A}$, Lemma~\ref{lem:meancurvature} gives us \eqref{eq:m0A} for all $\varepsilon\in(0,1]$, $A\geq A_{o}$. For~\eqref{eq:grows}, we use Proposition~\ref{prop:Hawkingre}, Lemma~\ref{lem:h}, Lemma~\ref{lem:subapriori}, and Taylor's Formula to find
\begin{align*}
&\mathcal{M}_{\mathcal{S}_{\varepsilon,A}}(\lbrace 1\rbrace\times\mathbb{S}^n)=\frac{(\sqrt{1+\varepsilon})^{n-1}r_{o}^{n-1}}{2} \left(p_{0,q,\Lambda}(\sqrt{1+\varepsilon}\,r_{o})-\frac{\varepsilon^{2} r_{o}^{2} }{(1+\varepsilon)A^{2}}\right)\\
&\quad\geq \frac{r_{o}^{n-1}}{2} \left((\sqrt{1+\varepsilon})^{n-1} p_{0,q,\Lambda}(\sqrt{1+\varepsilon}\,r_{o})-\frac{\varepsilon^{2} (1+\varepsilon)^\frac{n-3}{2} r_{o}^{2}}{A_{o}^{2}}\right)\\
&\quad=\frac{r_{o}^{n-1}}{2} \left((\sqrt{1+\varepsilon})^{n-1}\left[ p_{m_{o},q,\Lambda}(\sqrt{1+\varepsilon}\,r_{o})+\frac{2m_{o}}{(\sqrt{1+\varepsilon}\,r_{o})^{n-1}}\right]-\frac{\varepsilon^{2} (1+\varepsilon)^\frac{n-3}{2} r_{o}^{2}}{A_{o}^{2}}\right)\\
&\quad=\frac{r_{o}^{n-1}}{4}(2+(n-1)\varepsilon)(\underbrace{p_{m_{o},q,\Lambda}(r_{o})}_{=0}+\frac{1}{2}\underbrace{p_{m_{o},q,\Lambda}'(r_{o})}_{=\tfrac{(n-1)h_{q,\Lambda}(r_{o})}{r_{o}^{2n-1}}}r_{o}\varepsilon)+m_{o}+\mathcal{O}(\varepsilon^{2})\\
&\quad=m_{o}+\frac{(n-1)h_{q,\Lambda}(r_{o})}{4r_{o}^{n-1}}\varepsilon+\mathcal{O}(\varepsilon^{2})>m_{o}
\end{align*}
as $\varepsilon\searrow0$ for all $A\geq A_{o}$. The constant in the $\mathcal{O}$-notation depends only on $n$, $r_{o}$, $q$, $\Lambda$, and $A_{o}$. In other words, there exists $0<\varepsilon_{o}\leq1$ depending only on $n$, $r_{o}$, $q$, $\Lambda$, and $A_{o}$ such that $\mathcal{M}_{\mathcal{S}_{\varepsilon,A}}(\lbrace 1\rbrace\times\mathbb{S}^n)>m_{o}$ holds for all $0<\varepsilon\leq\varepsilon_{o}$ and all $A\geq A_{o}$.

To prove the upper estimate, we use again Proposition~\ref{prop:Hawkingre} to compute
\begin{align*}
\mathcal{M}_{\mathcal{S}_{\varepsilon,A}}(\lbrace 1\rbrace\times\mathbb{S}^n)&\leq\frac{(1+\varepsilon)^\frac{n-1}{2}r_{o}^{n-1}}{2} \,p_{0,q,\Lambda}(\sqrt{1+\varepsilon}\,r_{o})\\
&=\frac{(1+\varepsilon)^{\frac{n+1}{2}}r_o^{n-1}}{2}\left(\frac{1}{1+\varepsilon}+\frac{q^2}{r_o^{2(n-1)}(1+\varepsilon)^{n-2}}-\frac{2\Lambda}{n(n+1)} r_o^2 \right) \\
&<\frac{(1+\varepsilon)^{\frac{n+1}{2}}r_o^{n-1}}{2}\left(1+\frac{q^2}{r_o^{2(n-1)}}-\frac{2\Lambda}{n(n+1)} r_o^2 \right)\\
&=\frac{(1+\varepsilon)^{\frac{n+1}{2}}r_o^{n-1}}{2}\,p_{0,q,\Lambda}(r_{o})=(1+\varepsilon)^{\frac{n+1}{2}}m_{o}.\qedhere
\end{align*}
\end{proof}

\subsection{The Generalized Hawking Mass along Collar Extensions}\label{sec:Hawkingcollar}
The generalized Hawking mass naturally simplifies when studied in rotationally symmetric \ECR s, a case that will be useful for our analysis.

\begin{proposition}[Generalized Hawking Mass in Rotational Symmetry]\label{prop:sphereHawking}
Let $n\geq2$, let $q,\Lambda\in\mathbb{R}$ be constants with $\Lambda\leq0$, and let $f\colon I\to\R^{+}$ be a smooth, positive function defined on an interval $I$. Consider the rotationally symmetric \ECR\ given by $(I \times \mathbb{S}^n , \gamma=ds^2 +f(s)^2 g_* , E=\frac{q}{f(s)^{n}}\partial_{s}, \Lambda)$, where $g_{*}$ denotes the canonical metric on $\mathbb{S}^{n}$. Then every coordinate sphere $\{s\} \times \mathbb{S}^n$, $s\in I$, has generalized Hawking mass
\begin{align}\label{eq:Mp}
\mathcal{M}(\{s\} \times \mathbb{S}^n) &=m+\frac{f(s)^{n-1}}{2}  \left(p_{m,q,\Lambda}(f(s))-f'(s)^2 \right)
\end{align}
for every $m\in\mathbb{R}$.
\end{proposition}
\begin{proof}
This follows from a straightforward computation using Proposition~\ref{prop:Hawkingre} and exploiting the more general formula~\eqref{eq:meancurvature} which shows that the mean curvature $H(s)$ of $\Sigma_{s}$ is given by $H(s)=\frac{n f'(s)}{f(s)}$. More specifically, using~\eqref{eq:Hawkingsimple}, we find
\begin{align*}
\mathcal{M}(\{s\} \times \mathbb{S}^n) &=\frac{f(s)^{n-1}}{2}  \left(p_{0,q,\Lambda}(f(s))- f'(s)^2\right)\\
&=\frac{f(s)^{n-1}}{2}  \left(p_{m,q,\Lambda}(f(s))+\frac{2m}{f(s)^{n-1}}-f'(s)^{2}\right)\\
&=m+\frac{f(s)^{n-1}}{2}  \left(p_{m,q,\Lambda}(f(s))-f'(s)^{2}\right)
\end{align*}
for all $m\in\mathbb{R}$ which implies the claim.
\end{proof}

\begin{remark}[Generalized Hawking Mass for $\Lambda$-Reissner--Nordstr\"om Coordinate Spheres]
For $(n+1)$-dimensional sub-extremal \RN s of mass $m$ and charge $q$, it follows from Lemma~\ref{lem:diver}, Proposition~\ref{prop:sphereHawking}, and Proposition~\ref{prop:extendRN} that $\mathcal{M}(\lbrace s\rbrace\times\mathbb{S}^{n})=m$ and $\mathcal{Q}(\lbrace s\rbrace\times\mathbb{S}^{n})=q$ for all $s\geq 0$. These identities also hold for extremal and super-extremal \RN s as can be seen by a direct computation using that the mean curvature $H(r)$ of $\lbrace r\rbrace\times\mathbb{S}^{n}$, $r> r_{+}$ is given by $H(r)=\tfrac{n\,\sqrt{p_{m,q,\Lambda}(r)}}{r}$, a consequence of~\eqref{eq:meancurvature}.
\end{remark}

We will now discuss monotonicity of the generalized Hawking mass along collar extensions. Recall that the Hawking mass is known to be increasing along Inverse Mean Curvature Flow and its weak version devised by Huisken and Ilmanen~\cite{RPI} in $n+1=3$-dimensional Riemannian manifolds $(M,\gamma)$ satisfying the (uncharged) \DEC\ with $\Lambda=0$, and strictly increasing provided $(M,\gamma)$ satisfies the (uncharged)  \sDEC\ with $\Lambda=0$; this was originally observed by Geroch~\cite{Geroch} in the smooth case. Monotonicity of the asymptotically hyperbolic Hawking mass --- see Remark~\ref{rem:consistent} --- under (weak) Inverse Mean Curvature Flow has been noted to follow from the same considerations by Neves~\cite{Neves} in $n+1=3$ dimensions, assuming the (uncharged) \DEC\ with $\Lambda<0$, with strict monotonicity provided that it holds strictly. Similarly, monotonicity of the charged Hawking mass --- see Remark~\ref{rem:consistent} --- under Inverse Mean Curvature Flow has been established by Jang~\cite{Jang} and extended to weak Inverse Mean Curvature Flow by Disconzi and Khuri~\cite{Discon} in $n+1=3$ dimensions, assuming the charged \DEC\ with $\Lambda=0$, with strict monotonicity provided that it holds strictly.

As all these montonicity proofs require the Gauss--Bonnet Theorem, they do not naturally lend themselves to be adjusted to a higher dimensional scenario. 

For the purposes of this paper, it will be very useful to have the following monotonicity lemma for the $n=2$-generalized Hawking mass at hand. In this case, it is a simplified combination of the known monotonicity properties discussed above. For $n\geq3$, we can only expect and show monotonicity in a slight generalization of the rotationally symmetric scenario, allowing to circumvent the lack of a Gauss--Bonnet Theorem by adding a suitable assumption on the integral of the scalar curvature.

\begin{proposition}[Generalized Hawking Mass Monotonicity]\label{prop:monotonicity}
Let $n\geq2$ and let $\{g(t) \}_{t\in[a,b]}$ be a smooth path of metrics on $\mathbb{S}^{n}$ satisfying $\tr_{g(t)}g'(t)=0$ for all $t\in[a,b]$, with $r_{a}$ denoting the volume radius of $g(a)$, i.e., $\omega_{n} r_{a}^{2}\definedas\vert\lbrace a\rbrace\times\mathbb{S}^{n}\vert_{g(a)}$. Set $M\definedas[a,b]\times\mathbb{S}^{n}$ and let $v\colon M \to\mathbb{R}^{+}$ and $F\colon [a,b] \to \mathbb{R}^{+}$  be smooth, positive functions with $F'\geq0$. Let $\gamma$ be the smooth Riemannian metric given by $\gamma = v(t,\cdot)^2 dt^2 +F(t)^2 g(t)$ on $M$. Set $E\definedas \frac{q}{v(t,\cdot)r_a^n F(t)^n} \partial_t$ for some constant $q\in\R$, and let $\Lambda\leq0$. Assume that the \ECR\ $(M,\gamma,E,\Lambda)$ satisfies the \DEC~\eqref{eq:DEC}. Then for $n=2$, the generalized Hawking mass functional $\mathcal{M}$ in $(M,\gamma,E,\Lambda)$ satisfies
\begin{align}\label{eq:monotone}
\frac{d}{dt}\mathcal{M}(\lbrace t\rbrace\times\mathbb{S}^{2})&\geq0
\end{align}
for all $t\in[a,b]$. For $n\geq3$, the same conclusion follows provided that 
\begin{align}\label{eq:R}
\int_{\lbrace t\rbrace\times\mathbb{S}^{n}}\scal(g(t))dV_{g(t)}&\leq n(n-1)\omega_{n}r_{a}^{n-2}
\end{align}
holds for all $t\in[a,b]$. Furthermore, for all $n\geq2$, if $F'(t_{o})=0$ for some $t_{o}\in[a,b]$ then $\tfrac{d}{dt}\mathcal{M}(\lbrace t\rbrace\times\mathbb{S}^{n})\vert_{t_{o}}=0$. Conversely, if $F'>0$ on $[a,b]$ and \eqref{eq:R} holds if $n\geq3$, then $\tfrac{d}{dt}\mathcal{M}(\lbrace t\rbrace\times\mathbb{S}^{n})=0$ holds on $[a,b]$ if and only if $(M,\gamma,E,\Lambda)$ is isometric to a suitable piece of an $(n+1)$-dimensional \ECR\ of the form $(I\times\mathbb{S}^{n},\gamma=\frac{1}{p_{m,q,\Lambda}(r)}dr^{2}+r^{2}h_{*},E_{m,q,\Lambda},\Lambda)$, where $I\subset(0,\infty)$ is some interval with $p_{m,q,\Lambda}>0$ on $I$, $m\definedas\mathcal{M}(\lbrace a\rbrace\times\mathbb{S}^{n})$, and $h_{*}$ is a smooth metric of volume radius $1$ and constant scalar curvature $\scal(h_{*})=n(n-1)$ on $\mathbb{S}^{n}$. 

Moreover, for all $n\geq2$, inequality~\eqref{eq:monotone} holds strictly provided that the \sDEC~\eqref{eq:sDEC} holds and that $F$ satisfies $F'>0$ on $[a,b]$.
\end{proposition}
\begin{proof}
Let us first compute $\mathcal{M}(\lbrace t\rbrace\times\mathbb{S}^{n})$ for an arbitrary $t\in[a,b]$ and $n\geq2$, with $dV_{F(t)^{2}g(t)}$ and $dV_{g(t)}$ denoting the area forms induced on $\lbrace t\rbrace\times\mathbb{S}^{n}$ by $F(t)^{2}g(t)$ and by $g(t)$, respectively. It is a simple observation that $\vert\lbrace t\rbrace\times\mathbb{S}^n\vert_{F(t)^{2}g(t)}=\omega_{n} r_{a}^{n}F(t)^{n}$. By Lemma~\ref{lem:diver}, $E$ is divergence free and $\mathcal{Q}(\lbrace t\rbrace\times\mathbb{S}^{n})= q$ 
for all $t\in[a,b]$. Moreover, from Lemma~\ref{lem:meancurvature}, we know that $H(\lbrace t\rbrace\times\mathbb{S}^{n})=\frac{nF'(t)}{v(t,\cdot)F(t)}$ for all $t\in[a,b]$, where $ H(\lbrace t\rbrace\times\mathbb{S}^{n})$ denotes the mean curvature of $\lbrace t\rbrace\times\mathbb{S}^{n}$ with respect to the outward unit normal $\nu$. Using~\eqref{eq:Hawkingsimple} we can thus compute
\begin{align*}
&\frac{2\mathcal{M}(\lbrace t\rbrace\times\mathbb{S}^{n})}{r_{a}^{n-1}} = F(t)^{n-1} \left( p_{0,q,\Lambda}(r_aF(t))- \frac{1}{n^{2}\omega_{n}r_{a}^{n-2}F(t)^{n-2}} \int_{\lbrace t\rbrace\times\mathbb{S}^{n}}\!\!\!H(\lbrace t\rbrace\times\mathbb{S}^{n})^2\,dV_{F(t)^{2}g(t)} \right)\\
&=F(t)^{n-1} \left( p_{0,q,\Lambda}(r_aF(t))- \frac{F'(t)^{2}}{\omega_{n}r_{a}^{n-2}} \int_{\lbrace t\rbrace\times\mathbb{S}^{n}}\frac{1}{v(t,\cdot)^{2}}\,dV_{g(t)} \right)
\end{align*}
for all $t\in[a,b]$. Hence $t\mapsto\mathcal{M}(\lbrace t\rbrace\times\mathbb{S}^{n})$ is smooth and we have
\begin{align*}
&\frac{2}{r_{a}^{n-1}}\frac{d}{dt}\mathcal{M}(\lbrace t\rbrace\times\mathbb{S}^{2})\\
&= (n-1)F(t)^{n-2}F'(t)\left( p_{0,q,\Lambda}(r_aF(t))- \frac{F'(t)^{2}}{\omega_{n}r_{a}^{n-2}} \int_{\lbrace t\rbrace\times\mathbb{S}^{n}}\frac{1}{v(t,\cdot)^{2}}\,dV_{g(t)} \right)\\
&\quad\,+F(t)^{n-1}\left( p_{0,q,\Lambda}'(r_aF(t))r_aF'(t)- \frac{2F'(t)F''(t)}{\omega_{n}r_{a}^{n-2}} \int_{\lbrace t\rbrace\times\mathbb{S}^{n}}\frac{1}{v(t,\cdot)^{2}}\,dV_{g(t)}\right.\\
&\phantom{F(t)}\quad +\left.\frac{2F'(t)^{2}}{\omega_{n}r_{a}^{n-2}}\int_{\lbrace t\rbrace\times\mathbb{S}^{n}}\frac{v'(t,\cdot)}{v(t,\cdot)^{3}}\,dV_{g(t)}\right)
\end{align*}
because $\tfrac{d}{dt}dV_{g(t)}=0$ follows from our assumption $\tr_{g(t)}g'(t)=0$, see Remark~\ref{rem:vol0}. Thus 
\begin{align*}
\frac{2}{r_{a}^{n-1}}\frac{d}{dt}\mathcal{M}(\lbrace t\rbrace\times\mathbb{S}^{n})=\, &  (n-1)F(t)^{n-2}F'(t)\left(1-\frac{q^{2}}{r_{a}^{2(n-1)}F(t)^{2(n-1)}}-\frac{2\Lambda r_{a}^{2}F(t)^{2}}{n(n-1)}\right)\\
& +\frac{2F(t)^{n-1}F'(t)^{2}}{\omega_{n} r_{a}^{n-2}}\int_{\lbrace t\rbrace\times\mathbb{S}^{n}}\frac{v'(t,\cdot)}{v(t,\cdot)^{3}}\,dV_{g(t)}\\
& -\frac{F(t)^{n-2}F'(t)}{\omega_{n} r_{a}^{n-2}}\int_{\lbrace t\rbrace\times\mathbb{S}^{n}}\frac{(n-1)F'(t)^{2}+2F(t)F''(t)}{v(t,\cdot)^{2}}\,dV_{g(t)}
\end{align*}
holds for all $t\in[a,b]$. By definition of $\gamma$ and $E$ we know $\vert E\vert_{\gamma}^{2}=\frac{q^{2}}{r_{a}^{2n}F(t)^{2n}}$, so that the \DEC~\eqref{eq:DEC}, its evaluation on the collar via~\eqref{eq:scal}, and $F'\geq0$ give
\begin{align*}
&\frac{2}{r_{a}^{n-1}}\frac{d}{dt}\mathcal{M}(\lbrace t\rbrace\times\mathbb{S}^{n})\\
&\,\geq (n-1)F(t)^{n-2}F'(t)\left(1\;\cancel{-\frac{q^{2}}{r_{a}^{2(n-1)}F(t)^{2(n-1)}}-\frac{2\Lambda r_{a}^{2}F(t)^{2}}{n(n-1)}}\right)\\
&\quad\,+\cancel{\frac{2F(t)^{n-1}F'(t)^{2}}{\omega_{n} r_{a}^{n-2}}\int_{\lbrace t\rbrace\times\mathbb{S}^{n}}\frac{v'(t,\cdot)}{v(t,\cdot)^{3}}\,dV_{g(t)}}\\
&\quad\, + \frac{F(t)^{n}F'(t)}{n\omega_{n} r_{a}^{n-2}}\int_{\lbrace t\rbrace\times\mathbb{S}^{n}}\left[-\frac{\scal(g(t))}{F(t)^{2}}+\frac{2\triangle_{g(t)}v(t,\cdot)}{v(t,\cdot)F(t)^{2}}\cancel{-\frac{2nv'(t,\cdot)F'(t)}{v(t,\cdot)^{3}F(t)}}+\frac{\vert g'(t)\vert_{g(t)}^{2}}{4v(t,\cdot)^{2}}\right]dV_{g(t)}\\
&\quad\,  +\frac{F(t)^{n}F'(t)}{n\omega_{n} r_{a}^{n-2}}\cancel{\left[n(n-1)\frac{q^{2}}{r_{a}^{2n}F(t)^{2n}}+2\Lambda\right]}\int_{\lbrace t\rbrace\times\mathbb{S}^{n}} dV_{g(t)}\\
&\,= (n-1)F(t)^{n-2}F'(t)\left(1-\int_{\lbrace t\rbrace\times\mathbb{S}^{n}}\frac{ \scal(g(t))}{n(n-1)\omega_{n}r_{a}^{n-2}}\,dV_{g(t)}\right)\\
&\quad\,+\frac{2F(t)^{n-2}F'(t)}{n\omega_{n} r_{a}^{n-2}}\int_{\lbrace t\rbrace\times\mathbb{S}^{n}} \frac{\triangle_{g(t)}v(t,\cdot)}{v(t,\cdot)}\, dV_{g(t)}+\frac{F(t)^{n}F'(t)}{4n\omega_{n} r_{a}^{n-2}}\int_{\lbrace t\rbrace\times\mathbb{S}^{n}} \frac{\vert g'(t)\vert_{g(t)}^{2}}{v(t,\cdot)^{2}}\,dV_{g(t)}\\
&\,\geq \frac{2F(t)^{n-2}F'(t)}{n\omega_{n} r_{a}^{n-2}}\int_{\lbrace t\rbrace\times\mathbb{S}^{n}} \frac{\vert \nabla_{g(t)}v(t,\cdot)\vert_{g(t)}^{2}}{v(t,\cdot)^{2}} \,dV_{g(t)}\geq0,
\end{align*}
where $\nabla_{g(t)}v(t,\cdot)$ denotes the (partial) $\mathbb{S}^{n}$-gradient of $v(t,\cdot)$ with respect to the metric $g(t)$. For the last two steps, we have used the Gauss--Bonnet Theorem to see that the term in brackets vanishes for $n=2$ and Assumption~\eqref{eq:R} for $n\geq3$, the Divergence Theorem together with closedness of $\lbrace t\rbrace\times\mathbb{S}^{n}$, and $F'\geq0$ on $[a,b]$. This shows the desired monotonicity of the generalized Hawking mass. In particular, if $F'(t_{o})=0$ for some $t_{o}\in[a,b]$, we have $\frac{d}{dt}\mathcal{M}(\lbrace t\rbrace\times\mathbb{S}^{n})\vert_{t_{o}}=0$ as desired. Conversely, if $\frac{d}{dt}\mathcal{M}(\lbrace t\rbrace\times\mathbb{S}^{n})\vert_{t_{o}}=0$ and $F'(t_{o})>0$ then equality holds in the \DEC~\eqref{eq:DEC} on $\lbrace t\rbrace\times\mathbb{S}^{n}$, $v(t_{o},\cdot)=\text{const}$ on $\mathbb{S}^{n}$, and $g'(t_{o})=0$, and for $n\geq3$ equality must hold in Assumption~\eqref{eq:R} at $t_{o}$. In particular, if $F'>0$ on $[a,b]$ then $\tfrac{d}{dt}\mathcal{M}(\lbrace t\rbrace\times\mathbb{S}^{n})=0$ holds if and only if the \DEC~\eqref{eq:DEC} is saturated on $M$, $v(t,\cdot)=v(t)$ on $\mathbb{S}^{n}$ for all $t\in [a,b]$, and $g'=0$ on $M$, plus equality holds in Assumption~\eqref{eq:R} for all $t\in[a,b]$ for $n\geq3$. Using~\eqref{eq:scal}, this is equivalent to $\gamma=v(t)^{2}dt^{2}+F(t)^{2}g$ for some fixed metric $g$ with $\int_{\mathbb{S}^{n}}\scal(g)dV_{g} = n(n-1)\omega_{n}r_{a}^{n-2}$ if $n\geq3$ and
\begin{align*}
0&=\scal(g)+\frac{n}{v(t)^2}\left(\frac{2v'(t)F(t)F'(t)}{v(t)}-\left((n-1)F'(t)^{2}+2F(t)F''(t)\right) \right)\\
&\quad-\frac{n(n-1)q^{2}}{r_{a}^{2n}F(t)^{2(n-1)}}-2\Lambda F(t)^{2}
\end{align*}
for all $t\in[a,b]$. Of course, pieces of the $(n+1)$-dimensional \RN\ of mass $m:\equiv\mathcal{M}(\lbrace t\rbrace\times\mathbb{S}^{n})$, charge $q$, and cosmological constant $\Lambda$ can be written in the above form by choosing $v(t)=\sqrt{p_{m,q,\Lambda}(t)}>0$, $F(t)=\tfrac{t}{r_{a}}$, and $g(t)=r_{a}^{2}g_{*}$ hence $g'=0$, where $g_{*}$ denotes the canonical metric on $\mathbb{S}^{2}$. Moreover, by Proposition~\ref{prop:RN} and~\eqref{eq:scal}, we know that they satisfy the above identity as well as equality in Asssumption~\eqref{eq:R} if $n\geq3$. Conversely, we can conclude that $\scal(g)$ must be constant on $\mathbb{S}^{n}$ so that $\scal(g)=\tfrac{n(n-1)}{r_{a}^{2}}$ follows from the Gauss--Bonnet Theorem if $n=2$ and from equality in Asssumption~\eqref{eq:R} if $n\geq3$. Setting $r(t)\definedas F(t)r_{a}$ for all $t\in[a,b]$ and using that $r'>0$ on $[a,b]$ follows from the assumption $F'>0$ on $[a,b]$, the implicitly defined function $\varphi\colon r([a,b])\to\R^{+}$, $\varphi(r(t))\definedas \tfrac{v(t)}{r'(t)}$ is well-defined and smooth. Constancy of $m:\equiv\mathcal{M}(\lbrace t\rbrace\times\mathbb{S}^{n})$ thus gives $p_{m,q,\Lambda}(r)=p_{o,q,\Lambda}(r)-\frac{2m}{r^{n-1}}=\frac{1}{\varphi(r)^{2}}>0$ for all $r\in r([a,b])$ via~\eqref{eq:Hawkingsimple}. This tells us that 
\begin{align*}
\gamma&=v(t)^{2}dt^{2}+F(t)^{2}g=\varphi(r)^{2}dr^{2}+\frac{r^{2}}{r_{a}^{2}}g=\frac{1}{p_{m,q,\Lambda}(r)}dr^{2}+\frac{r^{2}}{r_{a}^{2}}g=\frac{1}{p_{m,q,\Lambda}(r)}dr^{2}+r^{2}h_{*},
\end{align*}
with $h_{*}\definedas r_{a}^{-2}g$ having volume radius $1$ and scalar curvature $\scal(h_{*})=n(n-1)$ on $\mathbb{S}^{n}$. Moreover, 
\begin{align*}
E&=\frac{q}{v(t)r_{a}^{n}F(t)^{n}}\partial_{t}=\frac{q}{\varphi(r)r^{n}}\partial_{r}=\frac{q\sqrt{p_{m,q,\Lambda}(r)}}{r^{n}}\partial_{r}=E_{m,q,\Lambda}
\end{align*}
so that $(M,\gamma,E,\Lambda)$ is isometric to $(r([a,b])\times\mathbb{S}^{2},\tfrac{1}{p_{m,q,\Lambda}(r)}dr^{2}+r^{2}h_{*},E_{m,q,\Lambda},\Lambda)$ as claimed.

Using the \sDEC~\eqref{eq:sDEC} and its evaluation on the collar via~\eqref{eq:scal} as well as $F'>0$ on $[a,b]$, the first of the above sequence of inequalities is indeed strict, showing strict monotonicity of the generalized Hawking mass.
\end{proof}

\begin{remark}[Rotational Symmetry]
Proposition~\ref{prop:monotonicity} in particular applies to rotationally symmetric \ECR s, where~\eqref{eq:R} is automatically satisfied.
\end{remark}
\newpage
\begin{remark}
Let us briefly discuss which metrics $h_{*}$ can arise in the rigidity statement of Proposition~\ref{prop:monotonicity}. For $n=2$, we know for example from the Uniformization Theorem (or from Corollary 10 in~\cite{Petersen}) that the Riemannian manifold $(\mathbb{S}^{2},h_{*})$ of volume radius $1$ and scalar curvature $\scal(h_{*})=2$ from Proposition~\ref{prop:monotonicity} is isometric to $(\mathbb{S}^{2},g_{*})$. For $n\geq3$, in the special case of Yamabe metrics, Assumption~\eqref{eq:R} is automatically satisfied with equality only in the round case, see for example \cite[Lemma~3.4]{LP}. One can hence conclude that if $h_{*}$ is a Yamabe metric, $(\mathbb{S}^{n},h_{*})$ must also be isometric to $(\mathbb{S}^{n},g_{*})$ as for $n=2$.

However, it is well-known that there are many metrics on $\mathbb{S}^{n}$ which have constant positive scalar curvature and are not Yamabe metrics, see for example Brendle's survey~\cite{Brendlesurvey}. In this case, it is conceivable that equality in~\eqref{eq:R} can hold without forcing $(\mathbb{S}^{n},h_{*})$ to be isometric to $(\mathbb{S}^{n},g_{*})$.
\end{remark}

\begin{remark}[Similarity to Monotonicity Under Smooth Inverse Mean Curvature Flow]
The monotonicity of the generalized Hawking mass proved in Proposition~\ref{prop:monotonicity} is actually not very different from monotonicity under (smooth) Inverse Mean Curvature Flow. More precisely, let $(M,\gamma)$ be a smooth, connected Riemannian manifold (possibly with boundary $\partial M\subset M$). Assume that $(M,g)$ is regularly foliated by a smooth Inverse Mean Curvature Flow $\lbrace\Sigma_{t}\rbrace_{t\in I}$ with spherical leaves, so that in particular $M\approx I\times\mathbb{S}^{n}$ for some interval $I$ and $\frac{d}{dt}X_{t}=\frac{1}{H_{t}}\nu_{t}\vert_{X_{t}}$ for all $t\in I$, where $X_{t}\colon\mathbb{S}^{n}\to M$, $X_{t}(\mathbb{S}^{n})=\Sigma_{t}$, is a smooth embedding, $\nu_{t}$ is a smooth choice of unit normals to $\Sigma_{t}$, and $H_{t}>0$ denotes the mean curvature of $\Sigma_{t}$ in $M$ with respect to $\nu_{t}$ for all $t\in I$. Choosing an appropriate gauge, namely vanishing shift, one can see that $(M,\gamma)$ is isometric to $(I\times\mathbb{S}^{n},v(t,\cdot)^{2}dt^{2}+F(t)^{2}g(t))$, with $F\colon I\to\R^{+}\colon t\mapsto \left(\tfrac{\vert\Sigma_{t}\vert_{h_{t}}}{\omega_{n}}\right)^{\frac{1}{n}}$ the volume radius function, $\lbrace g(t)= \tfrac{1}{F(t)^{2}}(X_{t})^{*}h_{t}\rbrace_{t\in I}$ a smooth path of metrics on $\mathbb{S}^{n}$, and $v(t,\cdot)=\frac{1}{H_{t}}\left(1+\frac{1}{2}\tr_{g(t)}g'(t)\right)$ determined by the mean curvature along the flow. To see why this definition of $v$ works, we refer the reader to the proof of Lemma~\ref{lem:meancurvature} and recall that the evolution equation for the volume $\vert\Sigma_{t}\vert_{h_{t}}$ under Inverse Mean Curvature Flow gives $\frac{d}{dt}\vert\mathbb{S}^{n}\vert_{h_{t}}=\vert\mathbb{S}^{n}\vert_{h_{t}}$ for all $t\in I$ which shows that $F'(t)=\tfrac{F(t)}{n}$. 

Having rewritten $(M,\gamma)$ in the form $(I\times\mathbb{S}^{n},v(t,\cdot)^{2}dt^{2}+F(t)^{2}g(t))$, we are almost in the situation studied in Proposition~\ref{prop:monotonicity}, except that we do not have $\tr_{g(t)}g'(t)=0$ along the path of metrics $\lbrace g(t)\rbrace_{t\in I}$. The other main difference between the general (smooth) Inverse Mean Curvature Flow case and the case considered here is that we have assumed a very specific form of the electric field. Proving monotonicity of the generalized Hawking mass under smooth, let alone weak Inverse Mean Curvature Flow in full generality (meaning for arbitrary divergence-free electric field and allowing $\tr_{g(t)}g'(t)\neq0$) would lead too far here but should not in principle lead to any problems.
\end{remark}

The above remark is complemented by the following observation which generalizes an observation made in~\cite{CPCworkinprogress} for $n=2$, $q=0$, and $\Lambda=0$.

\begin{proposition}[Collar Extensions Are a Reparametrized Inverse Mean Curvature Flow]
Let $n\geq2$ and consider a collar extension $(I\times\mathbb{S}^{n},v(t,\cdot)^{2}dt^{2}+F(t)^{2}g(t))$, with $I$ an interval, and $v\colon I\times\mathbb{S}^{n}\to\R^{+}$, $F\colon I\to\R^{+}$ smooth, positive functions with $F'>0$, and $\lbrace g(t)\rbrace_{t\in I}$ a smooth path of metrics with $\tr_{g(t)}g'(t)=0$ for all $t\in I$. Then the canonical embeddings $X_{t}\colon\mathbb{S}^{n}\to I\times\mathbb{S}^{n}\colon p\mapsto (t,p)$ are a reparametrized smooth Inverse Mean Curvature Flow, i.e., solve $\frac{d}{ds}X_{t(s)}=\frac{1}{H_{t(s)}}\nu_{t(s)}\vert_{X_{t(s)}}$ for all $s\in J$ for some smooth, bijective reparametrization $t\colon J\to I$ with $\tfrac{dt}{ds}>0$, where $\nu_{t}=\tfrac{1}{v(t,\cdot)}\partial_{t}$ denotes a smooth choice of unit normal and $H_{t}$ denotes the mean curvature of $\lbrace t\rbrace\times\mathbb{S}^{n}$ in $(I\times\mathbb{S}^{n},v(t,\cdot)^{2}dt^{2}+F(t)^{2}g(t))$ with respect to $\nu_{t}$.
\end{proposition}
\begin{proof}
By definition, we have $\tfrac{d}{dt}X_{t}\vert_{p}=\partial_{t}\vert_{X_{t}(p)}=v(t,p)\nu_{t}\vert_{X_{t}(p)}$ for all $(t,p)\in I\times\mathbb{S}^{n}$. From Lemma~\ref{lem:meancurvature}, we know that $v(t,\cdot)=\tfrac{nF'(t)}{(H_{t}\circ X_{t})F(t)}$ so that $\tfrac{d}{dt}X_{t}=\tfrac{nF'(t)}{H_{t}F(t)}\nu_{t}\vert_{X_{t}}$ for all $t\in I$. Now let $s\colon I\to\R\colon t\mapsto n\ln F(t)$, $J\definedas s(I)$ and note that $\tfrac{ds}{dt}=\tfrac{nF'(t)}{F(t)}>0$ for all $t\in I$ by assumption. Hence the inverse of $s$, $t\colon J\to I$, has $\tfrac{dt}{ds}=\tfrac{F(t)}{nF'(t)}>0$ which gives $\tfrac{d}{ds}X_{t(s)}=\tfrac{nF'(t)}{H_{t(s)}F(t)}\tfrac{dt}{ds}(s)\nu_{t}\vert_{X_{t(s)}}=\tfrac{1}{H_{t(s)}}\nu_{t}\vert_{X_{t(s)}}$ for all $s\in J$ as desired.
\end{proof}

It is not clear whether and actually seems rather unlikely that Proposition~\ref{prop:monotonicity} applies to the collar extensions constructed in Lemma~\ref{lem:collar} when $n\geq3$, except of course in the spherically symmetric case. However, strict monotonicity of the generalized Hawking mass in collar extensions with minimal inner boundary and $v=\text{const}$ can be asserted, provided that $A$ is chosen suitably large and that the quasi-local sub-extremality condition $h_{q,\Lambda}>0$ holds at the inner boundary as the following proposition shows. Note that the \DEC~\eqref{eq:DEC} is not assumed explicitly. This can be used to obtain an alternative proof of the lower estimate on the generalized Hawking mass in Lemma~\ref{lem:Hawking*}.

\begin{proposition}[Generalized Hawking Mass Monotonicity Revisited]\label{prop:monotonicity2}
Let $n\geq2$ and let $\{g(t) \}_{t\in[a,b]}$ be a smooth path of metrics on $\mathbb{S}^{n}$ satisfying $\tr_{g(t)}g'(t)=0$ for all $t\in[a,b]$, with $r_{a}$ denoting the volume radius of $g(a)$, i.e., $\omega_{n} r_{a}^{2}\definedas\vert\lbrace a\rbrace\times\mathbb{S}^{n}\vert_{g(a)}$. Set $M\definedas[a,b]\times\mathbb{S}^{n}$ and let $F\colon [a,b] \to \mathbb{R}^{+}$  be a smooth, positive function with $F'(t)>0$ for $t\in(a,b]$ and $F'(a)=0$, $F''(a)>0$. Let $\gamma_{A}$ be the smooth Riemannian metric given by $\gamma_{A} = A^2 dt^2 +F(t)^2 g(t)$ on $M$ for some constant $A>0$. Set $E_{A}\definedas \frac{q}{A r_a^n F(t)^n} \partial_t$ for some constant $q\in\R$, and let $\Lambda\leq0$. Assume that $h_{q,\Lambda}(r_{a}F(a))>0$ and consider the \ECR\ $\mathcal{S}_{A}=(M,\gamma_{A},E_{A},\Lambda)$. Then there exists $A_{1}>0$ depending only on $n$, $q$, $\Lambda$, $r_{a}$, $F(a)$, $F''(a)$, $\max_{[a,b]}\lbrace F^{n-2}F'^{2}\rbrace$, and $\max_{[a,b]}\lbrace F^{n-1}\vert F''\vert\rbrace$ such that the generalized Hawking mass functional $\mathcal{M}_{\mathcal{S}_{A}}$ in $\mathcal{S}_{A}$ satisfies $\tfrac{d}{dt}\mathcal{M}_{\mathcal{S}_{A}}(\lbrace t\rbrace\times\mathbb{S}^{2})>0$ on $(a,b]$ for all $A\geq A_{1}$.
\end{proposition}
\begin{proof}
Recall from the proof of Proposition~\ref{prop:monotonicity} that
\begin{align*}
\frac{d}{dt}\mathcal{M}_{\mathcal{S}_{A}}(\lbrace t\rbrace\times\mathbb{S}^{n})= \,& \frac{r_{a}^{n-1}F'(t)}{2}\\
&\times\left\lbrace\frac{(n-1)h_{q,\Lambda}(r_{a}F(t))}{r_{a}^{2(n-1)}F(t)^{n}}-\frac{r_{a}^{2}F(t)^{n-2}\left[(n-1)F'(t)^{2}+2F(t)F''(t)\right]}{A^{2}}\right\rbrace
\end{align*}
holds for all $t\in[a,b]$, where we have used that $v(t,\cdot)=A$ so that $v'(t,\cdot)=0$ and have simplified the resulting expression accordingly. Hence, using $F'(a)=0$, we find
\begin{align*}
\frac{d^{2}}{dt^{2}}\mathcal{M}_{\mathcal{S}_{A}}(\lbrace t\rbrace\times\mathbb{S}^{n})\vert_{t=a}= \underbrace{\frac{r_{a}^{n-1}F''(a)}{2}}_{>0}\left\lbrace\underbrace{\frac{(n-1)h_{q,\Lambda}(r_{a}F(a))}{r_{a}^{2(n-1)}F(a)^{n}}}_{>0}-\frac{2r_{a}^{2}F(a)^{n-1}F''(a)}{A^{2}}\right\rbrace
\end{align*}
as we have assumed $F''(a)>0$ and $h_{q,\Lambda}(r_{a}F(a))>0$. Hence there exists $A_{2}>0$ depending only on $n$, $q$, $\Lambda$, $r_{a}$, $F(a)$, and $F''(a)$ such that $\frac{d^{2}}{dt^{2}}\mathcal{M}_{\mathcal{S}_{A}}(\lbrace t\rbrace\times\mathbb{S}^{n})\vert_{t=a}>0$ holds for all $A\geq A_{2}$. As $\frac{d}{dt}\mathcal{M}_{\mathcal{S}_{A}}(\lbrace t\rbrace\times\mathbb{S}^{n})\vert_{t=a}=0$ by $F'(a)$ follows from Proposition~\ref{prop:monotonicity}, this tells us that $\frac{d}{dt}\mathcal{M}_{\mathcal{S}_{A}}(\lbrace t\rbrace\times\mathbb{S}^{n})>0$ on $(a,a+\delta)$ for some $0<\delta\leq b-a$. Now suppose that $\frac{d}{dt}\mathcal{M}_{\mathcal{S}_{A}}(\lbrace t\rbrace\times\mathbb{S}^{n})\vert_{t=c}=0$ for some $a+\delta\leq c\leq b$. Then by strict monotonicity of $h_{q,\Lambda}$ (see Lemma~\ref{lem:h}) and $F'>0$ on $(a,b]$, we have
\begin{align*}
0&=\frac{(n-1)h_{q,\Lambda}(r_{a}F(c))}{r_{a}^{2(n-1)}F(c)^{n}}-\frac{r_{a}^{2}F(c)^{n-2}\left[(n-1)F'(c)^{2}+2F(c)F''(c)\right]}{A^{2}}\\
&> \frac{(n-1)h_{q,\Lambda}(r_{a}F(a))}{r_{a}^{2(n-1)}F(c)^{n}}-\frac{(n-1)r_{a}^{2}\max_{[a,b]}\lbrace F^{n-2}F'^{2}\rbrace+2r_{a}^{2}\max_{[a,b]}\lbrace F^{n-1}\vert F''\vert\rbrace}{A^{2}}.
\end{align*}
which can be ensured to be strictly positive provided that $A\geq A_{1}$, where $A_{1}\geq A_{2}$ can be chosen such that it depends only on the quantities mentioned in the statement of this lemma, where we have used the assumption that $h_{q,\Lambda}(r_{a}F(a))>0$. This is a contradiction showing that $\frac{d}{dt}\mathcal{M}_{\mathcal{S}_{A}}(\lbrace t\rbrace\times\mathbb{S}^{n})>0$ on $(a,b]$ for all $A\geq A_{1}$.
\end{proof}

\section{Gluing methods}\label{sec:gluing}
In this section, we will present the gluing methods necessary to glue the collar extension constructed in Section~\ref{sec:collar} to a \RN\ in a way that the resulting manifold is a smooth \ECR. We follow the general procedure of Mantoulidis and Schoen in \cite{pri}. The gluing results of this section are generalizations of results of \cite{impor} and \cite{hyper}. 

The following lemma provides a way to smoothly glue two rotationally symmetric \ECR s (with coinciding cosmological constants).
	
\begin{lemma}[Gluing Lemma]\label{lem:glue}
Let $n\geq2$ and let $\Lambda\leq0$, $q_i  \in \mathbb{R}$, $i=1,2$, be constants. Let $f_i \colon[a_i,b_i] \to  \mathbb{R}^{+}$, $i=1,2$, be two smooth, positive functions and consider the smooth Riemannian metrics $\gamma_i=ds^2 +f_i (s)^2 g_*$ on $[a_i, b_i] \times \mathbb{S}^n$, where $g_{*}$ denotes the canonical metric on $\mathbb{S}^{n}$, and the smooth electric vector fields $E_i=\frac{q_i}{f_i^n} \partial_s$, again for $i=1,2$. Assume that the \ECR s $([a_i, b_i] \times \mathbb{S}^n,\gamma_{i},E_{i},\Lambda)$, $i=1,2$, satisfy the strict Dominant Energy Condition~\eqref{eq:sDEC} and assume that
\begin{enumerate}
\item $f_1(b_1) < f_2(a_2)$ and $f'_1(b_1)\geq f'_2(a_2)$,
\item $\mathcal{M}_{1}(\lbrace b_{1}\rbrace\times\mathbb{S}^n) \geq \frac{q_1^{2}}{f_{1}(b_{1})^{n-1}}+\frac{2\Lambda f_1(b_1)^{n+1}}{n(n-1)(n+1)}$ and
\item $\mathcal{M}_{2}(\lbrace a_{2}\rbrace\times\mathbb{S}^n) \geq \frac{q_2^{2}}{f_{2}(a_{2})^{n-1}}+\frac{2\Lambda f_2(a_2)^{n+1}}{n(n-1)(n+1)}$,
\end{enumerate}
where $\mathcal{M}_{i}$ denotes the Hawking mass functional within the \ECR\ $([a_i, b_i] \times \mathbb{S}^n,\gamma_{i},E_{i},\Lambda)$ for $i=1,2$, respectively.

Then, after an appropriate translation of the interval $[a_2, b_2]$ so that
\begin{align}\label{inter}
\left\{
\begin{array}{ll}
(a_2 -b_1)f_1'(b_1)= f_2(a_2) -f_1(b_1), &\text{ if } \;  f_1'(b_1)=f'_2(a_2) \\
(a_2 -b_1)f_1'(b_1)> f_2(a_2) -f_1(b_1) > (a_2 -b_1)f_2'(a_2) &\text{ if }\; \, f_1'(b_1)\geq f'_2(a_2), \\
\end{array} 
\right. 
\end{align}
there exists a smooth, positive function $f\colon [a_1, b_2] \to \mathbb{R}^{+}$ such that 
\begin{enumerate}[label=(\roman*)]
\item $f \equiv f_1$ on $[a_1, \frac{a_1 +b_1}{2}]$,
\item $f \equiv f_2$ on $[\frac{a_2+b_2}{2}, b_2 ]$, and
\item for any $q\in\mathbb{R}$ satisfying $q^2 \leq \min \{q_1 ^2, q_2^2  \}$, with $q^2 < q_1^2$ if equality holds in Assumption~(2) and $q^2 < q_2^2$ if equality holds in Assumption~(3), the \ECR\ $([a_{1},b_{2}]\times\mathbb{S}^{n},\gamma \definedas ds^2 +f(s)^2 g_*,E\definedas \frac{q}{f(s)^n} \partial_s,\Lambda)$ satisfies the strict Dominant Energy Condition~\eqref{eq:sDEC}.
\end{enumerate}
If, in addition, $f'_i >0$ on $[a_i, b_i]$ for $i=1,2$, one can arrange that $f'>0$ on $[a_{1},b_{2}]$.
\end{lemma}

\begin{proof}
Using Assumption~(1), it is straightforward to see that Condition~\eqref{inter} can be achieved. It is then easy to show that there exists a smooth function $\zeta \colon[b_1, a_2] \to \mathbb{R} $ satisfying
\begin{enumerate}[label=(\alph*)]
\item $\zeta(b_1)= f_1'(b_1), $
\item $\zeta(a_2)= f_2'(a_2), $
\item $\zeta' \leq 0, $ and
\item $\int_{b_1}^{a_2} \zeta(x) dx = f_2(a_2) -f_1(b_1). $
\end{enumerate}

\noindent We define the function $\widehat{f} \colon[b_1,a_2] \to \mathbb{R}\colon s\mapsto f_1(b_1) + \int_{b_1}^{s} \zeta(x) dx$ and observe that it satisfies
\begin{itemize}
\item $\widehat{f}(b_1)=f_1(b_1)$ and $\widehat{f}(a_2)=f_2(a_2)$,
\item $\widehat{f}'(b_1)=f'_1(b_1)$ and $\widehat{f}'(a_2)=f'_2(a_2)$,
\item $f_1'(b_1)\geq\widehat{f}'\geq f_2'(a_2)$ on $(b_1,a_2)$, and
\item $\widehat{f}''=\zeta'\leq 0$ on $[b_1,a_2]$.
\end{itemize}

\noindent Next, we define the extension $\widetilde{f}$ of the functions $f_i$ by
\begin{align*}
\widetilde{f}&\definedas   \left\{
\begin{array}{ll}
f_1 & \text{ on }  \;  [a_1, b_1] \\
\widehat{f}  &\text{ on }\; [b_1, a_2] \\
f_2 & \text{  on } \;  [a_2, b_2] \\
\end{array} 
\right. .
\end{align*}
By construction, we have that $\widetilde{f}\in C^{1,1}([a_1,b_2])$ and that $\widetilde{f}$ is $C^2$ away from $b_1,a_2$.  In view of Lemma~\ref{le4.1}, we set
\begin{align*}
\Omega [f] &\definedas \frac{n-1}{2f}\left( \frac{h_{q,\Lambda}\circ f}{f^{2(n-1)}}- (f')^2 \right)
\end{align*}
and note that the \sDEC~\eqref{eq:sDEC} holds for some radial profile $f$ if and only if $\Omega[f]>f''$ (when keeping the dimensional parameter $n$, the charge $q$, and the cosmological constant $\Lambda$ fixed). Note that since $\widetilde{f} = f_i $ on $[a_i, b_i]$, $i=1,2$, we have $\Omega [\widetilde{f}\,] = \Omega[f_i] > f_i'' = \widetilde{f}\,''$ on $[a_1,b_1) \cup (a_2, b_2].$
		
Using Assumption~(2), $q^{2}\leq q_{1}^{2}$ with $q^{2}< q_{1}^{2}$ if equality holds in Assumption~(2), and Proposition~\ref{prop:sphereHawking}, in particular~\eqref{eq:Mp} with parameter $m=0$, we find
\begin{align*}
\frac{2 f_1(b_1)}{n-1} \Omega[\widetilde{f}\,](b_1)&=\frac{2 f_1(b_1)}{n-1} \Omega[f_1](b_1) = \frac{h_{q,\Lambda}(f_{1}(b_{1}))}{f_{1}^{2(n-1)}(b_{1})} -  f_1'(b_1)^2\\
&=\frac{h_{q,\Lambda}(f_{1}(b_{1}))}{f_{1}^{2(n-1)}(b_{1})}-p_{0,q_{1},\Lambda}(f_{1}(b_{1}))+\frac{2\mathcal{M}_{1}(\lbrace b_{1}\rbrace\times\mathbb{S}^{n})}{f_{1}(b_{1})^{n-1}}\\
&>\frac{h_{q_{1},\Lambda}(f_{1}(b_{1}))}{f_{1}^{2(n-1)}(b_{1})}-p_{0,q_{1},\Lambda}(f_{1}(b_{1}))+\frac{2q_{1}^{2}}{f_{1}(b_{1})^{2(n-1)}}+\frac{4\Lambda f_{1}(b_{1})^{2}}{n(n-1)(n+1)}=0
\end{align*}
so that $\Omega[\widetilde{f}\,](b_1)>0$. Similarly, using Assumption~(3) and the conditions on $q_{2}$, one can check that $\Omega[\widetilde{f}\,](a_2) >0$. 

As a next step, let us show that $\Omega[\widetilde{f}\,] >0$ on $(b_1, a_2)$ as well. On $(b_1, a_2)$, $\widetilde{f}\,''=\zeta' \leq 0 $ by Condition (c) so that $\widetilde{f}\vert_{[b_{1},a_{2}]}$ attains its minimum at $b_{1}$ by Assumption~(1) and attains its maximum at some $s_*$ in $(b_1, a_2]$. Moreover, $\widetilde{f}\,''\leq0$ implies that $\widetilde{f}\,'$ is non-increasing so  that $f'_1 (b_1) \geq \widetilde{f}\,'(s) \geq 0 $ for all $s\in(b_1, s_*]$. Using $\widetilde{f}(b_{1})=f_{1}(b_{1})$ and the positivity of $\Omega[\widetilde{f}\,](b_1)$, we get 
\begin{align*}
\frac{q^2}{\widetilde{f}(t)^{2(n-1)}} \leq \frac{q^2}{f_1 (b_1)^{2(n-1)}} < 1- f_1' (b_1)^2 - \frac{2\Lambda f_{1}(b_{1})^{2} }{n(n-1)}\leq 1- \widetilde{f}\,'(s)^2  - \frac{2\Lambda \widetilde{f}(s)^{2}}{n(n-1)}
\end{align*}
for all $b_1 < s,t \leq s_*$, where we have used that $\Lambda\leq0$. On the other hand, if $s_* < a_2$, we have $0 \geq \widetilde{f}\,'(s) \geq f_2 '(a_2) $ for all $s\in[s_*, a_2 )$, and hence using $\widetilde{f}(a_{2})=f_{2}(a_{2})$ and the positivity of $\Omega[\widetilde{f}\,](a_2)$, we get
\begin{align*}
\frac{q^2}{\widetilde{f}(t)^{2(n-1)}} \leq \frac{q^2}{f_2 (a_2)^{2(n-1)}}  < 1- f_2' (a_2)^2 - \frac{2\Lambda f_{2}(a_{2})^{2} }{n(n-1)}\leq 1- \widetilde{f}\,'(s)^2  - \frac{2\Lambda \widetilde{f}(s)^{2}}{n(n-1)}
\end{align*}
for all $s_* \leq s,t< a_2$, again using $\Lambda\leq0$. From these two inequalities we see that 
\begin{align*}
\Omega[\widetilde{f}\,] > 0\geq \widetilde{f}\,''
\end{align*}
on $(b_{1},a_{2})$. From this and the above considerations, we can deduce that $\Omega[\widetilde{f}\,] > \widetilde{f}\,''$ on $[a_1, b_2] \setminus \{b_1, a_2\} $. In order to deal with $s=b_{1}$ and $s=a_{2}$, let $d>0$ be defined by
\begin{align*}
3d &\definedas \inf_{[a_1, b_2] \setminus \{b_1, a_2\}}  \left( \Omega[\widetilde{f}\,] - \widetilde{f}\,'' \right)>0,
\end{align*}
so that $\widetilde{f}\,'' +3d \leq \Omega[\widetilde{f}\,]$ on $[a_1, b_2] \setminus \{b_1, a_2\} $. Let $\delta >0$ be such that $ \frac{a_1 +b_1}{2}  < b_1 -\delta $ and $a_2+ \delta <  \frac{a_2 +b_2}{2}$. Let $\eta_\delta$ be a smooth cut-off function such that $\eta_\delta=1$ on $[b_1-\delta, a_2+\delta ]$ and $\eta_\delta=0$ on $[a_1, \frac{a_1 +b_1}{2}] \cup[\frac{a_2 +b_2}{2},b_2 ]$. Define a mollification $f_{\varepsilon}\colon[a_{1},b_{2}]\to\mathbb{R}^{+}$ of $\widetilde{f}$ by
\begin{align}
f_\varepsilon(t) &\definedas \int_\mathbb{R} \widetilde{f}(t- \varepsilon \eta_\delta (t) s) \phi(s) ds,
\end{align}
where $\phi$ is a standard mollifier. The mollification $f_{\varepsilon}$ coincides with $\widetilde{f} $ on $[a_1, \frac{a_1 +b_1}{2}] \cup[\frac{a_2 +b_2}{2},b_2 ]$ and coincides with the standard mollification of radius $\varepsilon$ on an interval properly containing $[b_1,a_2]$. It is not hard to check that $f_\varepsilon \rightarrow \widetilde{f}$ in $C^1([a_1, b_2])$ as $\varepsilon \searrow 0$, which implies that $\Omega[f_\varepsilon] \rightarrow \Omega[\widetilde{f}\,] $  in  $C^0([a_1, b_2])$ as $\varepsilon\searrow 0$. Thus, for $\varepsilon$ sufficiently small, we have 
\begin{align*}
\sup_{[a_1, b_2]} \big| \Omega[\widetilde{f}\,]  -\Omega[f_\varepsilon]  \big|< d.
\end{align*}
Note that $\Omega[\widetilde{f}\,]$ is uniformly continuous; therefore, for small enough $\varepsilon>0$, we have that $\Omega[\widetilde{f}\,](t) \leq \Omega[\widetilde{f}\,](s) +d $ for $t\in [s-\varepsilon, s+\varepsilon]$. Abusing notation, set $\widetilde{f}\,''(b_1) =f_1''(b_1)$ and $\widetilde{f}\,''(a_2) =f_2''(a_2)$. Then for any sufficiently small $\varepsilon$ it follows for all $s\in[a_{1},b_{2}]$ that 
\begin{align*}
f''_\varepsilon(s) &\leq \sup_{[s-\varepsilon, s+ \varepsilon]} \widetilde{f}\,''(t)+d \leq \sup_{[s-\varepsilon, s+ \varepsilon]} \Omega[\widetilde{f}\,](t) -3d +d \leq  \Omega[\widetilde{f}\,](s) -d\\
&< \Omega[f_\varepsilon](s) =\frac{n-1}{2 f_\varepsilon (s)} \left(  \frac{h_{q,\Lambda}\circ f_{\varepsilon}}{f_{\varepsilon}^{2(n-1)}}- f_\varepsilon'(s)^2  \right).
\end{align*}
In other words, the \ECR\ $([a_{1},b_{2}]\times\mathbb{S}^{n},\gamma\definedas ds^2 +f(s)^2 g_*,E\definedas \frac{q}{f(s)^n} \partial_s,\Lambda))$, where $q$ satisfies $q^2 \leq \min \{q_1 ^2, q_2^2  \}$ satisfies the \sDEC~\eqref{eq:sDEC} if $f$ is chosen as $f=f_{\varepsilon}$ for suitably small fixed $\varepsilon>0$.

If $f'_i >0$ on $[a_i, b_i]$ for $i=1,2$ then $\widetilde{f}'>0$ by construction, hence its mollification $f_{\varepsilon}$ also has $f_{\varepsilon}'>0$ for suitably small $\varepsilon$ as $f_\varepsilon \rightarrow \widetilde{f}$ in $C^1([a_1, b_2])$ as $\varepsilon\searrow 0$.
\end{proof}

Lemma~\ref{lem:glue} is slightly more general than the corresponding Lemma 4.2 in \cite{impor} as it does not impose lower bounds on $f_{1}(b_{1})$ in terms of $q_{1}$ and on $f_{2}(a_{2})$ in terms of $q_{2}$. Moreover, it allows equality in the Hawking mass estimate, paying the price of asking that the charge of the glued solution is strictly smaller in this case. 

\begin{remark}
Taking a more geometric perspective on Assumptions~(1) through~(3) in Lemma~\ref{lem:glue}, we see that (1) corresponds to asking that the volume radius of $\lbrace b_{1}\rbrace\times\mathbb{S}^{n}$ should be strictly smaller than that of $\lbrace a_{2}\rbrace\times\mathbb{S}^{n}$, while the mean curvature $H_{b_{1}}$ of $\lbrace b_{1}\rbrace\times\mathbb{S}^{n}$ shall be at least as big as that of $\lbrace a_{2}\rbrace\times\mathbb{S}^{n}$, $H_{b_{1}}\geq H_{a_{2}}$, (both with respect to the outward pointing normals $\nu=\partial_{s}$), see Lemma~\ref{lem:meancurvature}. By  definition of the generalized Hawking mass and again Lemma~\ref{lem:meancurvature}, Assumptions~(2) and~(3) can be viewed as assuming upper bounds on $H_{b_{1}}^{2}$ and $H_{a_{2}}^{2}$ in terms of the volume radii $f_{1}(b_{1})$, $f_{2}(a_{2})$, the charges $q_{i}$, and $\Lambda$, respectively. 

Moreover, note that the lower bounds on the generalized Hawking masses assumed in (2) and (3) are manifestly positive when $\Lambda=0$ but can be non-positive when $\Lambda<0$.
\end{remark}

Our next objective is to glue the collar extension constructed in Lemma~\ref{lem:collar} to a \RN\ via Lemma~\ref{lem:glue}. However, the \RN s saturate the \DEC~\eqref{eq:DEC} (see Proposition~\ref{prop:RN}) and hence do not satisfy it strictly, an assumption necessary to apply Lemma~\ref{lem:glue} and indeed necessary to allow for the mollification performed in its proof. We will hence first need to slightly bend the region of the \RN s we would like to glue to in order to ensure that their bent version satisfies the \sDEC~\eqref{eq:sDEC}, while keeping them the same outside the range needed for the gluing.

\begin{lemma}[Bending Lemma]\label{lem:bend}
Let $n\geq2$, let $a\geq0$ and $q\in\mathbb{R}$ be real constants, and let $f\colon[a,\infty)\to\mathbb{R}^{+}$ be a smooth, positive function. Consider the rotationally symmetric \ECR\ $\mathcal{S}\definedas([a, \infty )\times \mathbb{S}^n,\gamma = ds^2 +f(s)^2 g_{*},E= \frac{q}{f(s)^n} \partial_s,\Lambda)$ with radial profile $f$ and assume it satisfies the Dominant Energy Condition~\eqref{eq:DEC}.

Then for any $s_0>a$ with $f'(s_0) >0$, there exist $\delta >0$ and a smooth, positive function $\widetilde{f}\colon[s_0-\delta,\infty)\to\mathbb{R}^{+}$ such that the rotationally symmetric \ECR\ $\widetilde{\mathcal{S}}\definedas([s_0-\delta, \infty )\times \mathbb{S}^n,\widetilde{\gamma} = ds^2 +\widetilde{f}(s)^2 g_{*},\widetilde{E}= \frac{q}{{\widetilde{f}(s)}^n} \partial_s,\Lambda)$ satisfies the conditions
\begin{enumerate}
\item $\mathcal{S}$ and $\widetilde{\mathcal{S}}$ coincide when restricted to $[s_0,\infty)\times\mathbb{S}^{n}$, 
\item $\widetilde{f}\,'>0$ on $[s_0-\delta,s_0)$, and
\item $\widetilde{\mathcal{S}}$ restricted to $[s_0-\delta,s_0)\times\mathbb{S}^{n}$ satisfies the strict Dominant Energy Condition~\eqref{eq:sDEC}.
\end{enumerate}
Moreover, if $f(s_0) >\alpha>0$ for some constant $\alpha\in\R$ then one can ensure that $\widetilde{f}(s_0-\delta)> \alpha$. If $f''(s_0)>0$ then one can ensure that $\widetilde{f}\,'(s_0-\delta)< f'(s_0)$. 
\end{lemma}

\begin{proof}
Note that the claim holds trivially if the \sDEC\ already holds at $s=s_0$. Otherwise we will deform the metric $\gamma$ to increase its scalar curvature, adjusting the electric field accordingly in order to preserve its charge $q$ and the fact that it is divergence-free. To do so, for some $0<\delta \leq s_0-a$ to be determined, consider the smooth, positive function $\sigma\colon [s_0 -\delta,s_0)\to\mathbb{R}^{+}$ given by
\begin{align*}
\sigma(s) &\definedas \int_{s_0-\delta}^s \left( 1+ e^{-(t-s_0)^{-2}} \right) dt +K_\delta,\quad
K_\delta\definedas s_0 - \int_{s_0-\delta}^{s_0} \left( 1+ e^{-(t-s_0)^{-2}} \right) dt.
\end{align*}
 By construction, we have $\lim_{s\searrow s_0}\sigma(s) =s_0$ and therefore $\sigma$ can smoothly be extended by the identity $\sigma(s) =s$ for $s \geq s_0$. Note that $\sigma$ is strictly increasing on $(s_0-\delta,\infty)$. Now set $\widetilde{f}\definedas f\circ\sigma\colon(s_0-\delta,\infty)\to\mathbb{R}^{+}$ and denote the rotationally symmetric \ECR\ $\widetilde{\mathcal{S}}$ with radial profile $\widetilde{f}$ as in the statement of the lemma. Then clearly Claim~(1) has been achieved, regardless of our choice of $0<\delta \leq s_0-a$.
 
From Lemma~\ref{le4.1}, we know that $\widetilde{\mathcal{S}}$, restricted to $[s_0-\delta,s_0)\times\mathbb{S}^{n}$ satisfies the \sDEC~\eqref{eq:sDEC} if \eqref{s4.1} holds when applied to $\widetilde{f}\vert_{[s_0-\delta,s_0)\times\mathbb{S}^{n}}$ or in other words if
\begin{align*}
\widetilde{f}\,'' &< \frac{n-1}{2\widetilde{f}}\left(  \frac{h_{q,\Lambda}\circ \widetilde{f}}{\widetilde{f}^{2(n-1)}}- (\widetilde{f}\,')^2  \right).
\end{align*}
Denoting $\frac{d}{ds}$ by $\dot{} $ and $\frac{d}{d\sigma}$ by $'$, respectively, we have $\widetilde{f}\,'=(f'\circ\sigma)\dot{\sigma}$, $\widetilde{f}\,''=(f''\circ\sigma)\dot{\sigma}^{2}+(f'\circ\sigma)\ddot{\sigma}$, and hence the \sDEC\ holds on  $[s_0-\delta,s_0)\times\mathbb{S}^{n}$ provided we can show
\begin{align*}
(f''\circ\sigma)\dot{\sigma}^{2}+(f'\circ\sigma)\ddot{\sigma} &< \frac{n-1}{2(f\circ\sigma)}\left(  \frac{h_{q,\Lambda}\circ f\circ\sigma}{(f\circ\sigma)^{2(n-1)}}- (f'\circ\sigma)^{2}\dot{\sigma}^2 \right)
\end{align*}
on $[s_0-\delta,s_0)$. Using that $\mathcal{S}$ satisfies the \DEC, we know from Lemma~\ref{le4.1} that it satisfies~\eqref{4.1}, so that, using $\dot{\sigma}>0$, we get
\begin{align*}
(f''\circ\sigma)\dot{\sigma}^{2} &\leq \frac{n-1}{2(f\circ\sigma)}\left(  \frac{h_{q,\Lambda}\circ f\circ\sigma}{(f\circ\sigma)^{2(n-1)}}- (f'\circ\sigma)^{2}  \right)\dot{\sigma}^2\\
&= \frac{n-1}{2(f\circ\sigma)}\left(  \frac{h_{q,\Lambda}\circ f\circ\sigma}{(f\circ\sigma)^{2(n-1)}}- (f'\circ\sigma)^{2}\dot{\sigma}^2\right) -\left(1-\dot{\sigma}^{2}\right) \frac{(n-1)h_{q,\Lambda}\circ f\circ\sigma}{2(f\circ\sigma)^{2n-1}}.
\end{align*}
Hence we get
\begin{align*}
\widetilde{f}\,'' \leq&\, \frac{n-1}{2\widetilde{f}}\left( \frac{h_{q,\Lambda}\circ \widetilde{f}}{\widetilde{f}^{2(n-1)}}- (\widetilde{f}\,')^2 \right)-\left(1-\dot{\sigma}^{2}\right) \left[\frac{(n-1)h_{q,\Lambda}\circ \widetilde{f}}{2\widetilde{f}^{2n-1}}\,\right]+(f'\circ\sigma)\ddot{\sigma}.
\end{align*}
Now, as $s\to s_0$, $\left(1-\dot{\sigma}^{2}\right)(s)\approx-2e^{-(s-s_0)^{-2}}$, the term in square brackets $[\cdot]$ is bounded, and $\ddot{\sigma}(s)\approx\frac{2e^{-(s-s_0)^{-2}}}{(s-s_0)^{3}}<0$, while $f'\circ\sigma>0$ on $[s_0-\delta,s_0)$ for suitably small $0<\delta\leq s_0-a$ by assumption and by continuity of $f'$. Hence the middle term in the above equation vanishes faster than the last term which shows that the sum of both will converge to zero strictly from below as $s\nearrow s_0$. This shows that $\widetilde{\mathcal{S}}$ satisfies the \sDEC\ on $[s_0-\delta,s_0)\times\mathbb{S}^{n}$ by Lemma~\ref{le4.1}, for suitably small $0<\delta\leq s_0-a$. As $\dot{\sigma}>0$ and $f'\circ\sigma>0$ on $[s_0-\delta,s_0)$ by our choice of $\delta$, Claim~(2) follows as well.

Now assume that $f(s_0)> \alpha$ for some constant $\alpha>0$. Choosing $\delta>0$ sufficiently small gives $\widetilde{f}(s_0-\delta) > \alpha$ by continuity. Next, assume that $f''(s_0)>0$. By smoothness of $f$, $\sigma$, and the facts that $\lim_{s\to s_0}\dot{\sigma}=1$, $\lim_{s\to s_0}\ddot{\sigma}=0$, we deduce that $\widetilde{f}\,''=(f''\circ\sigma)\dot{\sigma}^{2}+(f'\circ\sigma)\ddot{\sigma}>0$ when $\delta>0$ is sufficiently small. Thus $\widetilde{f}\,'$ is strictly increasing on $[s_0-\delta, s_0]$ so that $ \widetilde{f}\,'(s_0-\delta)<\widetilde{f}\,'(s_0)=f'(s_0)$ as claimed. 
\end{proof}
\newpage
The following proposition combines the previous results on gluing (Lemma~\ref{lem:glue}) and bending (Lemma~\ref{lem:bend}) with estimates on the generalized Hawking mass $\mathcal{M}$ of the coordinate spheres. It allows to glue a given rotationally symmetric \ECR\ with suitable properties to a \RN.

\begin{proposition}[Gluing to $\Lambda$-Reissner--Nordstr\"om Manifolds]\label{prop:gluingRN}
Let $n\geq2$ and let $\Lambda\leq0$, $q \in \mathbb{R}$ be constants. Let $f \colon[a,b] \to  \mathbb{R}^{+}$ be a smooth, positive function. Consider the rotationally symmetric \ECR\ 
\begin{align*}
\mathcal{S}\definedas([a,b]\times \mathbb{S}^n,\gamma = ds^2 +f(s)^2 g_{*},E= \frac{q}{f(s)^n} \partial_s,\Lambda)
\end{align*}
with radial profile $f$ and assume that it satisfies the strict Dominant Energy Condition~\eqref{eq:sDEC}. Assume furthermore that
\begin{enumerate}[label=(\roman*)]
\item $m_{*}\definedas\mathcal{M}(\lbrace b\rbrace\times\mathbb{S}^n)\geq\frac{q^{2}}{f(b)^{n-1}}+\frac{2\Lambda f(b)^{n+1}}{n(n-1)(n+1)}$, and
\item $f'(b)>0$.
\end{enumerate} 
Then for any constants $m_e\geq m_{*}$, $q_{e}^{2}\leq q^{2}$, with $m_{e}>m_{*}$ or $q_{e}^{2}< q^{2}$, in particular $q_{e}^{2}<q^{2}$ if equality holds in Assumption~(i), there exists a rotationally symmetric \ECR\ $\widetilde{\mathcal{S}}=(\widetilde{M}=[a,\infty)\times\mathbb{S}^{n},\widetilde{\gamma}= ds^2 +\widetilde{f}(s)^2 g_{*},\widetilde{E}=\frac{q_{e}}{{\widetilde{f}(s)}^n} \partial_s,\Lambda)$ satisfying the Dominant Energy Condition~\eqref{eq:DEC} with the additional properties
\begin{enumerate}[label=(\Roman*)]
\item $\widetilde{f}=f$ on $[a,\frac{a+b}{2}]$,
\item $([c,\infty)\times\mathbb{S}^{n}\subset\widetilde{M},\widetilde{\gamma},\widetilde{E},\Lambda)$ is isometric to $([r,\infty)\times\mathbb{S}^n\subset M_{m_{e},q_{e},\Lambda},\gamma_{m_{e},q_{e},\Lambda},E_{m_{e},q_{e},\Lambda},\Lambda)$ for some $c>b$ and some $r>r_+=r_{+}(m_{e},q_{e},\Lambda)$.
\end{enumerate}

If, in addition, $f'>0$ on $[a,b]$ then $(\widetilde{M},\widetilde{\gamma})$ is foliated by strictly mean convex spheres that eventually coincide with coordinate spheres in $(M_{m_{e},q_{e},\Lambda},\gamma_{m_{e},q_{e},\Lambda})$ under this isometry.
\end{proposition}	

\begin{proof}
We want to apply Lemma~\ref{lem:glue} to glue $\mathcal{S}$ to the \ECR\ $\mathcal{S}_{m_{e},q_{e},\Lambda}\definedas(M_{m_{e},q_{e},\Lambda},\gamma_{m_{e},q_{e},\Lambda},E_{m_{e},q_{e},\Lambda},\Lambda)$. However, as $\mathcal{S}_{m_{e},q_{e},\Lambda}$ saturates the \DEC~\eqref{eq:DEC} by Proposition~\ref{prop:RN}, we will first need to apply Lemma~\ref{lem:bend} to it in order to enforce the \sDEC~\eqref{eq:sDEC} in the gluing region. If $\mathcal{S}_{m_{e},q_{e},\Lambda}$ is sub-extremal, we can apply Proposition~\ref{prop:extendRN} to rewrite and extend $\mathcal{S}_{m_{e},q_{e},\Lambda}$ to the rotationally symmetric \ECR\ 
\begin{align*}
\overline{\mathcal{S}}_{m_{e},q_{e},\Lambda}\definedas([0,\infty)\times\mathbb{S}^{n},\gamma_{m_{e},q_{e},\Lambda}=ds^{2}+u_{m_{e},q_{e},\Lambda}(s)^{2}g_{*},E_{m_{e},q_{e},\Lambda}=\frac{q_{e}}{u_{m_{e},q_{e},\Lambda}(s)^{n}}\partial_{s},\Lambda)
\end{align*}
with radial profile $u_{m_{e},q_{e},\Lambda}$ of which we know that $u_{m_{e},q_{e},\Lambda}'>0$ on $(0,\infty)$. If, on the other hand, $\mathcal{S}_{m_{e},q_{e},\Lambda}$ is extremal or super-extremal, then we can rewrite a suitable piece of it as the 
rotationally symmetric \ECR\ 
\begin{align*}
\mathcal{S}^{\mu}_{m_{e},q_{e},\Lambda}\definedas([0,\infty)\times\mathbb{S}^{n},\gamma_{m_{e},q_{e},\Lambda}=ds^{2}+u_{m_{e},q_{e},\Lambda}^{\mu}(s)^{2}g_{*},E_{m_{e},q_{e},\Lambda}=\frac{q_{e}}{u_{m_{e},q_{e},\Lambda}^{\mu}(s)^{n}}\partial_{s},\Lambda)
\end{align*}
by Remark~\ref{rem:notextend}, for any $\mu>r_{+}(m_{e},q_{e},\Lambda)$ in the extremal, and any $\mu>0$ in the super-extremal case, and with $u_{m_{e},q_{e},\Lambda}^{\mu}\geq\mu$ and ${u_{m_{e},q_{e},\Lambda}^{\mu}}'>0$ on $[0,\infty)$. To unify all cases, we write $\mathcal{S}_{m_{e},q_{e},\Lambda}^{\mu}\definedas\overline{\mathcal{S}}_{m_{e},q_{e},\Lambda}$ and $u_{m_{e},q_{e},\Lambda}^{\mu}\definedas u_{m_{e},q_{e},\Lambda}$ with $\mu\definedas r_{+}(m_{e},q_{e},\Lambda)$ in the sub-extremal case. 

In order to apply the Bending Lemma~\ref{lem:bend}, it will be convenient to pick $s_0>0$ such that $u_{m_{e},q_{e},\Lambda}^{\mu}(s_0)>f(b)$ and $0<{u_{m_{e},q_{e},\Lambda}^{\mu}}'(s_0)< f'(b)$. For suitable $\mu$, this can be achieved using Assumption~\emph{(ii)} and arguing separately in two cases.

\subsubsection*{In case $f(b)\in u_{m_{e},q_{e},\Lambda}^{\mu}([0,\infty))$,} we let $\varepsilon>0$ and pick $s_{\varepsilon}^{\mu}>0$ such that $u_{m_{e},q_{e},\Lambda}^{\mu}(s_{\varepsilon})=f(b)+\varepsilon$ which can be achieved with a unique $s_{\varepsilon}^{\mu}$ as $u_{m_{e},q_{e},\Lambda}^{\mu}$ is bijective onto its image. Then
\begin{align*}
{u_{m_e,q_{e},\Lambda}}'(s_\varepsilon^{\mu})^2&=p_{m_e,q_{e},\Lambda}(u_{m_e,q_{e},\Lambda}^{\mu}(s_\varepsilon^{\mu})) =p_{m_e,q_{e},\Lambda}(f(b)+\varepsilon) =p_{m_e,q_{e},\Lambda}(f(b)) +\mathcal{O}(\varepsilon)\\
& = p_{m_*,q_{e},\Lambda}(f(b)) -\frac{2(m_e-m_*)}{f(b)^{n-1}}+\mathcal{O}(\varepsilon)\\
& = p_{m_*,q,\Lambda}(f(b)) +\frac{q_{e}^{2}-q^{2}}{f(b)^{2(n-1)}}-\frac{2(m_e-m_*)}{f(b)^{n-1}}+\mathcal{O}(\varepsilon)
\end{align*} 
as $\varepsilon\searrow0$, using Taylor's Formula. Using Proposition~\ref{prop:sphereHawking} and in particular Formula~\eqref{eq:Mp} with $m=m_{*}$, we find that $p_{m_*,q,\Lambda}(f(b))= f'(b)^2$ so that $0<{u_{m_e,q_{e},\Lambda}^{\mu}}'(s_\varepsilon^{\mu})< f'(b)$ for suitably small $\varepsilon>0$ since $m_{e}\geq m_{*}$ and $q_{e}^{2}\leq q^{2}$ with $m_{e}>m_{*}$ or $q_{e}^{2}<q^{2}$. In the super-extremal case, $f(b)\in u_{m_{e},q_{e},\Lambda}^{\mu}([0,\infty))$ can be achieved by setting $\mu\definedas f(b)>0$. In the extremal case, $f(b)\in u_{m_{e},q_{e},\Lambda}^{\mu}([0,\infty))$ can be achieved provided that $f(b)>r_{+}(m_{e},q_{e},\Lambda)$, again by setting $\mu\definedas f(b)$. In the sub-extremal case, $f(b)\in u_{m_{e},q_{e},\Lambda}^{\mu}([0,\infty))$ can be achieved provided $f(b)\geq r_{+}(m_{e},q_{e},\Lambda)=\mu$. We pick one suitably small $\varepsilon>0$ and set $s_0\definedas s_{\varepsilon}^{\mu}$ for these choices of $\mu$ and $\varepsilon$. This shows that we can pick $s_{o}>0$ such that $u_{m_{e},q_{e},\Lambda}^{\mu}(s_0)>f(b)$ and ${u_{m_{e},q_{e},\Lambda}^{\mu}}'(s_0)<f'(b)$ for our choices of $\mu$ in the super-extremal case and conditionally in the extremal and sub-extremal cases when $f(b)> r_{+}(m_{e},q_{e},\Lambda)$ and $f(b)\geq r_{+}(m_{e},q_{e},\Lambda)$, respectively.

\subsubsection*{In case $f(b)\notin u_{m_{e},q_{e},\Lambda}^{\mu}([0,\infty))$,} we know that $f(b)<u_{m_{e},q_{e},\Lambda}^{\mu}(s_{o})$ holds automatically for any $s_{o}>0$. In the sub-extremal case, this happens precisely when $f(b)< r_{+}(m_{e},q_{e},\Lambda)$. Then of course by ${u_{m_{e},q_{e},\Lambda}^{\mu}}'(0)=0$ and because ${u_{m_{e},q_{e},\Lambda}^{\mu}}'$ is strictly increasing, $0<u_{m_{e},q_{e},\Lambda}'(s_0)<f'(b)$ will hold by continuity for suitably small $s_0>0$. In the extremal case, we can assume that $f(b)\leq r_{+}(m_{e},q_{e},\Lambda)$ as the other case has already been covered. Recall from Remark~\ref{rem:notextend} that ${u_{m_{e},q_{e},\Lambda}^{\mu}}'(0)=(p_{m_{e},q_{e},\Lambda}(\mu))^{\frac{1}{2}}\searrow0$ as $\mu\searrow r_{+}(m_{e},q_{e},\Lambda)$, hence  $0<{u_{m_{e},q_{e},\Lambda}^{\mu}}'(0)<f'(b)$ will hold for suitably small $\mu>r_{+}(m_{e},q_{e},\Lambda)$. Pick one such $\mu$. As ${u_{m_{e},q_{e},\Lambda}^{\mu}}'$ is strictly increasing, $0<{u_{m_{e},q_{e},\Lambda}^{\mu}}'(s_0)<f'(b)$ will hold by continuity for suitably small $s_0>0$. This shows that we can pick $s_{o}>0$ and $\mu>r_{+}(m_{e},q_{e},\Lambda)$ such that $u_{m_{e},q_{e},\Lambda}^{\mu}(s_0)>f(b)$ and ${u_{m_{e},q_{e},\Lambda}^{\mu}}'(s_0)<f'(b)$ also in this case.

By the Bending Lemma~\ref{lem:bend}, we know that for the parameters $s_0>0$ and $\mu$ selected above, there exist $\delta>0$ and a smooth, positive function $\widetilde{u}_{m_e,q_{e},\Lambda}\colon[s_0-\delta,\infty)\to\R^{+}$ such that the rotationally symmetric \ECR\ 
\begin{align*}
\widetilde{\mathcal{S}}_{m_e,q_{e},\Lambda}\definedas([s_0-\delta,\infty)\times\mathbb{S}^{n},\widetilde{\gamma}=ds^{2}+\widetilde{u}_{m_e,q_{e},\Lambda}(s)^{2}g_{*},\widetilde{E}=\frac{q}{\widetilde{u}_{m_e,q_{e},\Lambda}(s)^{n}}\partial_{s},\Lambda)
\end{align*}
with radial profile $\widetilde{u}_{m_e,q_{e},\Lambda}$ coincides with $\mathcal{S}_{m_{e},q_{e},\Lambda}^{\mu}$ when restricted to $[s_0,\infty)\times\mathbb{S}^{n}$, meaning that $\widetilde{u}_{m_e,q_{e},\Lambda}=u_{m_{e},q_{e},\Lambda}^{\mu}$ on $[s_0,\infty)$, and satisfies the \sDEC~\eqref{eq:sDEC} when restricted to the bending region $[s_0-\delta,s_0)\times\mathbb{S}^{n}$. Moreover, we know that $\widetilde{u}_{m_e,q_{e},\Lambda}'>0$ on $[s_0-\delta,\infty)$ as $\widetilde{u}_{m_e,q_{e},\Lambda}=u_{m_{e},q_{e},\Lambda}^{\mu}$ on $(s_0,\infty)$, see also Proposition~\ref{prop:extendRN} and Remark~\ref{rem:notextend}, respectively.

Furthermore, still by Lemma~\ref{lem:bend}, we get $\widetilde{u}_{m_e,q_{e},\Lambda}(s_0-\delta)>f(b)$ and $\widetilde{u}_{m_e,q_{e},\Lambda}'(s_0-\delta)<f'(b)$ from $u_{m_{e},q_{e},\Lambda}^{\mu}(s_0)>f(b)$, ${u_{m_{e},q_{e},\Lambda}^{\mu}}'(s_0)<f'(b)$, and the fact that ${u_{m_{e},q_{e},\Lambda}^{\mu}}''>0$ on $[0,\infty)$ (see Proposition~\ref{prop:extendRN} and Remark~\ref{rem:notextend}).

We now want to apply the Gluing Lemma~\ref{lem:glue} to the given \ECR\ $\mathcal{S}$ with radial profile $f_{1}=f \colon[a_{1}=a,b_{1}=b] \to  \mathbb{R}^{+}$ and the \ECR\ $\widetilde{\mathcal{S}}_{m_e,q_{e},\Lambda}$ with radial profile $f_{2}=\widetilde{u}_{m_e,q_{e},\Lambda}\vert_{[s_0-\delta,s_{1}]}\colon[a_{2}=s_0-\delta,b_{2}=s_{1}]\to\R^{+}$ for some $s_0-\delta<s_{1}<s_0$, e.g.~$s_{1}=s_0-\frac{\delta}{2}$. Note that we have restricted $\widetilde{u}_{m_e,q_{e},\Lambda}$ to a compact interval on which $\widetilde{\mathcal{S}}_{m_e,q_{e},\Lambda}$ satisfies the \sDEC\ so that radial profile $f_{2}$ matches the assumptions of Lemma~\ref{lem:glue}. We want the glued manifold to have charge $q_{e}$. To see that this is possible, we need to check Conditions~(1) through~(3) of Lemma~\ref{lem:glue}.  That is we will verify that
\begin{enumerate}
\item $f(b) < \widetilde{u}_{m_e,q_{e},\Lambda}(s_0-\delta)$ and $f'(b)\geq \widetilde{u}_{m_e,q_{e},\Lambda}'(s_0-\delta)$,
\item $\mathcal{M}_{\mathcal{S}}(\lbrace b\rbrace\times\mathbb{S}^n) \geq \frac{q^{2}}{f(b)^{n-1}}+\frac{2\Lambda f(b)^{n+1}}{n(n-1)(n+1)}$ with strict inequality if $q_{e}^{2}=q^{2}$, and
\item $\mathcal{M}_{\widetilde{\mathcal{S}}_{m_e,q_{e},\Lambda}}(\lbrace s_0-\delta\rbrace\times\mathbb{S}^n) > \frac{q_{e}^{2}}{\left(\widetilde{u}_{m_{e},q_{e},\Lambda}(s_{o}-\delta)\right)^{n-1}}+\frac{2\Lambda \left(\widetilde{u}_{m_e,q_{e},\Lambda}(s_0-\delta)\right)^{n+1}}{n(n-1)(n+1)} $,
\end{enumerate}
where $\mathcal{M}_{\mathcal{S}}$ and $\mathcal{M}_{\widetilde{\mathcal{S}}_{m_e,q_{e},\Lambda}}$ denote the generalized Hawking mass functionals withins the \ECR s $\mathcal{S}$ and $\widetilde{\mathcal{S}}_{m_e,q_{e},\Lambda}$, respectively.

We have ensured Condition~(1) by our choices of $s_0>0$, $0<\delta<s_0$ above. Condition~(2) follows directly from Assumption~\emph{(i)}, $q_{e}^{2}\leq q^{2}$, and $q_{e}^{2}<q^{2}$ if equality holds in Assumption~\emph{(i)}. Using Conditions~(1) and~(2), $q_{e}^{2}\leq q^{2}$, $\Lambda\leq0$, and Proposition~\ref{prop:sphereHawking}, in particular~\eqref{eq:Mp} with parameter $m=0$, we get
\begin{align*}
1-\left(\widetilde{u}_{m_{e},q_{e},\Lambda}'(s_{o}-\delta)\right)^{2}&\geq 1-f'(b)^{2}=1-p_{0,q,\Lambda}(f(b))+\frac{2\mathcal{M}_{\mathcal{S}}(\lbrace b\rbrace\times\mathbb{S}^{n})}{f(b)^{n-1}}\\
&>\frac{q_{e}^{2}}{f(b)^{2(n-1)}}+\frac{2\Lambda f(b)^{2}}{n(n-1)}\\
&\geq\frac{q_{e}^{2}}{\left(\widetilde{u}_{m_{e},q_{e},\Lambda}(s_{o}-\delta)\right)^{2(n-1)}}+\frac{2\Lambda \left(\widetilde{u}_{m_{e},q_{e},\Lambda}(s_{o}-\delta)\right)^{2}}{n(n-1)},
\end{align*}
so that Condition~(3) follows via the computation
\begin{align*}
&\mathcal{M}_{\widetilde{\mathcal{S}}_{m_e,q_{e},\Lambda}}(\lbrace s_0-\delta\rbrace\times\mathbb{S}^n)\\
&\quad\quad=\frac{\left(\widetilde{u}_{m_{e},q_{e},\Lambda}(s_{o}-\delta)\right)^{n-1}}{2}\times\\
&\quad\quad\quad\quad\left(1-\left(\widetilde{u}_{m_{e},q_{e},\Lambda}'(s_{o}-\delta)\right)^{2}+\frac{q_{e}^{2}}{\left(\widetilde{u}_{m_{e},q_{e},\Lambda}(s_{o}-\delta)\right)^{2(n-1)}}-\frac{2\Lambda \left(\widetilde{u}_{m_{e},q_{e},\Lambda}(s_{o}-\delta)\right)^{2}}{n(n+1)}\right)\\
&\quad\quad>\frac{q_{e}^{2}}{ \left(\widetilde{u}_{m_{e},q_{e},\Lambda}(s_{o}-\delta)\right)^{n-1}}+\frac{2\Lambda\left(\widetilde{u}_{m_{e},q_{e},\Lambda}(s_{o}-\delta)\right)^{n+1}}{n(n-1)(n+1)}.
\end{align*}

Having verified all assumptions of the Gluing Lemma~\ref{lem:glue}, we can now conclude that $\mathcal{S}$ and $\widetilde{\mathcal{S}}_{m_e,q_{e},\Lambda}$ can be smoothly glued to an \ECR\ $\widetilde{\mathcal{S}}$ as desired (after an appropriate translation of $[s_0-\delta,s_{1}]$ and re-attaching $(s_{1},\infty)\times\mathbb{S}^{n}$ after the application of the Gluing Lemma). If, in addition, $f'>0$ on $[a,b]$, then the Gluing Lemma~\ref{lem:glue}, together with the fact that $\widetilde{u}_{m_{e},q_{e},\Lambda}'>0$ on $[s_0-\delta,\infty)$, gives $\widetilde{f}'>0$ for the radial profile of $\widetilde{\mathcal{S}}$ implying the remaining claim about strict mean convexity of the coordinate spheres via Lemma~\ref{lem:meancurvature}.
\end{proof}

In previous work on charged manifolds~\cite{impor}, the corresponding Proposition 4.1 allowing to glue a given rotationally symmetric \ECR\ to a Reissner--Nordstr\"om manifold made stronger assumptions than our Proposition~\ref{prop:gluingRN}. More specifically, it required that $q_{e}=q$, $f(b)>\vert q\vert$, and $\mathcal{M}(\lbrace b\rbrace\times\mathbb{S}^{2})>\vert q\vert$. Together, these assumptions imply in particular that $(M_{m_{e},q,\Lambda},\gamma_{m_{e},q,\Lambda})$ is sub-extremal a priori. Allowing to glue to arbitrary \RN s subject to $m_{e}\geq m$, $q_{e}^{2}\leq q^{2}$, and $m_{e}>m$ or $q_{e}^{2}<q^{2}$ will allow us to deduce rather than enforce that the extensions to be constructed in Section~\ref{sec:main} are sub-extremal.

\section{Main Results}\label{sec:main}
We will now combine our results on the construction of collar extensions starting at minimal hypersurfaces with prescribed Riemannian metrics in Section~\ref{sec:collar} with the gluing results obtained in Section~\ref{sec:gluing} to obtain \ECR s with prescribed boundary geometry, satisfying the \DEC, and with total mass arbitrarily close to the optimal value of the Riemannian Penrose Inequality~\eqref{eq:RPI}.

\begin{theorem}[Existence of Extensions]\label{thm:main}
Let $n\geq2$, let $g_o$ be a smooth Riemannian metric on $\mathbb{S}^{n}$ and let $r_o$ denote the volume radius of $g_o$, $\abs{\mathbb{S}^n}_{g_0}\asdefined\omega_n r_o^n$, where $\omega_{n}=\vert \mathbb{S}^{n}\vert_{g_{*}}$ denotes the volume of the unit round sphere $(\mathbb{S}^{n},g_{*})$. Assume that one of the following holds
\begin{enumerate}
\item $n\geq2$, $\lambda_{1}(-\triangle_{g_o}+\tfrac{1}{2}\scal(g_{o}))>0$, and $(\mathbb{S}^{n},g_{o})$ is conformal to $(\mathbb{S}^{n},g_{*})$,\\ where $\lambda_{1}(-\triangle_{g_o}+\tfrac{1}{2}\scal(g_{o}))$ denotes the first eigenvalue of $-\triangle_{g_o}+\tfrac{1}{2}\scal(g_o)$, or
\item $n=2$ and $K(g_o)>-\tau$ on $\mathbb{S}^{2}$ for some $\tau>0$, or
\item $2\leq n\leq3$  and $\scal(g_o)>0$ on $\mathbb{S}^{n}$, where $\scal(g_o)$ denotes the scalar curvature of $g_o$, or
\item $n>3$, $\scal(g_o)>0$ on $\mathbb{S}^{n}$, and $(\mathbb{S}^n,g_o)$ has pointwise $\nicefrac{1}{4}$-pinched sectional curvature.
\end{enumerate}
Here, $\scal(g_{o})$ and $K(g_{o})$ denote the scalar and Gaussian curvatures of $g_{o}$, respectively. Then there exists a constant $\kappa \geq 0$ with $\kappa<\lambda_{1}(-\triangle_{g_o}+\tfrac{1}{2}\scal(g_{o}))$ in Case~(1), $\kappa\leq\tau$ in Case~(2), and $2\kappa<\scal(g_o)$ in Cases~(3) and~(4) such that, for any constants $\Lambda\leq0$ and $q\in\R$ satisfying
\begin{align}\tag{C}\label{eq:qcondthm}
\begin{split}
\left\{
\begin{array}{lll}
\frac{n(n-1)q^2}{r_{o}^{2n}}&<2(\kappa-\Lambda)&\text{ in Cases~(1), (3), and (4)}\\
\frac{q^2}{r_{o}^{4}}&<-(\kappa+\Lambda)&\text{ in Case~(2)},
\end{array}
\right.
\end{split}
\end{align}
and any constant $m\in\mathbb{R}$ satisfying
\begin{align}\label{eq:mcondthm}
m>m_o&\definedas\frac{r_o^{n-1}}{2}\left(1+\frac{q^2}{r_o^{2(n-1)}} -\frac{2\Lambda r_o^{2}}{n(n+1)}\right),
\end{align}
there exists an $(n+1)$-dimensional \ECR\ $(M,\gamma,E,\Lambda)$ which is geodesically complete up to its strictly outward minimizing minimal boundary $\partial M$, satisfies the Dominant Energy Condition
\begin{align*}
\scal(\gamma)\geq 2\Lambda+n(n-1)\abs{E}_\gamma^2,
\end{align*}
and has the following properties
\begin{enumerate}[label=(\roman*)]
\item $(\partial M,g)$ is isometric to $(\mathbb{S}^n,g_o)$, where $g$ is the metric induced on $\partial M$ by $\gamma$, and
\item there exists a compact subset $C\subset M$ and a radius $r_{C}>r_{+}$ such that $(M\setminus C,\gamma,E,\Lambda)$ is isometric to the subset $((r_{C},\infty)\times\mathbb{S}^n,\gamma_{m,q,\Lambda},E_{m,q,\Lambda},\Lambda)$ of the sub-extremal \RN\ of mass $m$, charge $q$, and cosmological constant $\Lambda$, with $r_{+}=r_{+}(m,q,\Lambda)$ denoting the volume radius of its horizon inner boundary.
\end{enumerate}
In particular, the (conjectured) Riemannian Penrose Inequality~\eqref{eq:RPI} in $(M,\gamma,E,\Lambda)$ reduces to $m_{o}\leq m$. In other words, the total mass $m$ can be chosen arbitrarily close to the optimal value in the (conjectured) Riemannian Penrose Inequality~\eqref{eq:RPI}.
.\end{theorem}

\begin{proof}
From Lemma~\ref{le1.2} (in Case~(1) with $n=2$, appealing to the Uniformization Theorem), Lemma~\ref{lem:uniform} (in Case~(1) with $n>2$), Lemma~\ref{path2} (in Case~(2) or Case~(3) for $n=2$), Corollary~\ref{coro:paths3} (in Case~(3)), and Corollary~\ref{coro:paths} (in Case~(4)), respectively, we obtain smooth paths of metrics $\{g(t)\}_{t\in[0,1]}$ on $\mathbb{S}^{n}$ satisfying Conditions \hyperref[page:conditions]{(E1-E3)} and 
\begin{enumerate}
\item[(E4)]\label{E4}  $\lambda_{1}(-\triangle_{g(t)}+\tfrac{1}{2}\scal(g(t)))>0$ for all $t\in[0,1]$ in Case~(1), $K(g(t))>-\tau$ on $\mathbb{S}^{2}$ for all $t\in[0,1]$ in Case~(2), and $\scal(g(t))>0$ on $\mathbb{S}^{n}$ for all $t\in[0,1]$ in Cases~(3) and (4).
\end{enumerate}
Next, we set $\kappa\geq0$ to be any constant satisfying
\begin{align*}
0\leq\kappa&<\left\{
\begin{array}{ll}
\min_{t\in[0,1]}\lambda_{1}(-\triangle_{g(t)}+\tfrac{1}{2}\scal(g(t)))&\text{ in Case~(1)}\\
\tfrac{1}{2}\min_{t\in[0,1]}\min_{\mathbb{S}^{n}}\scal(g(t))&\text{ in Cases~(3) and (4), and}
\end{array}
\right.\\
-\tau\leq-\kappa&<\quad\min_{t\in[0,1]}\min_{\mathbb{S}^{2}}K(g(t))\quad\quad\quad\quad\quad\quad\quad\;\text{ in Case~(2)}.
\end{align*}
Hence we can replace Condition~\hyperref[E4]{(E4)} above by 
\begin{enumerate}\label{E4'}
\item[(E4)'] $\lambda_{1}(-\triangle_{g(t)}+\tfrac{1}{2}\scal(g(t)))>\kappa$ for all $t\in[0,1]$ in Case~(1), $K(g(t))>-\kappa$ on $\mathbb{S}^{2}$ for all $t\in[0,1]$ in Case~(2), and $\scal(g(t))>\kappa$ on $\mathbb{S}^{n}$ for all $t\in[0,1]$ in Cases~(3) and~(4).
\end{enumerate}

For any $\Lambda\leq0$, $q\in\R$ subject to Condition~\eqref{eq:qcondthm} with respect to this choice of $\kappa$, Lemma~\ref{lem:subapriori} ensures the a priori sub-extremality condition $h_{q,\Lambda}(r_{o})>0$. Of course, by Conditions~\hyperref[page:conditions]{(E1)} and~\hyperref[page:conditions]{(E3)}, $(\mathbb{S}^{n},g(1))$ must be isometric to $(\mathbb{S}^{n},r_{o}^{2}g_{*})$ via an isometry $\Psi\colon(\mathbb{S}^{n},g(1))\to(\mathbb{S}^{n},r_{o}^{2}g_{*})$ (see for example Corollary 10 in~\cite{Petersen}). As it will be simpler to deal with $r_{o}^{2}g_{*}$ than with an abstractly given round metric, we change the path $\lbrace g(t)\rbrace_{t\in[0,1]}$ to $\lbrace\widetilde{g}(t)\definedas\Psi_{*}g(t)\rbrace_{t\in[0,1]}$ and note that the new smooth path $\lbrace\widetilde{g}(t)\rbrace_{t\in[0,1]}$ automatically satisfies Conditions~\hyperref[page:conditions]{(E2)}, \hyperref[page:conditions]{(E3)} and \hyperref[E4']{(E4)'} as well as the alternative condition
\begin{enumerate}
\item[(E1)']\label{E1'} $(\mathbb{S}^{n},\widetilde{g}(0))$ is isometric to $(\mathbb{S}^{n},g_o)$ and $\widetilde{g}(1)=r_{o}^{2}g_{*}$.
\end{enumerate}
Having said this, we will now drop the explicit mention of the isometry $\Psi$ and keep calling the path $\lbrace g(t)\rbrace_{t\in[0,1]}$. Lemma~\ref{lem:collar} then tells us that for this choice of $\kappa$ and for $q$ and $\Lambda$ subject to Condition~\eqref{eq:qcondthm}, there exist a constant $A_{o}>0$ and \ECR s 
\begin{align*}
\mathcal{S}_{\varepsilon,A}=([0,1]\times\mathbb{S}^{n},\gamma_{\varepsilon,A}=(Au(t,\cdot))^{2}dt^{2}+(1+\varepsilon t^{2})g(t),E_{\varepsilon,A}=\frac{q}{Au(t,\cdot)r_{o}^{n}(1+\varepsilon t^{2})^{\frac{n}{2}}}\partial_{t},\Lambda)
\end{align*}
parametrized by $\varepsilon\in(0,1]$ and $A\geq A_{o}$, with $u(t,\cdot)>0$ the positive $L^{2}(g(t))$-normalized first eigenfunction of $-\triangle_{g(t)}+\tfrac{1}{2}\scal(g(t))$ in Case~(1) (see Remark~\ref{rem:normalized}) and $u(t,\cdot)=1$ in Cases~(2) through (4) which satisfy the \sDEC~\eqref{eq:sDEC}, have minimal inner boundary $\partial M= \lbrace0\rbrace\times\mathbb{S}^{n}$ and strictly mean convex coordinate spheres $\{t \} \times \mathbb{S}^n$ for $t\in (0,1]$ with respect to the outward unit normal $\nu=\frac{1}{Au(t,\cdot)}\partial_t$. 

Before we can attempt at applying Proposition~\ref{prop:gluingRN} to glue the collar extensions $\mathcal{S}_{\varepsilon,A}$ to a \RN, we will need to make a change of coordinates on the collar extension, bringing it to the form $\gamma_{\varepsilon,A}=ds^{2}+f_{\varepsilon,A}(s)^{2}g_{*}$, $E_{\varepsilon,A}=\frac{q}{f_{\varepsilon,A}(s)^{n}}\partial_{s}$ in its rotationally symmetric part. Recalling that $g(t)=g(1)=r_{o}^{2}g_{*}$ for $t\in[\theta,1]$ and noting that this gives $u(t,\cdot)=1$ for $t\in[\theta,1]$, this can be achieved by setting $s=At$ and $f_{\varepsilon,A}(s)\definedas \sqrt{1+\frac{\varepsilon s^{2}}{A^{2}}}\,r_{o}$ for $s\in[\theta A,A]$. Hence, $\mathcal{S}_{\varepsilon,A}$ restricted to $[\theta,1]\times\mathbb{S}^{n}$ is isometric to the rotationally symmetric \ECR\ 
\begin{align*}
\widetilde{\mathcal{S}}_{\varepsilon,A}\definedas([\theta A,A]\times\mathbb{S}^{n},\gamma_{\varepsilon,A}=ds^{2}+f_{\varepsilon,A}(s)^{2}g_{*},E_{\varepsilon,A}=\frac{q}{f_{\varepsilon,A}(s)^{n}}\partial_{s},\Lambda)
\end{align*}
which satisfies the \sDEC~\eqref{eq:sDEC} and has strictly mean convex coordinate spheres with respect to the outward unit normal $\nu=\partial_s$. In view of Proposition~\ref{prop:gluingRN}, we will verify that
\begin{enumerate}[label=(\alph*)]
\item $m_{*}(\varepsilon,A)\definedas\mathcal{M}_{\widetilde{\mathcal{S}}_{\varepsilon,A}}(\lbrace A\rbrace\times\mathbb{S}^n)>\frac{q^{2}}{f_{\varepsilon,A}(A)^{n-1}}+\frac{2\Lambda f_{\varepsilon,A}(A)^{n+1}}{n(n-1)(n+1)}$, and
\item $f_{\varepsilon,A}'(A)>0$,
\end{enumerate} 
where $\mathcal{M}_{\widetilde{\mathcal{S}}_{\varepsilon,A}}$ denotes the generalized Hawking mass functional in $\widetilde{\mathcal{S}}_{\varepsilon,A}$. 

Condition~(b) follows immediately from Lemma~\ref{lem:meancurvature} and the fact that $\lbrace A\rbrace\times\mathbb{S}^{n}$ has strictly positive mean curvature in $\widetilde{\mathcal{S}}_{\varepsilon,A}$ by construction for all $\varepsilon\in(0,1]$ and all $A\geq A_{o}$. It remains to prove Condition~(a). Using Lemma~\ref{lem:Hawking*}, we know that $m_{*}(\varepsilon,A)> m_{o}$ which implies
\begin{align*}
m_{*}(\varepsilon,A)&> m_{o}=\frac{r_{o}^{n-1} p_{0,q,\Lambda}(r_{o})}{2}=\frac{r_{o}^{n-1}}{2}\left(p_{0,q,\Lambda}(r_{o})-1+\frac{h_{q,\Lambda}(r_{o})}{r_{o}^{2(n-1)}}+\frac{q^{2}}{r_{o}^{2(n-1)}}+\frac{2\Lambda r_{o}^{2}}{n(n-1)}\right)\\
&>\frac{q^2}{r_o^{n-1}}+\frac{2\Lambda r_o^{n+1}}{n(n-1)(n+1)}\geq\frac{q^2}{f_{\varepsilon,A}(A)^{n-1}}+\frac{2\Lambda f_{\varepsilon,A}(A)^{n+1}}{n(n-1)(n+1)}
\end{align*}
for all $0<\varepsilon\leq1$  and all $A\geq A_{o}$. This asserts Condition~(a).

Hence by Proposition~\ref{prop:gluingRN} and by construction, we know that $\widetilde{\mathcal{S}}_{\varepsilon,A}$ and thus $\mathcal{S}_{\varepsilon,A}$ can smoothly be glued to any \RN\ of mass $m>m_{*}(\varepsilon,A)>m_{o}$, charge $q$, and cosmological constant~$\Lambda$ for any $0<\varepsilon\leq1$ and $A\geq A_{o}$, such that the glued \ECR\ $(M,\gamma,E,\Lambda)$ satisfies the \DEC\ and Claims~\emph{(i)} and~\emph{(ii)}. Moreover, again by Proposition~\ref{prop:gluingRN} and by construction as well as by Lemma~\ref{lem:meancurvature}, we know that $(M,\gamma)$ is foliated by mean convex coordinate spheres $\lbrace \cdot\rbrace\times\mathbb{S}^{n}$ which are all strictly mean convex except $\partial M$ which is minimal. Hence $\partial M$ is a strictly outward minimizing minimal hypersurface in $(M,\gamma)$ by Theorem~\ref{thm:outwardminfoliation}. This proves all claims in case $m>m_{*}(\varepsilon,A)$ for some $0<\varepsilon\leq1$ and $A\geq A_{o}$.

As the theorem makes these claims for any $m>m_{o}$, we again apply Lemma~\ref{lem:Hawking*} to see that $m_{*}(\varepsilon,A)\searrow m_{o}$ as $\varepsilon\searrow0$. Thus $m_{o}<m_{*}(\varepsilon,A)<m$ for all suitably small $0<\varepsilon<\varepsilon_{1}$ with threshold $\varepsilon_{1}=\varepsilon_{1}(n,m,r_{o},q,\Lambda)$ and all $A\geq A_{o}$ which proves the claim. 

Finally, by definition of the generalized Hawking mass and by $m_{o}=\mathcal{M}_{\mathcal{S}_{\varepsilon,A}}(\lbrace0\rbrace\times\mathbb{S}^{n})$ for all $\varepsilon\in(0,\varepsilon_{1})$, $A\geq A_{o}$, we get the desired interpretation in view of the (conjectured) Riemannian Penrose Inequality.
\end{proof}

\begin{remark}\label{rem:extremality}
We point out that we do not need an a priori sub-extremality assumption; sub-extremality indeed follows from our assumptions via Lemma~\ref{lem:subapriori}. In this sense, our result somewhat simplifies that of~\cite{impor} where it is assumed for $n=2$ and $\Lambda=0$ that $q^{2}<r_{o}^{2}$ and $\kappa>\tfrac{q^{2}}{r_{o}^{4}}$; arguing as in the proof of Lemma~\ref{lem:subapriori}, we can indeed conclude that $q^{2}<r_{o}^{2}$ must hold provided that $\kappa>\tfrac{q^{2}}{r_{o}^{4}}$, which coincides with our Condition~\eqref{eq:qcondthm}.
\end{remark}

\begin{remark}\label{metrirema}
Recall from Section~\ref{sec:collar} that the assumption that $(\mathbb{S}^{n},g_{o})$ has pointwise $\nicefrac{1}{4}$-pinched sectional curvature for $n\geq4$ is a technical assumption ensuring that positivity of scalar curvature is preserved along a certain flow of metrics on $\mathbb{S}^{n}$. In fact, Theorem~\ref{thm:main} will apply to any metric $g_{o}$ on $\mathbb{S}^{n}$ for which one can find a path of metrics $\lbrace g(t)\rbrace_{t\in[0,1]}$ on $\mathbb{S}^{n}$, $n\geq4$, with $g(0)=g_{o}$, satisfying the Conditions \hyperref[page:conditions]{(E1-E3)}, and having $\scal(g(t))>0$ for all $t\in[0,1]$. For example, one can also use the path of metrics constructed by Cabrera~Pacheco and Miao~\cite{sec} using extrinsic geometric flows, where it is assumed that $(\mathbb{S}^{n},g_o)$, $n\geq4$, with $\scal(g_{o})>0$ isometrically embeds into $(n+1)$-dimensional Euclidean space as a star-shaped hypersurface.
\end{remark}

\section{Application to an Ad Hoc Generalization of Bartnik Mass}\label{sec:discussion}
In Theorem~\ref{thm:main}, we have seen that many metrics $g_{o}$ on $\mathbb{S}^{n}$ can arise as the metric induced on the minimal, strictly outward minimizing inner boundary of an \ECR\ satisfying the \DEC~\eqref{eq:DEC} which eventually coincides with a sub-extremal \RN\ of prescribed mass $m$, charge $q$, and cosmological constant $\Lambda\leq0$. We have seen that the possible choices of the cosmological constant $\Lambda$ are constrained by the geometry of $g_{o}$ (and a choice of path of metrics on $\mathbb{S}^{n}$). The range of possible charges $q$ was constrained by both the choice of cosmological constant $\Lambda$ and by the geometry of $g_{o}$ (and the chosen path). Finally, the mass $m$ can be prescribed freely, as long as it is larger than the generalized Hawking mass of $g_{o}$ which depends on the chosen charge $q$ and cosmological constant $\Lambda$ --- in consistency with the (conjectured) Riemannian Penrose Inequality (with charge)~\eqref{eq:RPI}, see Conjecture~\ref{conj:RPI}.

To see this, we have constructed so-called \emph{extensions}, namely \ECR s $(M,\gamma,E,\Lambda)$ satisfying the \DEC~\eqref{eq:DEC}, isometric to a piece of a suitable sub-extremal \RN\ $(M_{m,q,\Lambda},\gamma_{m,q,\Lambda},E_{m,q,\Lambda},\Lambda)$ of mass $m$, charge $q$, and cosmological constant $\Lambda$ outside some compact set, and such that $\partial M$ is minimal and strictly outward minimizing in $(M,\gamma)$, isometric to $(\mathbb{S}^{n},g_{o})$ with respect to the induced metric $g$, and has quasi-local charge $\mathcal{Q}(\partial M)=q$. We have discussed that the asymptotic behavior of $(M,\gamma,E,\Lambda)$ is precisely that of $(M_{m,q,\Lambda},\gamma_{m,q,\Lambda},E_{m,q,\Lambda},\Lambda)$, and in particular the asymptotic mass equals $m$ and the asymptotic charge equals $q$, see Remarks~\ref{rem:mass}, \ref{rem:charge}. If $\Lambda=0$, this means that $(M,\gamma,E,\Lambda)$ is asymptotically Euclidean with ADM mass $m$, while if $\Lambda<0$, $(M,\gamma,E,\Lambda)$ is asymptotically hyperbolic with asymptotic hyperbolic radius $\sqrt{-\tfrac{n(n+1)}{2\Lambda}}$ and has hyperbolic mass $m$, see Remark~\ref{rem:asymptotics}.

The constructed extensions hence fall in the class $\mathcal{A}(\mathbb{S}^{n},g_{o},H_{o},q,\Lambda)$ of \emph{(generalized) admissible extensions} of given \emph{(generalized) Bartnik data $(\mathbb{S}^{n},g_{o},H_{o},q,\Lambda)$} consisting of a metric $g_{o}$ on $\mathbb{S}^{n}$, a smooth, non-negative function $H_{o}\colon\mathbb{S}^{n}\to\R^{+}_{0}$, and parameters $q\in\R$ and $\Lambda\leq0$. Generalizing the standard definition by Bartnik~\cite{Bartnik-89} and the separate asymptotically hyperbolic~\cite{hyper} and electrically charged~\cite{impor} generalizations, we define this class $\mathcal{A}(\mathbb{S}^{n},g_{o},H_{o},q,\Lambda)$ of (generalized) admissible extensions of given Bartnik data $(\mathbb{S}^{n},g_{o},H_{o},q,\Lambda)$ as
\begin{align*}
\begin{split}
&\mathcal{A}(\mathbb{S}^{n},g_{o},H_{o},q,\Lambda)\definedas\\[-0.2ex]
&\quad\lbrace \text{\ECR s }(M,\gamma,E,\Lambda)\text{ which}\\[-0.2ex]
&\quad\text{ have strictly outward minimizing boundary }(\partial M,g)\text{ isometric to }(\mathbb{S}^{n},g_{o})\\
&\quad\text{ and with mean curvature }H=H_{o},\\[-0.2ex]
&\quad \text{ satisfy the Dominant Energy Condition }\operatorname{R}(\gamma)\geq 2\Lambda+n(n-1)\abs{E}_\gamma^2\text{ on }M,\\[-0.2ex]
&\quad\text{ are asymptotically Euclidean if }\Lambda=0,\\[-1.2ex]
&\quad\text{ are asymptotically hyperbolic with asymptotic hyperbolic radius }\sqrt{-\tfrac{n(n+1)}{2\Lambda}}\text{ if }\Lambda<0,\\[-0.2ex]
&\quad\text{ and have total charge }q\rbrace.
\end{split}
\end{align*}

Before we proceed to defining the (generalized) Bartnik mass of given Bartnik data, we should make a couple of comments. First of all, we could have given an alternative definition of generalized admissible extensions, replacing the condition that $\partial M$ should be strictly outward minimizing by the condition that the extension must not contain any minimal surfaces enclosing the boundary, except possibly the boundary itself. Both definitions are standard in the asymptotically Euclidean $n+1=3$-dimensional Riemannian case, i.e., without electric field. In this scenario, Jauregui~\cite{Jauregui1} and McCormick~\cite{McC} have analyzed when the two definitions and the Bartnik masses defined on their bases (see below) coincide. Both conditions are designed to rule out that the (generalized) Bartnik data ``hide'' behind a black hole horizon, i.e., behind a minimal surface, in an admissible extension. Note that the extensions we have constructed in this paper (see Theorem~\ref{thm:main}) naturally satisfy both conditions. For the strictly outward minimizing condition, see the statement of Theorem~\ref{thm:main}; the no minimal surface condition is satisfied because the asserted foliation by strictly mean convex hypersurfaces away from the boundary and the elliptic maximum principle exclude the existence of such minimal surfaces, see for example Chapter 4 in~\cite{DLee}.

Second, it is of course not clear a priori whether $\mathcal{A}(\mathbb{S}^{n},g_{o},H_{o},q,\Lambda)\neq\emptyset$ for given Bartnik data $(\mathbb{S}^{n},g_{o},H_{o},q,\Lambda)$. Originally, Bartnik only considered Bartnik data which were a priori embedded in a geodesically complete, asymptotically Euclidean Riemannian manifold satisfying the \DEC\ and an appropriate condition ensuring that the Bartnik data do not hide behind a horizon. If one does so, of course the set of admissible extensions is automatically non-empty. Here, however, we do not put such an assumption so in principle the (generalized) admissible class for given (generalized) Bartnik data can be empty. It is a consequence of the special case of (generalized) Bartnik data we consider here that they do possess (generalized) admissible extensions, see below.

Third, it is worth pointing out that we have not made an assumption nor a conclusion about physicality of (generalized) admissible extensions. In the asymptotically Euclidean $n+1=3$-dimensional Riemannian case (and with the strictly outward minimizing condition in the definition), physicality in the form of positive ADM mass follows from Huisken and Ilmanen’s proof of the Riemannian Penrose Inequality~\cite{RPI} via Inverse Mean Curvature Flow. As this approach doesn't lend itself to generalizations to higher dimensions due to the lack of a Gauss--Bonnet Theorem (see Section~\ref{sec:Hawkingcollar}), such a physicality assertion cannot (yet?) be made in higher dimensions (except when $H_{o}=0$, see below). Similarly, physicality in the sense of positive hyperbolic mass has not been shown in the asymptotically hyperbolic case (except in certain cases with $H_{o}=0$, see below). In an electrically charged scenario such as the one considered in this paper, physicality should be interpreted as not only asking that the total mass be positive but a fortiori that the total mass $m$, total charge $q$, and cosmological constant $\Lambda\leq0$ are subject to a sub-extremality condition. One possibility could be to request that the \RN\ of mass $m$, charge $q$, and cosmological constant $\Lambda$ be sub-extremal, see Definition~\ref{def:extremality}. Again, this is not known in general for admissible extensions. 

In this paper we only consider the case of \emph{minimal (generalized) Bartnik data}, i.e., of Bartnik data $(\mathbb{S}^{n},g_{o},H_{o},q,\Lambda)$ with $H_{o}=0$. In this case, physicality of (generalized) admissible extensions (with respect to the definition given here) in the sense of positivity of total mass $m$ is guaranteed by the Riemannian Penrose Inequality (with charge)~\eqref{eq:RPI} in all cases in which it has been proven, see Remark~\ref{rem:RPI} which can be seen as follows: Recall that, for minimal hypersurfaces, the left hand side of the (conjectured) Riemannian Penrose Inequality (with charge)~\eqref{eq:RPI} coincides with the generalized Hawking mass $\mathcal{M}(\partial M)$ of the boundary, see Definition~\ref{def:Hawking}. From Proposition~\ref{prop:Hawkingre}, we know that $\mathcal{M}(\partial M)=\tfrac{r(\partial M)^{n-1}}{2}\,p_{0,\mathcal{Q}(\partial M),\Lambda}(r(\partial M))$, where $r(\partial M)$ denotes the volume radius of $\partial M$ with respect to the induced metric and $\mathcal{Q}(\partial M)$ denotes the quasi-local charge of $\partial M$, see Definition~\ref{def:QLcharge}. Hence $\mathcal{M}(\partial M)>0$ is manifest for minimal boundaries. In a given admissible extension $(M,\gamma,E,\Lambda)$ of minimal Bartnik data $(\mathbb{S}^{n},g_{o},H_{o},q,\Lambda)$, we in fact have $\mathcal{Q}(\partial M)=q$ as $E$ is divergence-free. Hence, the total mass $m$ of any (generalized) admissible extension of minimal (generalized) Bartnik data is necessarily positive.

Let us now turn our attention to the question of sub-extremality of (generalized) admissible extensions $(M,\gamma,E,\Lambda)$ of minimal (generalized) Bartnik data $(\mathbb{S}^{n},g_{o},H_{o},q,\Lambda)$. By Proposition~\ref{prop:QLsub}, the quasi-local sub-extremality condition $h_{\mathcal{Q}(\partial M),\Lambda}(r(\partial M))>0$ ensures that the \RN\ of mass $\mathcal{M}(\partial M)$, charge $\mathcal{Q}(\partial M)$, and cosmological constant $\Lambda\leq0$ is sub-extremal. Lemma~\ref{lem:largerm} tells us that the \RN\ of mass $m\geq\mathcal{M}(\partial M)$, charge $\mathcal{Q}(\partial M)$, and cosmological constant $\Lambda$ must also be sub-extremal. Hence the Riemannian Penrose Inequality (with charge)~\eqref{eq:RPI} --- in the cases in which it is proven --- together with the quasi-local sub-extremality condition $h_{\mathcal{Q}(\partial M),\Lambda}(r(\partial M))>0$ in fact implies physicality of (generalized) admissible extensions of minimal (generalized) Bartnik data in the sense that the \RN\ of the same total mass, total charge, and cosmological constant is sub-extremal. This proves the following proposition.
 
 \begin{proposition}[A Priori Sub-Extremality for Minimal Bartnik Data]\label{prop:Bartnikapriori}
 Let $(\mathbb{S}^{n},g_{o},H_{o}=0,q,\Lambda)$ be minimal (generalized) Bartnik data, i.e., let $g_{o}$ be a smooth Riemannian metric on $\mathbb{S}^{n}$, and let $q\in\R$ and $\Lambda\leq0$ be constants, and let $r_{o}$ denote the volume radius of $g_{o}$. Assume that $h_{q,\Lambda}(r_{o})>0$. Let $(M,\gamma,E,\Lambda)$ be a (generalized) admissible extension of $(\mathbb{S}^{n},g_{o},H_{o}=0,q,\Lambda)$, i.e., an \ECR\ satisfying the \DEC~\eqref{eq:DEC}, having minimal, strictly outward minimizing boundary $\partial M$ with induced metric $g$ such that $(\partial M,g)$ is isometric to $(\mathbb{S}^{n},g_{o})$, and being asymptotically Euclidean of ADM mass $m$ and total charge $q$ if $\Lambda=0$ or asymptotically hyperbolic with asymptotic hyperbolic radius $\sqrt{-\tfrac{n(n+1)}{2\Lambda}}$ of hyperbolic mass $m$ and total charge $q$ if $\Lambda<0$. Then the \RN\ of mass $m$, charge $q$, and cosmological constant $\Lambda$ is sub-extremal provided that the (conjectured) Riemannian Penrose Inequality (with charge)~\eqref{eq:RPI} applies to a class of \ECR s containing $(M,\gamma,E,\Lambda)$.
 \end{proposition}
\vfill\newpage
Let us now turn our attention to the \emph{(generalized) Bartnik mass} functional $\mathcal{B}$. For given (generalized) Bartnik data $(\mathbb{S}^{n},g_{o},H_{o},q,\Lambda)$, we set
\begin{align}\label{def:Bartnik}
\mathcal{B}(\mathbb{S}^{n},g_{o},H_{o},q,\Lambda)&\definedas \inf\lbrace m(M,\gamma,E,\Lambda)\,\vert\,(M,\gamma,E,\Lambda)\in\mathcal{A}(\mathbb{S}^{n},g_{o},H_{o},q,\Lambda)\rbrace,
\end{align}
where $m(M,\gamma,E,\Lambda)$ denotes the ADM mass of $(M,\gamma)$ if $\Lambda=0$ and the hyperbolic mass of $(M,\gamma)$ if $\Lambda<0$. Note that we deviate from the original definition of Bartnik mass and admissible extensions here which did not consider manifolds with boundary but geodesically complete manifolds containing a given domain $\Omega$ with compact boundary $\partial\Omega$. This has to do with regularity issues across $\partial\Omega$ in potential mass-minimizing admissible extensions and has become rather standard in the literature. See the survey by Anderson~\cite{survey2} and the article by Bartnik~\cite{BartnikTsingHua} for more information.

As the generalized Hawking mass of given (generalized) Bartnik data $(\mathbb{S}^{n},g_{o},H_{o},q,\Lambda)$ in a given (generalized) admissible extension only depends on the (generalized) Bartnik data themselves, we can abuse notation and denote the generalized Hawking mass of (generalized) Bartnik data as a functional of those data as well, i.e., as $\mathcal{M}(\mathbb{S}^{n},g_{o},H_{o},q,\Lambda)$. Then the (conjectured) Riemannian Penrose Inequality (with charge)~\eqref{eq:RPI} can be rewritten as
\begin{align}\label{eq:RPIrewritten}
\mathcal{M}(\mathbb{S}^{n},g_{o},H_{o}=0,q,\Lambda)&\leq\mathcal{B}(\mathbb{S}^{n},g_{o},H_{o}=0,q,\Lambda)
\end{align}
for all minimal (generalized) Bartnik data $(\mathbb{S}^{n},g_{o},H_{o},q,\Lambda)$ with $\mathcal{A}(\mathbb{S}^{n},g_{o},H_{o},q,\Lambda)\neq\emptyset$.\\

Complementing this, Theorem~\ref{thm:main} can be rewritten as giving the following upper estimate on the generalized Bartnik mass.
\begin{theorem}[Generalized Bartnik Mass Estimate]\label{thm:Bartnikestimate}
Let $n\geq2$, let $g_o$ be a smooth Riemannian metric on $\mathbb{S}^{n}$ and let $r_o$ denote the volume radius of $g_o$, $\abs{\mathbb{S}^n}_{g_0}\asdefined\omega_n r_o^n$, where $\omega_{n}=\vert \mathbb{S}^{n}\vert_{g_{*}}$ denotes the volume of the unit round sphere $(\mathbb{S}^{n},g_{*})$. Assume that one of the following holds
\begin{enumerate}
\item $n\geq2$, $\lambda_{1}(-\triangle_{g_o}+\tfrac{1}{2}\scal(g_{o}))>0$, and $(\mathbb{S}^{n},g_{o})$ is conformal to $(\mathbb{S}^{n},g_{*})$,\\ where $\lambda_{1}(-\triangle_{g_o}+\tfrac{1}{2}\scal(g_{o}))$ denotes the first eigenvalue of $-\triangle_{g_o}+\tfrac{1}{2}\scal(g_o)$, or
\item $n=2$ and $K(g_o)>-\tau$ on $\mathbb{S}^{2}$ for some $\tau>0$ or
\item $2\leq n\leq3$ and $\scal(g_o)>0$ on $\mathbb{S}^{n}$, where $\scal(g_o)$ denotes the scalar curvature of $g_o$, or
\item $n>3$, $\scal(g_o)>0$ on $\mathbb{S}^{n}$, and $(\mathbb{S}^n,g_o)$ has pointwise $\nicefrac{1}{4}$-pinched sectional curvature.
\end{enumerate}
Here, $\scal(g_{o})$ and $K(g_{o})$ denote the scalar and Gaussian curvatures of $g_{o}$, respectively. Then there exists a constant $\kappa \geq 0 $ with $\kappa<\lambda_{1}(-\triangle_{g_o}+\tfrac{1}{2}\scal(g_{o}))$ in Case~(1), $\kappa\leq\tau$ in Case~(2), and $2\kappa<\scal(g_o)$ in Cases~(3) and~(4) such that, for any constants $\Lambda\leq0$, $q\in\R$ satisfying
\begin{align}\tag{C}\label{eq:qcondthmB}
\begin{split}
\left\{
\begin{array}{lll}
\frac{n(n-1)q^2}{r_{o}^{2n}}&<2(\kappa-\Lambda)&\text{ in Cases~(1), (3), and (4)}\\
\frac{q^2}{r_{o}^{4}}&<-(\kappa+\Lambda)&\text{ in Case~(2)},
\end{array}
\right.
\end{split}
\end{align}
the minimal (generalized) Bartnik data $(\mathbb{S}^{n},g_{o},H_{o},q,\Lambda)$ have $\mathcal{A}(\mathbb{S}^{n},g_{o},H_{o},q,\Lambda)\neq\emptyset$ and
\begin{align}\label{eq:massestimate}
\mathcal{B}(\mathbb{S}^{n},g_{o},H_{o}=0,q,\Lambda)&\leq\mathcal{M}(\mathbb{S}^{n},g_{o},H_{o}=0,q,\Lambda).
\end{align}
Moreover, the \RN\ of mass $\mathcal{M}(\mathbb{S}^{n},g_{o},H_{o}=0,q,\Lambda)$, charge $q$, and cosmological constant $\Lambda$ is sub-extremal.
\end{theorem}
\begin{proof}
Under the given assumptions, Theorem~\ref{thm:main} tells us that for any constant $m>m_{o}\definedas\tfrac{r_{o}^{n-1}}{2}\,p_{0,q,\Lambda}(r_{o})$ there exists a (generalized) admissible extension $(M,\gamma,E,\Lambda)$ of $(\mathbb{S}^{n},g_{o},H_{o},q,\Lambda)$ of mass $m$ and total charge $q$. This proves the claim that $\mathcal{A}(\mathbb{S}^{n},g_{o},H_{o},q,\Lambda)\neq\emptyset$ and shows that $\mathcal{B}(\mathbb{S}^{n},g_{o},H_{o}=0,q,\Lambda)\leq m_{o}$. By construction, we have $r_{o}=r(\partial M)$ so that $m_{o}=\mathcal{M}(\mathbb{S}^{n},g_{o},H_{o}=0,q,\Lambda)$ by Proposition~\ref{prop:Hawkingre} which shows the desired inequality~\eqref{eq:massestimate}.

Moreover, Lemma~\ref{lem:subapriori} tells us that $h_{q,\Lambda}(r_{o})>0$ and thus Proposition~\ref{prop:Bartnikapriori} tells us that the \RN\ of mass $m_{o}=\mathcal{M}(\mathbb{S}^{n},g_{o},H_{o}=0,q,\Lambda)$, charge $q$, and cosmological constant $\Lambda$ is sub-extremal.
\end{proof}

\bibliographystyle{amsplain}
\bibliography{Lit_new}
\vfill
\end{document}